\documentclass[A4j,11pt]{article}
\usepackage{amsfonts,amsmath,amscd}
\usepackage{graphics}
\usepackage{amssymb}
\begin{document}

\title{{\bf The mean curvature flow for\\
invariant hypersurfaces in a Hilbert space\\
with an almost free group action}}
\author{{\bf Naoyuki Koike}}
\date{}
\maketitle

\begin{abstract}
In this paper, we study the regularized mean curvature flow starting from 
invariant hypersurfaces in a Hilbert space equipped with an isometric 
almost free Hilbert Lie group action whose orbits are minimal regularizable submanifolds, 
where ``almost free'' means that the stabilizers of the group action are finite.  
First we obtain the evolution equations for some geometric quantities along the regularized mean curvature flow.  
Next, by using the evolution equations, we prove a horizontally strongly 
convexity preservability theorem for the regularized mean curvature flow.  
From this theorem, we derive the strongly convexity preservability theorem for 
the mean curvature flow starting from compact Riemannian suborbifolds in the orbit space 
(which is a Riemannian orbifold) of the Hilbert Lie group action.  
\end{abstract}

\vspace{0.5truecm}





\section{Introduction}
R. S. Hamilton ([Ha]) proved the existenceness and the uniqueness 
(in short time) of solutions satisfying any initial condition of 
a weakly parabolic 
equation for sections of a finite dimensional vector bundle.  
The Ricci flow equation for 
Riemannian metrics on a fixed compact manifold $M$ is a weakly parabolic 
equation, where we note that the Riemannian metrics are sections of the 
$(0,2)$-tensor bundle $T^{(0,2)}M$ of $M$.  
Let $f_t$ ($0\leq t<T$) be a $C^{\infty}$-family of immersions of $M$ into 
the $m$-dimensional Euclidean space ${\Bbb R}^m$.  Define a map 
$F:M\times [0,T)\to{\Bbb R}^m$ by 
$F(x,t):=f_t(x)$ ($(x,t)\in M\times[0,T)$).  
The mean curvature flow equation is described as 
$$\frac{\partial F}{\partial t}=\triangle_tf_t,$$
where $\triangle_t$ is the Laplacian operator of the metric $g_t$ on $M$ 
induced from the Euclidean metric of ${\Bbb R}^m$ by $f_t$.  Here we note that 
$\triangle_tf_t$ is equal to the mean curvature vector of $f_t$.  
This evolution equation also is a weakly parabolic equation, where we note that 
the immersions $f_t$'s are regarded as sections of the trivial bundle 
$M\times{\Bbb R}^m$ over $M$ under the identification of 
$f_t$ and its graph immersion ${\rm id}_M\times f:M\to M\times{\Bbb R}^m$ 
(${\rm id}_M\,:\,$the identity map of $M$).  
Hence we can apply the Hamilton's result to 
this evolution equation and hence can show the existenceness and the 
uniqueness (in short time) of solution of this evolution equation satisfying 
any initial condition.  
In this paper, we consider the case where the ambient space is a 
(separable infinite dimensional) Hilbert space $V$.  Let $M$ be a Hilbert manifold and $f_t$ ($0\leq t<T$) be 
a $C^{\infty}$-family of immersions of $M$ into $V$.  
Assume that $f_t$ is regularizable, where "regularizability" means 
that the codimension of $f$ is finite, for each normal vector $v$ of $M$, 
the shape operator $A_v$ is a compact operator, and that the regularized trace 
${\rm Tr}_r\,A_v$ of $A_v$ and the trace ${\rm Tr}\,A_v^2$ of $A_v^2$ exist.  
Note that the notions of the regularized trace and the regularized mean curvature vector were introduced in [HLO] 
(see the next section about the definitions of these notions).  
Denote by $H_t$ the regularized mean curvature vector of $f_t$.  
Define a map $F:M\times[0,T)\to V$ as above in terms of $f_t$'s.  
We call $f_t$'s ($0\leq t<T$) the {\it regularized mean curvature flow} 
if the following evolution equation holds:
$$\frac{\partial F}{\partial t}=\triangle^r_tf_t.\leqno{(1.1)}$$
Here $\triangle^r_tf_t$ is defined as the vector field along $f_t$ satisfying 
$$\langle\triangle^r_tf_t,v\rangle:={\rm Tr}_r\langle(\nabla^tdf_t)
(\cdot,\cdot),v\rangle^{\sharp}\,\,\,\,(\forall\,v\in V),$$
where $\nabla^t$ is the Riemannian connection of the metric $g_t$ on $M$ induced from the metric 
$\langle\,\,,\,\,\rangle$ of $V$ by $f_t$, $\langle(\nabla^tdf_t)(\cdot,\cdot),v\rangle^{\sharp}$ is the 
$(1,1)$-tensor field on $M$ defined by $g_t(\langle(\nabla^tdf_t)(\cdot,\cdot),v\rangle^{\sharp}(X),Y)
=\langle(\nabla^tdf_t)(X,Y),v\rangle\,\,(X,Y\in TM)$ and ${\rm Tr}_r(\cdot)$ is the regularized trace 
of $(\cdot)$.  Note that $\triangle^r_tf_t$ is equal to $H_t$.  
In general, the existenceness and the uniqueness (in short time) of solutions 
of this evolution equation satisfying any initial condition has not been 
shown yet.  For we cannot apply the Hamilton's result to this evolution 
equation because it is regarded as the evolution equation for sections of 
the {\it infinite} dimensional vector bundle $M\times V$ over $M$.  
However we can show the existenceness and the uniqueness (in short time) of 
solutions of this evolution equation in special case.  
In this paper, we consider a isometric almost free action of 
a Hilbert Lie group $G$ on a Hilbert space $V$ whose orbits are regularized minimal, that is, 
they are regularizable submanifold and their regularized mean curvature vectors vanish, where ``almost free" 
means that the stabilizers of the action are finite.  
Let $M(\subset V)$ be a $G$-invariant submanifold in $V$.  
Assume that the image of $M$ by the orbit map of the $G$-action is compact.  Let $f$ be the inclusion 
map of $M$ into $V$.  We first show that the regularized mean curvature flow starting from $M$ exists 
uniquely in short time (see Proposition 4.1).  
In particular, we consider the case where $M$ is a hypersurface.  
The first purpose of this paper is to obtain the evolution equations for 
various geometrical quantities along the regularized mean curvature flow starting from $G$-invariant 
hypersurfaces (see Section 4).  The second purpose is to prove a maximum principal for an 
evolution equation related to a $G'$-invariant symmetric $(0,2)$-tensor fields 
$S_t$'s on a Hilbert manifold $M$ equipped with an isometric almost free Hilbert Lie 
group action $G'$ such that $M/G'$ is a finite dimensional compact Riemannian orbifold 
(see Section 5).  
The third purpose is to prove a horizontally strongly convexity preservability theorem for 
the regularized mean curvature flow starting from the above invariant hypersurface 
by using the evolution equations in Section 4 
and imitating the discussion in the proof of a maximum principal in Section 5 (see Section 6).  
From this theorem, we derive the strongly convexity preservability theorem for 
the mean curvature flow starting from compact Riemannian suborbifolds 
in the orbit space $V/G$ (which is a Riemannian orbifold) (see Section 7).  

\section{The regularized mean curvature flow} 
Let $f_t$ ($0\leq t<T$) be a one-parameter $C^{\infty}$-family of immersions 
of a manifold $M$ into a (finite dimensional) Riemannian manifold $N$, where 
$T$ is a positive constant or $T=\infty$.  
Denote by $H_t$ the mean curvature vector of $f_t$.  Define a map 
$F:M\times[0,T)\to N$ by $F(x,t)=f_t(x)$ ($(x,t)\in M\times[0,T)$).  
If, for each $t\in[0,T)$, 
$\displaystyle{\frac{\partial F}{\partial t}=H_t}$ holds, then $f_t$ ($0\leq t<T$) is called 
a {\it mean curvature flow}.  

Let $f$ be an immersion of an (infinite dimensional) Hilbert manifold $M$ 
into a Hilbert space $V$ and $A$ the shape tensor of $f$.  If ${\rm codim}\,M\,
<\,\infty$ and $A_v$ is a compact operator for each normal vector $v$ of $f$, 
then $M$ is called a {\it Fredholm submanifold}.  
In this paper, we then call $f$ a {\it Fredholm immersion}.  
Furthermore, if, for each normal vector $v$ of $M$, the regularized trace 
${\rm Tr}_r\,A_v$ and ${\rm Tr}\,A_v^2$ exist, then $M$ is called 
{\it regularizable submanifold}, where 
${\rm Tr}_r\,A_v$ is defined by 
${\rm Tr}_r\,A_v:=\sum\limits_{i=1}^{\infty}(\mu^+_i+\mu^-_i)$ 
($\mu^-_1\leq\mu^-_2\leq\cdots\leq 0\leq\cdots\leq\mu^+_2\leq\mu^+_1\,:\,$
the spectrum of $A_v$).  
Note that the notion of the regularized trace was defined in [HLO] and that 
it differs from the trace defined in terms of the zeta function in [KT].  
In this paper, we then call $f$ {\it regularizable immersion}.  
If $f$ is a regulalizable immersion, then the {\it regularized mean 
curvature vector} $H$ of $f$ is defined by 
$\langle H,v\rangle={\rm Tr}_r\,A_v\,\,(\forall\,v\in T^{\perp}M)$, where 
$\langle\,\,,\,\,\rangle$ is the inner product of $V$ and $T^{\perp}M$ is 
the normal bundle of $f$.  
If $H=0$, then $f$ is said to be {\it minimal}.  
In particular, if $f$ is of codimension one, then we call the norm 
$\vert\vert H\vert\vert$ of $H$ 
the {\it regularized mean curvature function} of $f$.  

Let $f_t$ ($0\leq t<T$) be a $C^{\infty}$-family of 
regularizable immersions of $M$ into $V$.  Denote by $H_t$ the 
regularized mean curvature vector of $f_t$.  
Define a map $F:M\times[0,T)\to V$ by 
$F(x,t):=f_t(x)$ ($(x,t)\in M\times[0,T)$).  If 
$\frac{\partial F}{\partial t}=H_t$ holds, then we call $f_t$ ($0\leq t<T$) 
the {\it regularized mean curvature flow}.  
It has not been known whether the regulalized mean curvature flow starting 
from any regularizable hypersurface exists uniquely in short time.  However 
its existence and uniqueness (in short time) is shown in a special case 
(see Proposition 4.1).  

\section{The mean curvature flow in Riemannian orbifolds}
In this section, we shall define the notion of the mean curvaure flow starting from 
a suborbifold in a Riemannian orbifold.  
First we recall the notions of a Riemannian orbifold and a suborbifold following to 
[AK,BB,GKP,Sa,Sh,Th].  
Let $M$ be a paracompact Hausdorff space and $(U,\phi,\widetilde U/\Gamma)$ a triple 
satisfying the following conditions:

\vspace{0.2truecm}

(i) $U$ is an open set of $M$, 

(ii) $\widehat U$ is an open set of ${\Bbb R}^n$ and $\Gamma$ is a finite subgroup of 
the $C^k$-diffeomorphism 

group ${\rm Diff}^k(\widehat U)$ of $\widehat U$,

(iii) $\phi$ is a homeomorphism of $U$ onto $\widehat U/\Gamma$.

\vspace{0.2truecm}

\noindent
Such a triple $(U,\phi,\widehat U/\Gamma)$ is called an $n$-{\it dimensional orbifold chart}.  
Let ${\cal O}:=\{(U_{\lambda},\phi_{\lambda},\widehat U/\Gamma_{\lambda})
\,\vert\,\lambda\in\Lambda\}$ be a family of $n$-dimensional orbifold charts of $M$ satisfying 
the following conditions:

\vspace{0.2truecm}

(O1) $\{U_{\lambda}\,\vert\,\lambda\in\Lambda\}$ is an open covering of $M$,

(O2) For $\lambda,\mu\in\Lambda$ with $U_{\lambda}\cap U_{\mu}\not=\emptyset$, 
there exists an $n$-dimensional orbifold chart 

$(W,\psi,\widehat W/\Gamma')$ 
such that $C^k$-embeddings $\rho_{\lambda}:\widehat W\hookrightarrow\widehat U_{\lambda}$ 
and $\rho_{\mu}:\widehat W\hookrightarrow\widehat U_{\mu}$ 

satisfying $\phi_{\lambda}^{-1}\circ\pi_{\Gamma_{\lambda}}\circ\rho_{\lambda}
=\psi^{-1}\circ\pi_{\Gamma'}$ and 
$\phi_{\mu}^{-1}\circ\pi_{\Gamma_{\mu}}\circ\rho_{\mu}
=\psi^{-1}\circ\pi_{\Gamma'}$, where 

$\pi_{\Gamma_{\lambda}},\,\pi_{\Gamma_{\mu}}$ and $\pi_{\Gamma'}$ are 
the orbit maps of $\Gamma_{\lambda},\,\Gamma_{\mu}$ and $\Gamma'$, respectively.  

\vspace{0.2truecm}

\noindent
Such a family ${\cal O}$ is called an $n$-{\it dimensional} $C^k$-{\it orbifold atlas} of $M$ and 
the pair $(M,{\cal O})$ is called an $n$-{\it dimensional} $C^k$-{\it orbifold}.  
Let $(U_{\lambda},\phi_{\lambda},\widehat U_{\lambda}/\Gamma_{\lambda})$ be 
an $n$-dimensional orbifold chart around $x\in M$.  Then the group 
$(\Gamma_{\lambda})_{\widehat x}:=\{b\in\Gamma_{\lambda}\,\vert\,b(\widehat x)=\widehat x\}$ is 
unique for $x$ up to the conjugation, where $\widehat x$ is a point of $\widehat U_{\lambda}$ 
with $(\phi_{\lambda}^{-1}\circ\pi_{\Gamma_{\lambda}})(\widehat x)=x$.  
Denote by $(\Gamma_{\lambda})_x$ the conjugate class of this group $(\Gamma_{\lambda})_{\widehat x}$, 
This conjugate class is called the {\it local group at} $x$.  
If the local group at $x$ is not trivial, then $x$ is called a {\it singular point} of 
$(M,{\cal O})$.  
Denote by ${\rm Sing}(M,{\cal O})$ (or ${\rm Sing}(M)$) 
the set of all singular points of $(M,{\cal O})$.  
This set ${\rm Sing}(M,{\cal O})$ is called the {\it singular set} of $(M,{\cal O})$.  
Let $x\in M$ and $(U_{\lambda},\phi_{\lambda},\widehat U_{\lambda}/\Gamma_{\lambda})$ an orbifold chart 
around $x$.  Take $\widehat x_{\lambda}\in \widehat U_{\lambda}$ with $\pi_{\Gamma_{\lambda}}
(\widehat x_{\lambda})=x$.  
The group $(\Gamma_{\lambda})_{\widehat x_{\lambda}}$ acts on 
$T_{\widehat x_{\lambda}}\widehat U_{\lambda}$ naturally.  
Denote by ${\cal O}_{M,x}$ the subfamily of ${\cal O}_M$ 
consisting of all orbifold charts around $x$.  
Give $\displaystyle{{\cal T}_x
:=\mathop{\oplus}_{(U_{\lambda},\phi_{\lambda},\widehat U_{\lambda}/\Gamma_{\lambda})
\in{\cal O}_{M,x}}
T_{\widehat x_{\lambda}}\widehat U_{\lambda}/(\Gamma_{\lambda})_{\widehat x_{\lambda}}}$
an equivalence relation $\sim$ as follows.  
Let $(U_{\lambda_1},\phi_{\lambda_1},\widehat U_{\lambda_1}/\Gamma_{\lambda_1})$ and 
$(U_{\lambda_2},\phi_{\lambda_2},\widehat U_{\lambda_2}/\Gamma_{\lambda_2})$ be 
members of ${\cal O}_{M,x}$.  
Let $\eta$ be the diffeomorphism of a sufficiently small neighborhood of $\widehat x_{\lambda_1}$ in 
$\widehat U_{\lambda_1}$ into $\widehat U_{\lambda_2}$ satisfying 
$\phi_{\lambda_2}^{-1}\circ\pi_{\Gamma_{\lambda_2}}\circ\eta=\phi_{\lambda_1}^{-1}\circ
\pi_{\Gamma_{\lambda_1}}$.  Define an equivalence relatiton $\sim$ in ${\cal T}_x$
as the relation generated by 
$$\begin{array}{c}
\displaystyle{[v_1]\sim[v_2]\,\,\,\mathop{\Longleftrightarrow}_{\rm def}\,\,\,
[v_2]=[\eta_{\ast}(v_1)]}\\
\displaystyle{([v_1]\in 
T_{\widehat x_{\lambda_1}}\widehat U_{\lambda_1}/(\Gamma_{\lambda_1})_{\widehat x_{\lambda_1}},
\,\,[v_2]\in 
T_{\widehat x_{\lambda_2}}\widehat U_{\lambda_2}/(\Gamma_{\lambda_2})_{\widehat x_{\lambda_2}}),}
\end{array}$$
where $[v_i]$ ($i=1,2$) is the $(\Gamma_{\lambda_i})_{\widehat x_{\lambda_i}}$-orbits through 
$v_i\in T_{\widehat x_{\lambda_i}}\widehat U_{\lambda_i}$.  
We call the quotient space ${\cal T}_x/\sim$ the {\it orbitangent space of} $M$ {\it at} $x$ and 
denote it by $T_xM$.  If $(M,{\cal O})$ is of class $C^k$ ($r\geq1$), then 
$\displaystyle{TM:=\mathop{\oplus}_{x\in M}T_xM}$ is a $C^{r-1}$-orbifold in a natural manner.  
We call $TM$ the {\it orbitangent bundle of} $M$.  
The group $(\Gamma_{\lambda})_{\widehat x_{\lambda}}$ acts on the $(r,s)$-tensor space 
$T_{\widehat x_{\lambda}}^{(r,s)}\widehat U_{\lambda}$ of 
$T_{\widehat x_{\lambda}}\widehat U_{\lambda}$ naturally.  
Give $\displaystyle{{\cal T}^{(r,s)}_x
:=\mathop{\oplus}_{(U_{\lambda},\phi_{\lambda},\widehat U_{\lambda}/\Gamma_{\lambda})
\in{\cal O}_{M,x}}
T_{\widehat x_{\lambda}}^{(r,s)}\widehat U_{\lambda}/(\Gamma_{\lambda})_{\widehat x_{\lambda}}}$
an equivalence relation $\sim$ as above.  
We call the quotient space ${\cal T}^{(r,s)}_x/\sim$ the $(r,s)$-{\it orbitensor space of} $M$ {\it at} 
$x$ and denote it by $T_x^{(r,s)}M$.  If $(M,{\cal O})$ is of class $C^k$ ($r\geq1$), then 
$\displaystyle{T^{(r,s)}M:=\mathop{\oplus}_{x\in M}T^{(r,s)}_xM}$ is a $C^{r-1}$-orbifold 
in a natural manner.  
We call $T^{(r,s)}M$ the $(r,s)$-{\it orbitensor bundle of} $M$.  

Let $(M,{\cal O}_M)$ and $(N,{\cal O}_N)$ be orbifolds, and 
$f$ a map from $M$ to $N$.  
If, for each $x\in M$ and each pair of an orbifold chart 
$(U_{\lambda},\phi_{\lambda},\widehat U_{\lambda}/\Gamma_{\lambda})$ of $(M,{\cal O}_M)$ 
around $x$ and an orbifold chart 
$(V_{\mu},\psi_{\mu},\widehat V_{\mu}/\Gamma'_{\mu})$ of $(N,{\cal O}_N)$ 
around $f(x)$ ($f(U_{\lambda})\subset V_{\mu}$), 
there exists a $C^k$-map $\widehat f_{\lambda,\mu}:\widehat U_{\lambda}\to\widehat V_{\mu}$ with 
$f\circ\phi_{\lambda}^{-1}\circ\pi_{\Gamma_{\lambda}}
=\psi_{\mu}^{-1}\circ\pi_{\Gamma'_{\mu}}\circ\widehat f_{\lambda,\mu}$, then 
$f$ is called a $C^k$-{\it orbimap} (or simply a $C^k$-{\it map}).  Also 
$\widehat f_{\lambda,\mu}$ is called a {\it local lift} of $f$ with respect to 
$(U_{\lambda},\phi_{\lambda},\widehat U_{\lambda}/\Gamma_{\lambda})$ and 
$(V_{\mu},\psi_{\mu},\widehat V_{\mu}/\Gamma'_{\mu})$.  
Furthermore, if each local lift $\widehat f_{\lambda,\mu}$ is an immersion, then 
$f$ is called a $C^k$-{\it orbiimmersion} (or simply a $C^k$-{\it immersion}) and 
$(M,{\cal O}_M)$ is called a $C^k$-{\it (immersed) suborbifold} in $(N,{\cal O}_N,g)$.  
Similarly, if each local lift $\widehat f_{\lambda,\mu}$ is a submersion, then 
$f$ is called a $C^k$-{\it orbisubmersion}.  

In the sequel, we assume that $r=\infty$.  
Denote by ${\rm pr}_{TM}$ and ${\rm pr}_{T^{(r,s)}M}$ the natural projections 
of $TM$ and $T^{(r,s)}M$ onto $M$, respectively.  
These are $C^{\infty}$-orbimaps.  
We call a $C^k$-orbimap $X:M\to TM$ with ${\rm pr}_{TM}\circ X={\rm id}$ a 
$C^k$-{\it orbitangent vector field} on $(M,{\cal O}_M)$ and 
a $C^k$-orbimap $S:M\to T^{(r,s)}M$ with ${\rm pr}_{T^{(r,s)}M}\circ S={\rm id}$ 
a $(r,s)$-{\it orbitensor field} of class $C^k$ on $(M,{\cal O}_M)$.  
If a $(r,s)$-orbitensor field $g$ of class $C^k$ on $(M,{\cal O}_M)$ is positive definite and 
symmetric, then we call $g$ a $C^k$-{\it Riemannian orbimetric} and $(M,{\cal O_M},g)$ 
a $C^k$-{\it Riemannian orbifold}.  

Let $f$ be a $C^{\infty}$-{\it orbiimmersion} of an $C^{\infty}$-orbifold $(M,{\cal O}_M)$ into 
$C^{\infty}$-Riemannian orbifold $(N,{\cal O}_N,g)$.  
Take an orbifold chart $(U_{\lambda},\phi_{\lambda},\widehat U_{\lambda}/\Gamma_{\lambda})$ of $M$ 
around $x$ and an orbifold chart $(V_{\mu},\psi_{\mu},\widehat V_{\mu}/\Gamma'_{\mu})$ of $N$ around 
$f(x)$ with $f(U_{\lambda})\subset V_{\mu}$.  Let $\widehat f_{\lambda,\mu}$ be the local lift of $f$ 
with respect to these orbifold charts and $\widehat g_{\mu}$ that of $g$ to $\widehat V_{\mu}$.  
Denote by $(T^{\perp}_{\widehat x_{\lambda}}\widehat U_{\lambda})_{\mu}$ 
the orthogonal complement of 
$(\widehat f_{\lambda,\mu})_{\ast}(T_{\widehat x_{\lambda}}\widehat U_{\lambda})$ in 
$(T_{\widehat{f(x)}_{\mu}}\widehat V_{\mu},(\widehat g_{\mu})_{\widehat{f(x)}_{\mu}})$.  
The group $(\Gamma'_{\mu})_{\widehat{f(x)}_{\mu}}$ acts on 
$(T^{\perp}_{\widehat x_{\lambda}}\widehat U_{\lambda})_{\mu}$ naturally.  
Give 
$${\cal T}^{\perp}_x
:=\mathop{\oplus}_{(U_{\lambda},\phi_{\lambda},\widehat U_{\lambda}/\Gamma_{\lambda})
\in{\cal O}_{M,x}}
\mathop{\oplus}_{(V_{\mu},\psi_{\mu},\widehat V_{\mu}/\Gamma_{\mu})\in{\cal O}_{N,f(x)}}
(T_{\widehat x_{\lambda}}^{\perp}\widehat U_{\lambda})_{\mu}/(\Gamma'_{\mu})_{\widehat{f(x)}_{\mu}}$$
an equivalence relation $\sim$ as follows.  
Let $(U_{\lambda_i},\phi_{\lambda_i},\widehat U_{\lambda_i}/\Gamma_{\lambda_i})$ ($i=1,2$) be 
members of ${\cal O}_{M,x}$ and 
$(V_{\mu_i},\psi_{\mu_i},\widehat V_{\mu_i}/\Gamma'_{\mu_i})$ ($i=1,2$) 
members of ${\cal O}_{N,f(x)}$ with $f(U_{\lambda_i})\subset V_{\mu_i}$.  
Let $\eta_{\mu_1,\mu_2}$ be the diffeomorphism of a sufficiently small neighborhood of 
$\widehat{f(x)}_{\mu_1}$ in $\widehat V_{\mu_1}$ into $\widehat V_{\mu_2}$ satisfying 
$\psi_{\mu_2}^{-1}\circ\pi_{\Gamma'_{\mu_2}}\circ\eta_{\mu_1,\mu_2}=\psi_{\mu_1}^{-1}\circ
\pi_{\Gamma'_{\mu_1}}$.  Define an equivalence relatiton $\sim$ in 
${\cal T}^{\perp}_x$ 
as the relation generated by 
$$\begin{array}{c}
\displaystyle{[\xi_1]\sim[\xi_2]\,\,\,\mathop{\Longleftrightarrow}_{\rm def}\,\,\,
[\xi_2]=[(\eta_{\mu_1,\mu_2})_{\ast}(\xi_1)]}\\
\displaystyle{([\xi_1]\in (T_{\widehat x_{\lambda_1}}^{\perp}\widehat U_{\lambda_1})_{\mu_1}
/(\Gamma'_{\mu_1})_{\widehat{f(x)}_{\mu_1}},\,\,[\xi_2]\in 
(T_{\widehat x_{\lambda_2}}\widehat U_{\lambda_2})_{\mu_2}/(\Gamma'_{\mu_2})_{\widehat{f(x)}_{\mu_2}}),}
\end{array}$$
where $[\xi_i]$ ($i=1,2$) is the $(\Gamma'_{\mu_i})_{\widehat{f(x)}_{\mu_i}}$-orbits through 
$\xi_i\in(T_{\widehat x_{\lambda_i}}^{\perp}\widehat U_{\lambda_i})_{\mu_i}$.  
We call the quotient space ${\cal T}^{\perp}_x/\sim$ the {\it orbinormal space of} $M$ {\it at} $x$ and 
denote it by $T^{\perp}_xM$.  If $f$ is of class $C^{\infty}$, then 
$\displaystyle{T^{\perp}M:=\mathop{\oplus}_{x\in M}T^{\perp}_xM}$ is a $C^{\infty}$-orbifold 
in a natural manner.  
We call $T^{\perp}M$ the {\it orbinormal bundle} of $M$.  
Denote by ${\rm pr}_{T^{\perp}M}$ the natural projection of $T^{\perp}M$ onto $M$.  
This is $C^{\infty}$-orbisubmersion.  
We call a $C^k$-orbimap $\xi:M\to T^{\perp}M$ with ${\rm pr}_{T^{\perp}M}\circ \xi={\rm id}$ a 
$C^k$-{\it orbinormal vector field} of $(M,{\cal O}_M)$ in $(N,{\cal O}_N,g)$.  

Take an orbifold chart $(U_{\lambda},\phi_{\lambda},\widehat U_{\lambda}/\Gamma_{\lambda})$ of $M$ 
around $x$ and an orbifold chart $(V_{\mu},\psi_{\mu},\widehat V_{\mu}/\Gamma'_{\mu})$ of $N$ around 
$f(x)$ with $f(U_{\lambda})\subset V_{\mu}$.  
Denote by $(T^{\perp}_{\widehat x_{\lambda}}\widehat U_{\lambda})_{\mu}^{(r,s)}$ 
the $(r,s)$-tensor space of $(T^{\perp}_{\widehat x_{\lambda}}\widehat U_{\lambda})_{\mu}$.  
The group $(\Gamma'_{\mu})_{\widehat{f(x)}_{\mu}}$ acts on 
$(T^{\perp}_{\widehat x_{\lambda}}\widehat U_{\lambda})_{\mu}^{(r,s)}$ naturally.  
Give 
$$({\cal T}^{\perp}_x)^{(r,s)}
:=\mathop{\oplus}_{(U_{\lambda},\phi_{\lambda},\widehat U_{\lambda}/\Gamma_{\lambda})
\in{\cal O}_{M,x}}
\mathop{\oplus}_{(V_{\mu},\psi_{\mu},\widehat V_{\mu}/\Gamma_{\mu})\in{\cal O}_{N,f(x)}}
(T_{\widehat x_{\lambda}}^{\perp}\widehat U_{\lambda})_{\mu}^{(r,s)}
/(\Gamma'_{\mu})_{\widehat{f(x)}_{\mu}}$$
an equivalence relation $\sim$ as follows.  
Let $(U_{\lambda_i},\phi_{\lambda_i},\widehat U_{\lambda_i}/\Gamma_{\lambda_i}),\,
(V_{\mu_i},\psi_{\mu_i},\widehat V_{\mu_i}/\Gamma'_{\mu_i})$ ($i=1,2$) 
and $\eta_{\mu_1,\mu_2}$ be as above.  
Define an equivalence relatiton $\sim$ in 
$({\cal T}^{\perp}_x)^{(r,s)}$ 
as the relation generated by 
$$\begin{array}{c}
\displaystyle{[S_1]\sim[S_2]\,\,\,\mathop{\Longleftrightarrow}_{\rm def}\,\,\,
[S_1]=[(\eta_{\mu_1,\mu_2})^{\ast}(S_1)]}\\
\displaystyle{([S_1]\in(T_{\widehat x_{\lambda_1}}^{\perp}\widehat U_{\lambda_1})_{\mu_1}^{(r,s)}
/(\Gamma'_{\mu_1})_{\widehat{f(x)}_{\mu_1}},\,\,[S_2]\in 
(T_{\widehat x_{\lambda_2}}\widehat U_{\lambda_2})_{\mu_2}^{(r,s)}/
(\Gamma'_{\mu_2})_{\widehat{f(x)}_{\mu_2}}),}
\end{array}$$
where $[S_i]$ ($i=1,2$) is the $(\Gamma'_{\mu_i})_{\widehat{f(x)}_{\mu_i}}$-orbits through 
$S_i\in(T_{\widehat x_{\lambda_i}}^{\perp}\widehat U_{\lambda_i})_{\mu_i}^{(r,s)}$.  
We denote the quotient space $({\cal T}^{\perp}_x)^{(r,s)}/\sim$ by 
$(T^{\perp}_xM)^{(r,s)}$.  If $f$ is of class $C^{\infty}$, then 
$\displaystyle{(T^{\perp}M)^{(r,s)}:=\mathop{\oplus}_{x\in M}(T^{\perp}_xM)^{(r,s)}}$ is a 
$C^{\infty}$-orbifold in a natural manner.  
We call $(T^{\perp}M)^{(r,s)}$ the $(r,s)$-{\it orbitensor bundle} of $T^{\perp}M$.  
Denote by ${\rm pr}_{(T^{\perp}M)^{(r,s)}}$ the natural projection of 
$(T^{\perp}M)^{(r,s)}$ onto $M$.  This is $C^{\infty}$-orbisubmersion.  

Next we shall define the tensor product 
$T^{(r,s)}M\otimes(T^{\perp}M)^{(s',t')}$ of $T^{(r,s)}M$ and $(T^{\perp}M)^{(s',t')}$.  
Take an orbifold chart $(U_{\lambda},\phi_{\lambda},\widehat U_{\lambda}/\Gamma_{\lambda})$ of $M$ 
around $x$ and an orbifold chart $(V_{\mu},\psi_{\mu},\widehat V_{\mu}/\Gamma'_{\mu})$ of $N$ around 
$f(x)$ with $f(U_{\lambda})\subset V_{\mu}$.  
The group $(\Gamma_{\lambda})_{\widehat x_{\lambda}}\times(\Gamma'_{\mu})_{\widehat{f(x)}_{\mu}}$ 
acts on $(T^{(r,s)}_{\widehat x_{\lambda}}\widehat U_{\lambda})\otimes
(T^{\perp}_{\widehat x_{\lambda}}\widehat U_{\lambda})_{\mu}^{(s',t')}$ naturally.  
Give 
$$\begin{array}{l}
\hspace{0.5truecm}\displaystyle{{\cal T}^{(r,s)}_x\otimes({\cal T}^{\perp}_x)^{(s',t')}
:=\mathop{\oplus}_{(U_{\lambda},\phi_{\lambda},\widehat U_{\lambda}/\Gamma_{\lambda})\in{\cal O}_{M,x}}
\mathop{\oplus}_{(V_{\mu},\psi_{\mu},\widehat V_{\mu}/\Gamma_{\mu})\in{\cal O}_{N,f(x)}}}\\
\hspace{4.3truecm}\displaystyle{\left((T^{(r,s)}_{\widehat x_{\lambda}}\widehat U_{\lambda})\otimes
(T^{\perp}_{\widehat x_{\lambda}}\widehat U_{\lambda})_{\mu}^{(s',t')})\right)
/((\Gamma_{\lambda})_{\widehat x_{\lambda}}\times(\Gamma'_{\mu})_{\widehat{f(x)}_{\mu}})}
\end{array}$$
an equivalence relation $\sim$ as follows.  
Let $(U_{\lambda_i},\phi_{\lambda_i},\widehat U_{\lambda_i}/\Gamma_{\lambda_i}),\,
(V_{\mu_i},\psi_{\mu_i},\widehat V_{\mu_i}/\Gamma'_{\mu_i})$ ($i=1,2$) and 
$\eta_{\mu_1,\mu_2}$ be as above.  Also let $\eta_{\lambda_1,\lambda_2}$ be a diffeomorphism 
defined in similar to $\eta_{\mu_1,\mu_2}$.  
Define an equivalence relation $\sim$ in 
${\cal T}^{(r,s)}_x\otimes({\cal T}^{\perp}_x)^{(s',t')}$ 
as the relation generated by 
$$\begin{array}{c}
\displaystyle{[S_1]\sim[S_2]\,\,\,\mathop{\Longleftrightarrow}_{\rm def}\,\,\,
[S_1]=[(\eta_{\lambda_1,\lambda_2}^{\ast}\otimes\eta_{\mu_1,\mu_2}^{\ast})S_2]}\\
\displaystyle{([S_1]\in
(T^{(r,s)}_{\widehat x_{\lambda_1}}\widehat U_{\lambda_1})\otimes
(T^{\perp}_{\widehat x_{\lambda_1}}\widehat U_{\lambda_1})_{\mu_1}^{(s',t')}),\,\,\,\,
[S_2]\in(T^{(r,s)}_{\widehat x_{\lambda_2}}\widehat U_{\lambda_2})\otimes
(T^{\perp}_{\widehat x_{\lambda_2}}\widehat U_{\lambda_2})_{\mu_2}^{(s',t')})),}
\end{array}$$
where $[S_i]$ ($i=1,2$) is the 
$((\Gamma_{\lambda_i})_{\widehat x_{\lambda_i}}\times(\Gamma'_{\mu_i})_{\widehat{f(x)}_{\mu_i}})$-orbits 
through $S_i\in
((T^{(r,s)}_{\widehat x_{\lambda_i}}\widehat U_{\lambda_i})\otimes
(T^{\perp}_{\widehat x_{\lambda_i}}\widehat U_{\lambda_i})_{\mu_i}^{(s',t')}))$.  
We denote the quotient space $({\cal T}^{(r,s)}_x\otimes({\cal T}^{\perp}_x)^{(s',t')})/\sim$ by 
$T^{(r,s)}_xM\otimes(T^{\perp}_xM)^{(s',t')}$.  
Set $\displaystyle{T^{(r,s)}M\otimes(T^{\perp}M)^{(s',t')}
:=\mathop{\oplus}_{x\in M}\left(T^{(r,s)}_xM\otimes(T^{\perp}_xM)^{(s',t')}\right)}$.  
If $f$ is of class $C^{\infty}$, then 
$T^{(r,s)}M\otimes(T^{\perp}M)^{(s',t')}$ is a $C^{\infty}$-orbifold in a natural manner.  
We call $T^{(r,s)}M\otimes(T^{\perp}M)^{(s',t')}$ the {\it orbitensor product bundle} of 
$T^{(r,s)}M$ and $(T^{\perp}M)^{(s',t')}$.  
Denote by ${\rm pr}_{T^{(r,s)}M\otimes(T^{\perp}M)^{(s',t')}}$ the natural projection of 
$T^{(r,s)}M\otimes(T^{\perp}M)^{(s',t')}$ onto $M$.  This is a $C^{\infty}$-orbisubmersion.  
We call a $C^k$-orbimap $S:M\to T^{(r,s)}M\otimes(T^{\perp}M)^{(s',t')}$ with 
${\rm pr}_{T^{(r,s)}M\otimes(T^{\perp}M)^{(s',t')}}\circ S={\rm id}$ a 
$C^k$-{\it section of} $T^{(r,s)}M\otimes(T^{\perp}M)^{(s',t')}$.  
Let $g,\,h,\,A,\,H$ and $\xi$ be the induced metric, the second fundamental form, the shape tensor, 
the mean curvature and a unit normal vector field of the immersion 
$f\vert_{M\setminus{\rm Sing}(M)}:M\setminus{\rm Sing}(M)\hookrightarrow N\setminus{\rm Sing}(N)$, 
respectively.  It is easy to show that $g,\,h,\,A$ and $H$ extend 
a $(0,2)$-orbitensor field of class $C^{\infty}$ on $(M,{\cal O}_M)$, 
a $C^k$-section of $T^{(0,2)}M\otimes T^{\perp}M$, a $C^k$-section of 
$T^{(1,1)}M\otimes(T^{\perp}M)^{(0,1)}$ and a $C^{\infty}$-orbinormal vector field on $(M,{\cal O}_M)$.  
We denote these extensions by the same symbols.  
We call these extensions $g,\,h,\,A$ and $H$ the {\it induced orbimetric}, the 
{\it second fundamental orbiform}, the {\it shape orbitensor} and the 
{\it mean curvature orbifunction} of $f$.  
Here we note that $\xi$ does not necessarily extend a $C^{\infty}$-orbinormal vector field on 
$(M,{\cal O})$ (see Figure 2).  

Now we shall define the notion of the mean curvature flow starting from 
a $C^{\infty}$-suborbifold in a $C^{\infty}$-Riemannian orbifold.  
Let $f_t$ ($0\leq t<T$) be a $C^{\infty}$-family of 
$C^{\infty}$-orbiimmersions of a $C^{\infty}$-orbifold $(M,{\cal O}_M)$ into 
a $C^{\infty}$-Riemannian orbifold $(N,{\cal O}_N,g)$.  
Assume that, for each $(x_0,t_0)\in M\times[0,T)$ and each pair of an orbifold chart 
$(U_{\lambda},\phi_{\lambda},\widehat U_{\lambda}/\Gamma_{\lambda})$ of $(M,{\cal O}_M)$ around 
$x_0$ and an orbifold chart $(V_{\mu},\phi_{\mu},\widehat V_{\mu}/\Gamma'_{\mu})$ of 
$(N,{\cal O}_N)$ around $f_{t_0}(x_0)$ such that $f_t(U_{\lambda})\subset V_{\mu}$ for any 
$t\in[t_0,t_0+\varepsilon)$ ($\varepsilon:$ a sufficiently small positive number), 
there exists local lifts $(\widehat f_t)_{\lambda,\mu}:\widehat U_{\lambda}\to\widehat V_{\mu}$ of 
$f_t$ ($t\in[t_0,t_0+\varepsilon)$) such that they give the mean curvature flow in 
$(\widehat V_{\mu},\widehat g_{\mu})$, 
where $\widehat g_{\mu}$ is the local lift of $g$ to $\widehat V_{\mu}$.  
Then we call $f_t$ ($0\leq t<T$) the {\it mean curvature flow} in $(N,{\cal O}_N,g)$.  

\vspace{0.5truecm}

\noindent
{\bf Theorem 3.1.} {\sl For any $C^{\infty}$-orbiimmersion $f$ of a compact $C^{\infty}$-orbifold 
into a $C^{\infty}$-Riemannian orbifold, the mean curvature flow starting from $f$ exists uniquely 
in short time.}

\vspace{0.5truecm}

\noindent
{\it Proof.} Let $f$ be a $C^{\infty}$-orbiimmersion of an $n$-dimensional compact 
$C^{\infty}$-orbifold $(M,{\cal O}_M)$ into an $(n+r)$-dimensional $C^{\infty}$-Riemmannian 
orbifold $(N,{\cal O}_N,g)$.  
Fix $x_0\in M$.  
Take an orbifold chart $(U_{\lambda},\phi_{\lambda},\widehat U_{\lambda}/\Gamma_{\lambda})$ of 
$(M,{\cal O}_M)$ around $x_0$ 
and an orbifold chart $(V_{\mu},\psi_{\mu},\widehat V_{\mu}/\Gamma'_{\mu})$ of $(N,{\cal O}_N)$ around 
$f(x_0)$ such that $f(U_{\lambda})\subset V_{\mu}$ and that $\widehat U_{\lambda}$ is relative compact.  
Also, let $\widehat f_{\lambda,\mu}:\widehat U_{\lambda}\hookrightarrow\widehat V_{\mu}$ be 
a local lift of $f$ and $\widehat g_{\mu}$ a local lift of $g$ (to $\widehat V$).  
Since $\widehat U_{\lambda}$ is relative compact, there exists the mean curvature flow 
$(\widehat f_{\lambda,\mu})_t:\widehat U_{\lambda}\hookrightarrow(\widehat V_{\mu},
\widehat g_{\mu})$ ($0\leq t<T)$ starting from $\widehat f_{\lambda,\mu}:\widehat U_{\lambda}
\hookrightarrow(\widehat V_{\mu},\widehat g_{\mu})$.  
Since $\widehat f_{\lambda,\mu}$ is projetable to $f\vert_{U_{\lambda}}$ and $\widehat g_{\mu}$ is 
$\Gamma'_{\mu}$-invariant, 
$(\widehat f_{\lambda,\mu})_t$ ($0\leq t<T$) also are projectable to maps of $U_{\lambda}$ into 
$V_{\mu}$.  Denote by $(f_{\lambda,\mu})_t$'s these maps of $U_{\lambda}$ into $V_{\mu}$.  
It is clear that $(f_{\lambda,\mu})_t$ ($0\leq t<T$) is the mean curvature flow starting from 
$f\vert_{U_{\lambda}}$.  
Hence, it follows from the arbitrariness of $x_0$ and the compactness of $M$ that 
the mean curvature flow starting from $f$ exists uniquely in short time.  \hspace{5.35truecm}q.e.d.

\newpage


\centerline{
\unitlength 0.1in
\begin{picture}( 38.4400, 23.9000)(  8.8000,-27.4000)
%
\special{pn 8}%
\special{pa 2580 930}%
\special{pa 4090 930}%
\special{fp}%
%
\special{pn 13}%
\special{ar 4222 890 1048 198  1.8169763 2.9073920}%
\put(33.1000,-5.2000){\makebox(0,0)[lb]{$\widehat V_{\mu}$}}%
%
\special{pn 8}%
\special{ar 2984 2492 848 434  4.8924697 4.9111905}%
\special{ar 2984 2492 848 434  4.9673527 4.9860735}%
\special{ar 2984 2492 848 434  5.0422357 5.0609565}%
\special{ar 2984 2492 848 434  5.1171187 5.1358394}%
\special{ar 2984 2492 848 434  5.1920017 5.2107224}%
\special{ar 2984 2492 848 434  5.2668847 5.2856054}%
\special{ar 2984 2492 848 434  5.3417677 5.3604884}%
\special{ar 2984 2492 848 434  5.4166507 5.4353714}%
\special{ar 2984 2492 848 434  5.4915337 5.5102544}%
\special{ar 2984 2492 848 434  5.5664167 5.5851374}%
\special{ar 2984 2492 848 434  5.6412997 5.6600204}%
\special{ar 2984 2492 848 434  5.7161827 5.7349034}%
\special{ar 2984 2492 848 434  5.7910657 5.8097864}%
\special{ar 2984 2492 848 434  5.8659487 5.8846694}%
\special{ar 2984 2492 848 434  5.9408316 5.9595524}%
\special{ar 2984 2492 848 434  6.0157146 6.0344354}%
\special{ar 2984 2492 848 434  6.0905976 6.1093184}%
\special{ar 2984 2492 848 434  6.1654806 6.1842014}%
\special{ar 2984 2492 848 434  6.2403636 6.2590844}%
%
\special{pn 8}%
\special{pa 2960 2490}%
\special{pa 3836 2490}%
\special{fp}%
%
\special{pn 8}%
\special{pa 2974 2492}%
\special{pa 3132 2060}%
\special{fp}%
%
\special{pn 13}%
\special{ar 3924 2670 1148 338  3.7016562 4.5879067}%
%
\special{pn 13}%
\special{pa 3280 2080}%
\special{pa 3272 2112}%
\special{pa 3264 2142}%
\special{pa 3252 2172}%
\special{pa 3238 2200}%
\special{pa 3222 2228}%
\special{pa 3206 2256}%
\special{pa 3186 2280}%
\special{pa 3166 2304}%
\special{pa 3146 2330}%
\special{pa 3122 2352}%
\special{pa 3100 2374}%
\special{pa 3078 2396}%
\special{pa 3052 2416}%
\special{pa 3026 2434}%
\special{pa 3000 2454}%
\special{pa 2974 2472}%
\special{pa 2958 2480}%
\special{sp}%
%
\special{pn 8}%
\special{pa 3070 2210}%
\special{pa 3108 2108}%
\special{fp}%
\special{sh 1}%
\special{pa 3108 2108}%
\special{pa 3066 2164}%
\special{pa 3090 2158}%
\special{pa 3104 2176}%
\special{pa 3108 2108}%
\special{fp}%
%
\special{pn 8}%
\special{pa 3460 2490}%
\special{pa 3564 2490}%
\special{fp}%
\special{sh 1}%
\special{pa 3564 2490}%
\special{pa 3496 2470}%
\special{pa 3510 2490}%
\special{pa 3496 2510}%
\special{pa 3564 2490}%
\special{fp}%
\put(32.9000,-27.4000){\makebox(0,0)[lt]{$V_{\mu}$}}%
%
\special{pn 8}%
\special{pa 3396 1450}%
\special{pa 3396 1826}%
\special{fp}%
\special{sh 1}%
\special{pa 3396 1826}%
\special{pa 3416 1760}%
\special{pa 3396 1774}%
\special{pa 3376 1760}%
\special{pa 3396 1826}%
\special{fp}%
\put(35.4500,-15.5400){\makebox(0,0)[lt]{$\psi_{\mu}^{-1}\circ\pi_{\Gamma'_{\mu}}$}}%
%
\special{pn 13}%
\special{pa 3290 2400}%
\special{pa 3350 2590}%
\special{fp}%
\special{sh 1}%
\special{pa 3350 2590}%
\special{pa 3350 2520}%
\special{pa 3334 2540}%
\special{pa 3312 2532}%
\special{pa 3350 2590}%
\special{fp}%
%
\special{pn 13}%
\special{pa 3630 2350}%
\special{pa 3602 2560}%
\special{fp}%
\special{sh 1}%
\special{pa 3602 2560}%
\special{pa 3630 2498}%
\special{pa 3608 2508}%
\special{pa 3590 2492}%
\special{pa 3602 2560}%
\special{fp}%
%
\special{pn 13}%
\special{pa 3130 2340}%
\special{pa 2940 2202}%
\special{fp}%
\special{sh 1}%
\special{pa 2940 2202}%
\special{pa 2982 2258}%
\special{pa 2984 2232}%
\special{pa 3006 2224}%
\special{pa 2940 2202}%
\special{fp}%
%
\special{pn 13}%
\special{pa 3270 2150}%
\special{pa 3074 2052}%
\special{fp}%
\special{sh 1}%
\special{pa 3074 2052}%
\special{pa 3126 2100}%
\special{pa 3122 2076}%
\special{pa 3144 2064}%
\special{pa 3074 2052}%
\special{fp}%
%
\special{pn 13}%
\special{pa 4280 830}%
\special{pa 4430 1010}%
\special{fp}%
\special{sh 1}%
\special{pa 4430 1010}%
\special{pa 4402 946}%
\special{pa 4396 968}%
\special{pa 4372 972}%
\special{pa 4430 1010}%
\special{fp}%
\put(44.8500,-8.7700){\makebox(0,0)[lt]{'s$\,\,:\,\,\widehat H_{\lambda,\mu}$}}%
%
\special{pn 13}%
\special{pa 4160 2260}%
\special{pa 4310 2440}%
\special{fp}%
\special{sh 1}%
\special{pa 4310 2440}%
\special{pa 4282 2376}%
\special{pa 4276 2398}%
\special{pa 4252 2402}%
\special{pa 4310 2440}%
\special{fp}%
\put(43.6500,-23.0700){\makebox(0,0)[lt]{'s$\,\,:\,\,H$}}%
%
\special{pn 13}%
\special{pa 3250 530}%
\special{pa 3234 558}%
\special{pa 3220 586}%
\special{pa 3206 616}%
\special{pa 3194 644}%
\special{pa 3184 674}%
\special{pa 3174 706}%
\special{pa 3164 736}%
\special{pa 3160 768}%
\special{pa 3156 800}%
\special{pa 3156 832}%
\special{pa 3164 862}%
\special{pa 3172 892}%
\special{pa 3190 918}%
\special{pa 3212 934}%
\special{sp}%
%
\special{pn 8}%
\special{pa 3400 2490}%
\special{pa 3504 2490}%
\special{fp}%
\special{sh 1}%
\special{pa 3504 2490}%
\special{pa 3436 2470}%
\special{pa 3450 2490}%
\special{pa 3436 2510}%
\special{pa 3504 2490}%
\special{fp}%
%
\special{pn 8}%
\special{pa 3050 2270}%
\special{pa 3088 2168}%
\special{fp}%
\special{sh 1}%
\special{pa 3088 2168}%
\special{pa 3046 2224}%
\special{pa 3070 2218}%
\special{pa 3084 2236}%
\special{pa 3088 2168}%
\special{fp}%
%
\special{pn 13}%
\special{pa 3210 930}%
\special{pa 3430 818}%
\special{fp}%
\special{sh 1}%
\special{pa 3430 818}%
\special{pa 3362 830}%
\special{pa 3382 842}%
\special{pa 3380 866}%
\special{pa 3430 818}%
\special{fp}%
%
\special{pn 13}%
\special{pa 3500 1030}%
\special{pa 3590 880}%
\special{fp}%
\special{sh 1}%
\special{pa 3590 880}%
\special{pa 3540 928}%
\special{pa 3564 926}%
\special{pa 3574 948}%
\special{pa 3590 880}%
\special{fp}%
\special{pa 3590 880}%
\special{pa 3590 880}%
\special{fp}%
\put(12.8000,-7.8000){\makebox(0,0)[lb]{$\widehat U_{\lambda}$}}%
\put(20.8000,-7.7000){\makebox(0,0)[lt]{$\displaystyle{\mathop{\hookrightarrow}^{\widehat f_{\lambda,\mu}}}$}}%
%
\special{pn 13}%
\special{ar 1366 820 540 200  0.4172545 2.6852442}%
%
\special{pn 8}%
\special{pa 1370 1500}%
\special{pa 1370 1876}%
\special{fp}%
\special{sh 1}%
\special{pa 1370 1876}%
\special{pa 1390 1810}%
\special{pa 1370 1824}%
\special{pa 1350 1810}%
\special{pa 1370 1876}%
\special{fp}%
\put(14.6500,-15.8400){\makebox(0,0)[lt]{$\phi_{\lambda}^{-1}\circ\pi_{\Gamma_{\lambda}}$}}%
\put(11.8000,-26.3000){\makebox(0,0)[lt]{$U_{\lambda}$}}%
\put(22.4000,-22.0000){\makebox(0,0)[lt]{$\displaystyle{\mathop{\hookrightarrow}^{f}}$}}%
%
\special{pn 13}%
\special{pa 3180 650}%
\special{pa 3384 690}%
\special{fp}%
\special{sh 1}%
\special{pa 3384 690}%
\special{pa 3322 658}%
\special{pa 3332 680}%
\special{pa 3314 696}%
\special{pa 3384 690}%
\special{fp}%
%
\special{pn 13}%
\special{pa 3800 1060}%
\special{pa 3850 890}%
\special{fp}%
\special{sh 1}%
\special{pa 3850 890}%
\special{pa 3812 948}%
\special{pa 3836 942}%
\special{pa 3850 960}%
\special{pa 3850 890}%
\special{fp}%
%
\special{pn 13}%
\special{pa 3150 810}%
\special{pa 3408 762}%
\special{fp}%
\special{sh 1}%
\special{pa 3408 762}%
\special{pa 3340 756}%
\special{pa 3356 772}%
\special{pa 3346 794}%
\special{pa 3408 762}%
\special{fp}%
%
\special{pn 8}%
\special{pa 3020 1180}%
\special{pa 3430 610}%
\special{fp}%
%
\special{pn 13}%
\special{pa 2980 2460}%
\special{pa 3270 2330}%
\special{fp}%
\special{sh 1}%
\special{pa 3270 2330}%
\special{pa 3202 2340}%
\special{pa 3222 2352}%
\special{pa 3218 2376}%
\special{pa 3270 2330}%
\special{fp}%
\put(47.2400,-10.8700){\makebox(0,0)[lt]{(a local lift of $H$)}}%
%
\special{pn 13}%
\special{ar 1410 2230 540 200  0.4172545 2.6852442}%
\end{picture}%
\hspace{3truecm}}

\vspace{0.5truecm}

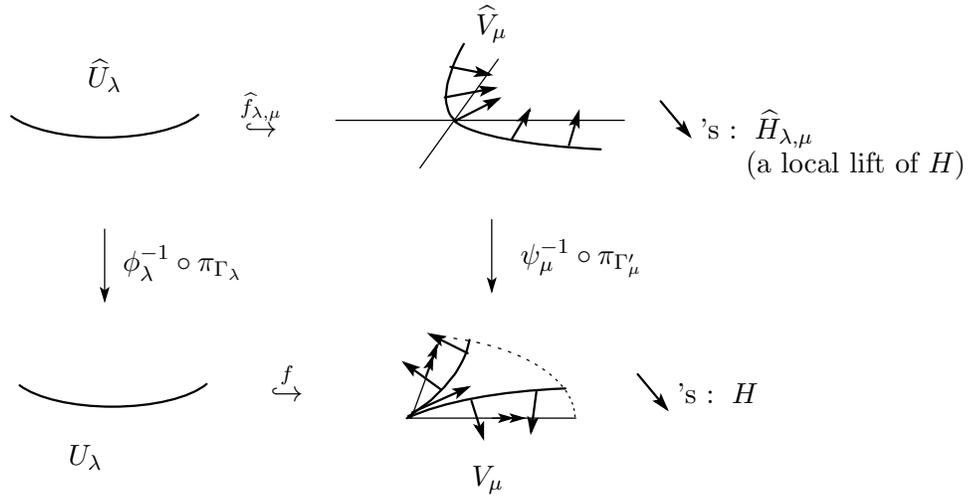
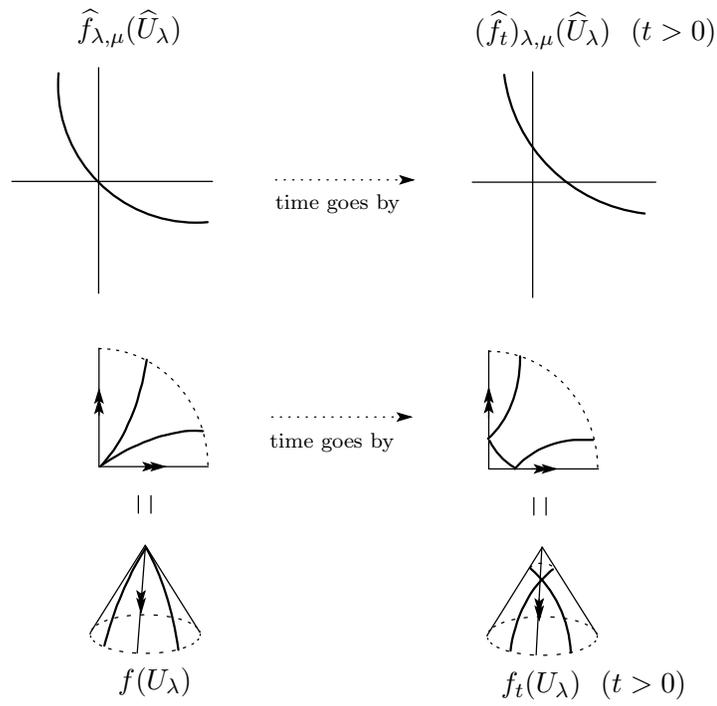
\centerline{{\bf Figure 1.}}

\vspace{0.5truecm}

\centerline{
\unitlength 0.1in
\begin{picture}( 33.7000, 34.2000)(  6.9000,-37.7000)
%
\special{pn 8}%
\special{pa 1146 2720}%
\special{pa 1146 2100}%
\special{fp}%
%
\special{pn 8}%
\special{pa 1146 2720}%
\special{pa 1718 2720}%
\special{fp}%
%
\special{pn 8}%
\special{ar 1146 2720 570 618  4.7123890 4.7326080}%
\special{ar 1146 2720 570 618  4.7932651 4.8134842}%
\special{ar 1146 2720 570 618  4.8741413 4.8943603}%
\special{ar 1146 2720 570 618  4.9550175 4.9752365}%
\special{ar 1146 2720 570 618  5.0358936 5.0561127}%
\special{ar 1146 2720 570 618  5.1167698 5.1369888}%
\special{ar 1146 2720 570 618  5.1976459 5.2178650}%
\special{ar 1146 2720 570 618  5.2785221 5.2987411}%
\special{ar 1146 2720 570 618  5.3593982 5.3796173}%
\special{ar 1146 2720 570 618  5.4402744 5.4604934}%
\special{ar 1146 2720 570 618  5.5211506 5.5413696}%
\special{ar 1146 2720 570 618  5.6020267 5.6222458}%
\special{ar 1146 2720 570 618  5.6829029 5.7031219}%
\special{ar 1146 2720 570 618  5.7637790 5.7839981}%
\special{ar 1146 2720 570 618  5.8446552 5.8648742}%
\special{ar 1146 2720 570 618  5.9255314 5.9457504}%
\special{ar 1146 2720 570 618  6.0064075 6.0266266}%
\special{ar 1146 2720 570 618  6.0872837 6.1075027}%
\special{ar 1146 2720 570 618  6.1681598 6.1883789}%
\special{ar 1146 2720 570 618  6.2490360 6.2692550}%
%
\special{pn 8}%
\special{pa 3186 2732}%
\special{pa 3186 2112}%
\special{fp}%
%
\special{pn 8}%
\special{pa 3186 2732}%
\special{pa 3760 2732}%
\special{fp}%
%
\special{pn 8}%
\special{ar 3186 2732 570 616  4.7123890 4.7326251}%
\special{ar 3186 2732 570 616  4.7933333 4.8135694}%
\special{ar 3186 2732 570 616  4.8742777 4.8945138}%
\special{ar 3186 2732 570 616  4.9552220 4.9754581}%
\special{ar 3186 2732 570 616  5.0361664 5.0564025}%
\special{ar 3186 2732 570 616  5.1171107 5.1373468}%
\special{ar 3186 2732 570 616  5.1980551 5.2182912}%
\special{ar 3186 2732 570 616  5.2789994 5.2992355}%
\special{ar 3186 2732 570 616  5.3599438 5.3801799}%
\special{ar 3186 2732 570 616  5.4408881 5.4611242}%
\special{ar 3186 2732 570 616  5.5218325 5.5420686}%
\special{ar 3186 2732 570 616  5.6027768 5.6230129}%
\special{ar 3186 2732 570 616  5.6837212 5.7039573}%
\special{ar 3186 2732 570 616  5.7646655 5.7849016}%
\special{ar 3186 2732 570 616  5.8456099 5.8658460}%
\special{ar 3186 2732 570 616  5.9265542 5.9467903}%
\special{ar 3186 2732 570 616  6.0074986 6.0277347}%
\special{ar 3186 2732 570 616  6.0884429 6.1086790}%
\special{ar 3186 2732 570 616  6.1693873 6.1896234}%
\special{ar 3186 2732 570 616  6.2503316 6.2705677}%
%
\special{pn 13}%
\special{pa 1396 2160}%
\special{pa 1390 2192}%
\special{pa 1384 2222}%
\special{pa 1378 2254}%
\special{pa 1372 2286}%
\special{pa 1362 2316}%
\special{pa 1354 2346}%
\special{pa 1346 2378}%
\special{pa 1336 2408}%
\special{pa 1324 2438}%
\special{pa 1312 2468}%
\special{pa 1298 2496}%
\special{pa 1286 2526}%
\special{pa 1272 2554}%
\special{pa 1256 2582}%
\special{pa 1238 2610}%
\special{pa 1222 2636}%
\special{pa 1202 2660}%
\special{pa 1182 2686}%
\special{pa 1162 2710}%
\special{pa 1152 2720}%
\special{sp}%
%
\special{pn 13}%
\special{pa 1690 2534}%
\special{pa 1658 2532}%
\special{pa 1626 2532}%
\special{pa 1594 2534}%
\special{pa 1564 2540}%
\special{pa 1532 2546}%
\special{pa 1500 2554}%
\special{pa 1470 2560}%
\special{pa 1440 2570}%
\special{pa 1410 2582}%
\special{pa 1380 2592}%
\special{pa 1352 2606}%
\special{pa 1322 2618}%
\special{pa 1292 2632}%
\special{pa 1264 2648}%
\special{pa 1238 2664}%
\special{pa 1210 2680}%
\special{pa 1182 2698}%
\special{pa 1164 2710}%
\special{sp}%
%
\special{pn 13}%
\special{ar 3472 2420 320 346  2.0734485 2.6475073}%
%
\special{pn 13}%
\special{ar 3000 2150 350 496  6.2676717 6.2831853}%
\special{ar 3000 2150 350 496  0.0000000 1.0230866}%
%
\special{pn 13}%
\special{pa 3322 2732}%
\special{pa 3344 2708}%
\special{pa 3366 2686}%
\special{pa 3392 2668}%
\special{pa 3418 2650}%
\special{pa 3446 2632}%
\special{pa 3474 2620}%
\special{pa 3504 2608}%
\special{pa 3534 2598}%
\special{pa 3566 2590}%
\special{pa 3598 2584}%
\special{pa 3628 2580}%
\special{pa 3660 2580}%
\special{pa 3692 2580}%
\special{pa 3724 2580}%
\special{pa 3732 2580}%
\special{sp}%
%
\special{pn 8}%
\special{ar 1388 3598 288 100  6.2831853 6.3453615}%
\special{ar 1388 3598 288 100  6.5318900 6.5940661}%
\special{ar 1388 3598 288 100  6.7805946 6.8427708}%
\special{ar 1388 3598 288 100  7.0292993 7.0914755}%
\special{ar 1388 3598 288 100  7.2780040 7.3401801}%
\special{ar 1388 3598 288 100  7.5267086 7.5888848}%
\special{ar 1388 3598 288 100  7.7754133 7.8375895}%
\special{ar 1388 3598 288 100  8.0241179 8.0862941}%
\special{ar 1388 3598 288 100  8.2728226 8.3349988}%
\special{ar 1388 3598 288 100  8.5215273 8.5837034}%
\special{ar 1388 3598 288 100  8.7702319 8.8324081}%
\special{ar 1388 3598 288 100  9.0189366 9.0811128}%
\special{ar 1388 3598 288 100  9.2676413 9.3298174}%
%
\special{pn 8}%
\special{ar 1388 3598 288 102  3.1415927 3.2034483}%
\special{ar 1388 3598 288 102  3.3890153 3.4508710}%
\special{ar 1388 3598 288 102  3.6364380 3.6982937}%
\special{ar 1388 3598 288 102  3.8838607 3.9457164}%
\special{ar 1388 3598 288 102  4.1312834 4.1931390}%
\special{ar 1388 3598 288 102  4.3787061 4.4405617}%
\special{ar 1388 3598 288 102  4.6261287 4.6879844}%
\special{ar 1388 3598 288 102  4.8735514 4.9354071}%
\special{ar 1388 3598 288 102  5.1209741 5.1828298}%
\special{ar 1388 3598 288 102  5.3683968 5.4302524}%
\special{ar 1388 3598 288 102  5.6158195 5.6776751}%
\special{ar 1388 3598 288 102  5.8632421 5.9250978}%
\special{ar 1388 3598 288 102  6.1106648 6.1725205}%
%
\special{pn 8}%
\special{pa 1388 3132}%
\special{pa 1102 3588}%
\special{fp}%
%
\special{pn 8}%
\special{pa 1388 3132}%
\special{pa 1674 3588}%
\special{fp}%
%
\special{pn 8}%
\special{pa 1388 3148}%
\special{pa 1346 3696}%
\special{fp}%
%
\special{pn 8}%
\special{ar 3466 3608 288 100  6.2831853 6.3450410}%
\special{ar 3466 3608 288 100  6.5306080 6.5924637}%
\special{ar 3466 3608 288 100  6.7780307 6.8398863}%
\special{ar 3466 3608 288 100  7.0254533 7.0873090}%
\special{ar 3466 3608 288 100  7.2728760 7.3347317}%
\special{ar 3466 3608 288 100  7.5202987 7.5821544}%
\special{ar 3466 3608 288 100  7.7677214 7.8295771}%
\special{ar 3466 3608 288 100  8.0151441 8.0769997}%
\special{ar 3466 3608 288 100  8.2625668 8.3244224}%
\special{ar 3466 3608 288 100  8.5099894 8.5718451}%
\special{ar 3466 3608 288 100  8.7574121 8.8192678}%
\special{ar 3466 3608 288 100  9.0048348 9.0666905}%
\special{ar 3466 3608 288 100  9.2522575 9.3141131}%
%
\special{pn 8}%
\special{ar 3466 3608 288 102  3.1415927 3.2032893}%
\special{ar 3466 3608 288 102  3.3883793 3.4500759}%
\special{ar 3466 3608 288 102  3.6351659 3.6968626}%
\special{ar 3466 3608 288 102  3.8819526 3.9436492}%
\special{ar 3466 3608 288 102  4.1287392 4.1904358}%
\special{ar 3466 3608 288 102  4.3755258 4.4372225}%
\special{ar 3466 3608 288 102  4.6223124 4.6840091}%
\special{ar 3466 3608 288 102  4.8690991 4.9307957}%
\special{ar 3466 3608 288 102  5.1158857 5.1775824}%
\special{ar 3466 3608 288 102  5.3626723 5.4243690}%
\special{ar 3466 3608 288 102  5.6094590 5.6711556}%
\special{ar 3466 3608 288 102  5.8562456 5.9179423}%
\special{ar 3466 3608 288 102  6.1030322 6.1647289}%
%
\special{pn 8}%
\special{pa 3466 3142}%
\special{pa 3180 3600}%
\special{fp}%
%
\special{pn 8}%
\special{pa 3466 3142}%
\special{pa 3754 3600}%
\special{fp}%
%
\special{pn 8}%
\special{pa 3464 3142}%
\special{pa 3434 3706}%
\special{fp}%
%
\special{pn 13}%
\special{pa 1174 3658}%
\special{pa 1182 3628}%
\special{pa 1192 3596}%
\special{pa 1200 3566}%
\special{pa 1210 3536}%
\special{pa 1220 3506}%
\special{pa 1230 3476}%
\special{pa 1240 3446}%
\special{pa 1254 3416}%
\special{pa 1266 3386}%
\special{pa 1278 3356}%
\special{pa 1292 3328}%
\special{pa 1306 3300}%
\special{pa 1318 3270}%
\special{pa 1334 3242}%
\special{pa 1350 3214}%
\special{pa 1366 3186}%
\special{pa 1382 3160}%
\special{pa 1388 3150}%
\special{sp}%
%
\special{pn 13}%
\special{ar 1022 3976 560 1108  5.4319385 6.0123742}%
%
\special{pn 13}%
\special{pa 3300 3684}%
\special{pa 3304 3652}%
\special{pa 3310 3620}%
\special{pa 3318 3590}%
\special{pa 3328 3558}%
\special{pa 3336 3528}%
\special{pa 3348 3498}%
\special{pa 3360 3468}%
\special{pa 3376 3440}%
\special{pa 3390 3412}%
\special{pa 3406 3384}%
\special{pa 3424 3358}%
\special{pa 3444 3334}%
\special{pa 3466 3310}%
\special{pa 3488 3288}%
\special{pa 3510 3266}%
\special{pa 3522 3258}%
\special{sp}%
%
\special{pn 8}%
\special{ar 3464 3274 72 48  3.6787619 3.8804426}%
\special{ar 3464 3274 72 48  4.4854846 4.6871653}%
\special{ar 3464 3274 72 48  5.2922073 5.4938879}%
%
\special{pn 13}%
\special{ar 3048 3724 568 604  5.3890781 6.2349993}%
%
\special{pn 8}%
\special{pa 1146 2394}%
\special{pa 1146 2318}%
\special{fp}%
\special{sh 1}%
\special{pa 1146 2318}%
\special{pa 1126 2384}%
\special{pa 1146 2370}%
\special{pa 1166 2384}%
\special{pa 1146 2318}%
\special{fp}%
%
\special{pn 8}%
\special{pa 1390 2720}%
\special{pa 1446 2720}%
\special{fp}%
\special{sh 1}%
\special{pa 1446 2720}%
\special{pa 1380 2700}%
\special{pa 1394 2720}%
\special{pa 1380 2740}%
\special{pa 1446 2720}%
\special{fp}%
%
\special{pn 8}%
\special{pa 3186 2414}%
\special{pa 3186 2312}%
\special{fp}%
\special{sh 1}%
\special{pa 3186 2312}%
\special{pa 3166 2380}%
\special{pa 3186 2366}%
\special{pa 3206 2380}%
\special{pa 3186 2312}%
\special{fp}%
%
\special{pn 8}%
\special{pa 3444 2732}%
\special{pa 3536 2732}%
\special{fp}%
\special{sh 1}%
\special{pa 3536 2732}%
\special{pa 3470 2712}%
\special{pa 3484 2732}%
\special{pa 3470 2752}%
\special{pa 3536 2732}%
\special{fp}%
%
\special{pn 8}%
\special{pa 3450 3414}%
\special{pa 3450 3482}%
\special{fp}%
\special{sh 1}%
\special{pa 3450 3482}%
\special{pa 3470 3416}%
\special{pa 3450 3430}%
\special{pa 3430 3416}%
\special{pa 3450 3482}%
\special{fp}%
%
\special{pn 8}%
\special{pa 1366 3388}%
\special{pa 1366 3442}%
\special{fp}%
\special{sh 1}%
\special{pa 1366 3442}%
\special{pa 1386 3374}%
\special{pa 1366 3388}%
\special{pa 1346 3374}%
\special{pa 1366 3442}%
\special{fp}%
%
\special{pn 8}%
\special{pa 1352 2872}%
\special{pa 1352 2964}%
\special{fp}%
%
\special{pn 8}%
\special{pa 1418 2872}%
\special{pa 1418 2964}%
\special{fp}%
%
\special{pn 8}%
\special{pa 3422 2882}%
\special{pa 3422 2974}%
\special{fp}%
%
\special{pn 8}%
\special{pa 3486 2882}%
\special{pa 3486 2974}%
\special{fp}%
\put(12.4400,-39.2100){\makebox(0,0)[lb]{$f(U_{\lambda})$}}%
\put(32.4900,-39.4000){\makebox(0,0)[lb]{$f_t(U_{\lambda})\,\,\,\,(t>0)$}}%
%
\special{pn 8}%
\special{pa 2050 2460}%
\special{pa 2778 2460}%
\special{dt 0.045}%
\special{sh 1}%
\special{pa 2778 2460}%
\special{pa 2710 2440}%
\special{pa 2724 2460}%
\special{pa 2710 2480}%
\special{pa 2778 2460}%
\special{fp}%
%
\special{pn 8}%
\special{pa 1146 2450}%
\special{pa 1146 2372}%
\special{fp}%
\special{sh 1}%
\special{pa 1146 2372}%
\special{pa 1126 2438}%
\special{pa 1146 2424}%
\special{pa 1166 2438}%
\special{pa 1146 2372}%
\special{fp}%
\special{pa 1146 2372}%
\special{pa 1146 2372}%
\special{fp}%
%
\special{pn 8}%
\special{pa 1410 2720}%
\special{pa 1490 2720}%
\special{fp}%
\special{sh 1}%
\special{pa 1490 2720}%
\special{pa 1422 2700}%
\special{pa 1436 2720}%
\special{pa 1422 2740}%
\special{pa 1490 2720}%
\special{fp}%
%
\special{pn 8}%
\special{pa 3186 2466}%
\special{pa 3186 2366}%
\special{fp}%
\special{sh 1}%
\special{pa 3186 2366}%
\special{pa 3166 2434}%
\special{pa 3186 2420}%
\special{pa 3206 2434}%
\special{pa 3186 2366}%
\special{fp}%
%
\special{pn 8}%
\special{pa 3408 2732}%
\special{pa 3494 2732}%
\special{fp}%
\special{sh 1}%
\special{pa 3494 2732}%
\special{pa 3428 2712}%
\special{pa 3442 2732}%
\special{pa 3428 2752}%
\special{pa 3494 2732}%
\special{fp}%
%
\special{pn 8}%
\special{pa 1366 3420}%
\special{pa 1366 3482}%
\special{fp}%
\special{sh 1}%
\special{pa 1366 3482}%
\special{pa 1386 3414}%
\special{pa 1366 3428}%
\special{pa 1346 3414}%
\special{pa 1366 3482}%
\special{fp}%
%
\special{pn 8}%
\special{pa 3450 3460}%
\special{pa 3450 3522}%
\special{fp}%
\special{sh 1}%
\special{pa 3450 3522}%
\special{pa 3470 3454}%
\special{pa 3450 3468}%
\special{pa 3430 3454}%
\special{pa 3450 3522}%
\special{fp}%
%
\special{pn 8}%
\special{pa 690 1226}%
\special{pa 1740 1226}%
\special{fp}%
%
\special{pn 8}%
\special{pa 1144 1812}%
\special{pa 1144 630}%
\special{fp}%
%
\special{pn 8}%
\special{pa 3100 1230}%
\special{pa 4060 1230}%
\special{fp}%
%
\special{pn 8}%
\special{pa 3416 1834}%
\special{pa 3416 652}%
\special{fp}%
%
\special{pn 8}%
\special{pa 2070 1218}%
\special{pa 2790 1218}%
\special{dt 0.045}%
\special{sh 1}%
\special{pa 2790 1218}%
\special{pa 2724 1198}%
\special{pa 2738 1218}%
\special{pa 2724 1238}%
\special{pa 2790 1218}%
\special{fp}%
%
\special{pn 13}%
\special{ar 1654 720 722 722  1.4867395 3.2259088}%
%
\special{pn 13}%
\special{ar 4100 560 840 840  1.6876702 3.0132189}%
\put(10.3000,-5.2000){\makebox(0,0)[lb]{$\widehat f_{\lambda,\mu}(\widehat U_{\lambda})$}}%
\put(31.1000,-5.3000){\makebox(0,0)[lb]{$(\widehat f_t)_{\lambda,\mu}(\widehat U_{\lambda})\,\,\,\,(t>0)$}}%
\put(20.7000,-12.9000){\makebox(0,0)[lt]{{\scriptsize time goes by}}}%
\put(20.4000,-25.4000){\makebox(0,0)[lt]{{\scriptsize time goes by}}}%
\end{picture}%
\hspace{0.7truecm}}

\vspace{0.5truecm}

\centerline{{\bf Figure 1(Continued).}}

\newpage


\centerline{
\unitlength 0.1in
\begin{picture}( 39.8700, 23.9000)(  7.3700,-27.4000)
%
\special{pn 8}%
\special{pa 2570 940}%
\special{pa 4080 940}%
\special{fp}%
\put(33.1000,-5.2000){\makebox(0,0)[lb]{$\widehat V_{\mu}$}}%
%
\special{pn 8}%
\special{ar 2984 2492 848 434  4.8924697 4.9111905}%
\special{ar 2984 2492 848 434  4.9673527 4.9860735}%
\special{ar 2984 2492 848 434  5.0422357 5.0609565}%
\special{ar 2984 2492 848 434  5.1171187 5.1358394}%
\special{ar 2984 2492 848 434  5.1920017 5.2107224}%
\special{ar 2984 2492 848 434  5.2668847 5.2856054}%
\special{ar 2984 2492 848 434  5.3417677 5.3604884}%
\special{ar 2984 2492 848 434  5.4166507 5.4353714}%
\special{ar 2984 2492 848 434  5.4915337 5.5102544}%
\special{ar 2984 2492 848 434  5.5664167 5.5851374}%
\special{ar 2984 2492 848 434  5.6412997 5.6600204}%
\special{ar 2984 2492 848 434  5.7161827 5.7349034}%
\special{ar 2984 2492 848 434  5.7910657 5.8097864}%
\special{ar 2984 2492 848 434  5.8659487 5.8846694}%
\special{ar 2984 2492 848 434  5.9408316 5.9595524}%
\special{ar 2984 2492 848 434  6.0157146 6.0344354}%
\special{ar 2984 2492 848 434  6.0905976 6.1093184}%
\special{ar 2984 2492 848 434  6.1654806 6.1842014}%
\special{ar 2984 2492 848 434  6.2403636 6.2590844}%
%
\special{pn 8}%
\special{pa 2960 2490}%
\special{pa 3836 2490}%
\special{fp}%
%
\special{pn 8}%
\special{pa 2974 2492}%
\special{pa 3132 2060}%
\special{fp}%
%
\special{pn 8}%
\special{pa 3070 2210}%
\special{pa 3108 2108}%
\special{fp}%
\special{sh 1}%
\special{pa 3108 2108}%
\special{pa 3066 2164}%
\special{pa 3090 2158}%
\special{pa 3104 2176}%
\special{pa 3108 2108}%
\special{fp}%
%
\special{pn 8}%
\special{pa 3640 2490}%
\special{pa 3744 2490}%
\special{fp}%
\special{sh 1}%
\special{pa 3744 2490}%
\special{pa 3676 2470}%
\special{pa 3690 2490}%
\special{pa 3676 2510}%
\special{pa 3744 2490}%
\special{fp}%
\put(32.9000,-27.4000){\makebox(0,0)[lt]{$V_{\mu}$}}%
%
\special{pn 8}%
\special{pa 3396 1450}%
\special{pa 3396 1826}%
\special{fp}%
\special{sh 1}%
\special{pa 3396 1826}%
\special{pa 3416 1760}%
\special{pa 3396 1774}%
\special{pa 3376 1760}%
\special{pa 3396 1826}%
\special{fp}%
\put(35.4500,-15.5400){\makebox(0,0)[lt]{$\psi_{\mu}^{-1}\circ\pi_{\Gamma'_{\mu}}$}}%
%
\special{pn 13}%
\special{pa 4280 830}%
\special{pa 4430 1010}%
\special{fp}%
\special{sh 1}%
\special{pa 4430 1010}%
\special{pa 4402 946}%
\special{pa 4396 968}%
\special{pa 4372 972}%
\special{pa 4430 1010}%
\special{fp}%
\put(44.8500,-8.7700){\makebox(0,0)[lt]{'s$\,\,:\,\,\widehat H_{\lambda,\mu}$}}%
%
\special{pn 13}%
\special{pa 4160 2260}%
\special{pa 4310 2440}%
\special{fp}%
\special{sh 1}%
\special{pa 4310 2440}%
\special{pa 4282 2376}%
\special{pa 4276 2398}%
\special{pa 4252 2402}%
\special{pa 4310 2440}%
\special{fp}%
\put(43.6500,-23.0700){\makebox(0,0)[lt]{'s$\,\,:\,\,H$}}%
%
\special{pn 8}%
\special{pa 3590 2490}%
\special{pa 3694 2490}%
\special{fp}%
\special{sh 1}%
\special{pa 3694 2490}%
\special{pa 3626 2470}%
\special{pa 3640 2490}%
\special{pa 3626 2510}%
\special{pa 3694 2490}%
\special{fp}%
%
\special{pn 8}%
\special{pa 3050 2270}%
\special{pa 3088 2168}%
\special{fp}%
\special{sh 1}%
\special{pa 3088 2168}%
\special{pa 3046 2224}%
\special{pa 3070 2218}%
\special{pa 3084 2236}%
\special{pa 3088 2168}%
\special{fp}%
\put(12.8000,-7.8000){\makebox(0,0)[lb]{$\widehat U_{\lambda}$}}%
\put(20.7000,-7.8000){\makebox(0,0)[lt]{$\displaystyle{\mathop{\hookrightarrow}^{\widehat f_{\lambda,\mu}}}$}}%
%
\special{pn 13}%
\special{ar 1320 800 650 202  0.4185692 2.6833195}%
%
\special{pn 8}%
\special{pa 1370 1500}%
\special{pa 1370 1876}%
\special{fp}%
\special{sh 1}%
\special{pa 1370 1876}%
\special{pa 1390 1810}%
\special{pa 1370 1824}%
\special{pa 1350 1810}%
\special{pa 1370 1876}%
\special{fp}%
\put(14.6500,-15.8400){\makebox(0,0)[lt]{$\phi_{\lambda}^{-1}\circ\pi_{\Gamma_{\lambda}}$}}%
\put(11.8000,-26.3000){\makebox(0,0)[lt]{$U_{\lambda}$}}%
\put(23.4000,-21.9000){\makebox(0,0)[lt]{$\displaystyle{\mathop{\hookrightarrow}^{f}}$}}%
%
\special{pn 8}%
\special{pa 3020 1180}%
\special{pa 3430 610}%
\special{fp}%
\put(47.2400,-10.8700){\makebox(0,0)[lt]{(a local lift of $H$)}}%
%
\special{pn 13}%
\special{ar 1370 2120 808 182  0.5918768 1.5707963}%
%
\special{pn 20}%
\special{sh 1}%
\special{ar 1360 2290 10 10 0  6.28318530717959E+0000}%
\special{sh 1}%
\special{ar 1360 2290 10 10 0  6.28318530717959E+0000}%
%
\special{pn 13}%
\special{pa 3380 1040}%
\special{pa 3540 900}%
\special{fp}%
\special{sh 1}%
\special{pa 3540 900}%
\special{pa 3478 930}%
\special{pa 3500 936}%
\special{pa 3504 960}%
\special{pa 3540 900}%
\special{fp}%
%
\special{pn 13}%
\special{ar 3630 710 540 370  1.0675219 2.4619273}%
%
\special{pn 13}%
\special{ar 2790 1160 540 370  4.2091145 5.6035199}%
%
\special{pn 13}%
\special{pa 3630 1070}%
\special{pa 3770 840}%
\special{fp}%
\special{sh 1}%
\special{pa 3770 840}%
\special{pa 3718 888}%
\special{pa 3742 886}%
\special{pa 3752 908}%
\special{pa 3770 840}%
\special{fp}%
%
\special{pn 13}%
\special{pa 3280 980}%
\special{pa 3400 910}%
\special{fp}%
\special{sh 1}%
\special{pa 3400 910}%
\special{pa 3332 926}%
\special{pa 3354 938}%
\special{pa 3352 962}%
\special{pa 3400 910}%
\special{fp}%
%
\special{pn 13}%
\special{pa 2750 800}%
\special{pa 2580 1040}%
\special{fp}%
\special{sh 1}%
\special{pa 2580 1040}%
\special{pa 2636 998}%
\special{pa 2612 996}%
\special{pa 2602 974}%
\special{pa 2580 1040}%
\special{fp}%
%
\special{pn 13}%
\special{pa 2980 820}%
\special{pa 2810 990}%
\special{fp}%
\special{sh 1}%
\special{pa 2810 990}%
\special{pa 2872 958}%
\special{pa 2848 952}%
\special{pa 2844 930}%
\special{pa 2810 990}%
\special{fp}%
%
\special{pn 13}%
\special{pa 3110 870}%
\special{pa 2990 960}%
\special{fp}%
\special{sh 1}%
\special{pa 2990 960}%
\special{pa 3056 936}%
\special{pa 3034 928}%
\special{pa 3032 904}%
\special{pa 2990 960}%
\special{fp}%
%
\special{pn 13}%
\special{pa 2960 2500}%
\special{pa 2986 2480}%
\special{pa 3010 2460}%
\special{pa 3038 2444}%
\special{pa 3064 2426}%
\special{pa 3092 2412}%
\special{pa 3122 2398}%
\special{pa 3150 2386}%
\special{pa 3182 2376}%
\special{pa 3212 2366}%
\special{pa 3242 2356}%
\special{pa 3274 2352}%
\special{pa 3304 2344}%
\special{pa 3336 2338}%
\special{pa 3368 2336}%
\special{pa 3400 2332}%
\special{pa 3432 2332}%
\special{pa 3464 2330}%
\special{pa 3496 2330}%
\special{pa 3528 2332}%
\special{pa 3560 2332}%
\special{pa 3592 2334}%
\special{pa 3624 2340}%
\special{pa 3654 2346}%
\special{pa 3686 2350}%
\special{pa 3718 2358}%
\special{pa 3748 2364}%
\special{pa 3778 2374}%
\special{pa 3810 2384}%
\special{pa 3810 2384}%
\special{sp}%
%
\special{pn 13}%
\special{pa 3540 2320}%
\special{pa 3460 2610}%
\special{fp}%
\special{sh 1}%
\special{pa 3460 2610}%
\special{pa 3498 2552}%
\special{pa 3474 2560}%
\special{pa 3458 2540}%
\special{pa 3460 2610}%
\special{fp}%
%
\special{pn 13}%
\special{pa 3290 2350}%
\special{pa 3290 2560}%
\special{fp}%
\special{sh 1}%
\special{pa 3290 2560}%
\special{pa 3310 2494}%
\special{pa 3290 2508}%
\special{pa 3270 2494}%
\special{pa 3290 2560}%
\special{fp}%
%
\special{pn 13}%
\special{pa 3110 2410}%
\special{pa 3170 2530}%
\special{fp}%
\special{sh 1}%
\special{pa 3170 2530}%
\special{pa 3158 2462}%
\special{pa 3146 2482}%
\special{pa 3122 2480}%
\special{pa 3170 2530}%
\special{fp}%
\end{picture}%
\hspace{3truecm}}

\vspace{0.5truecm}

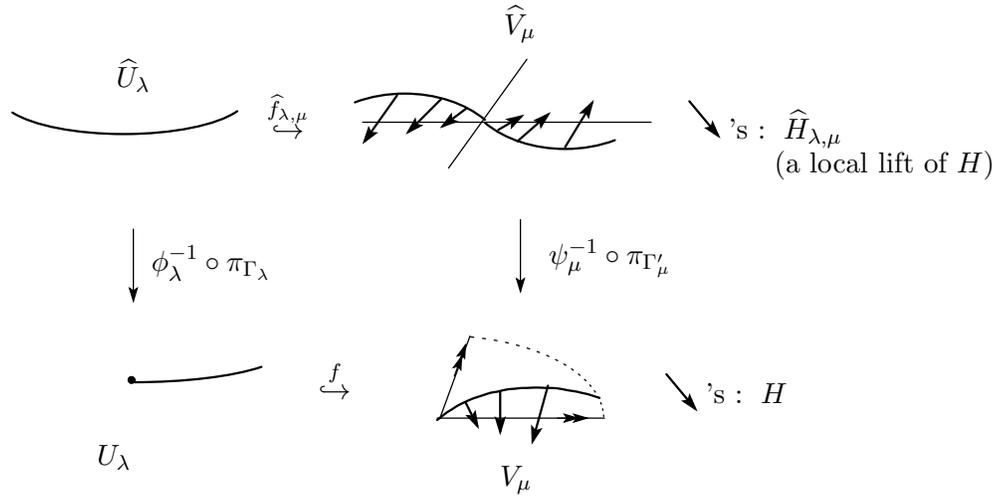
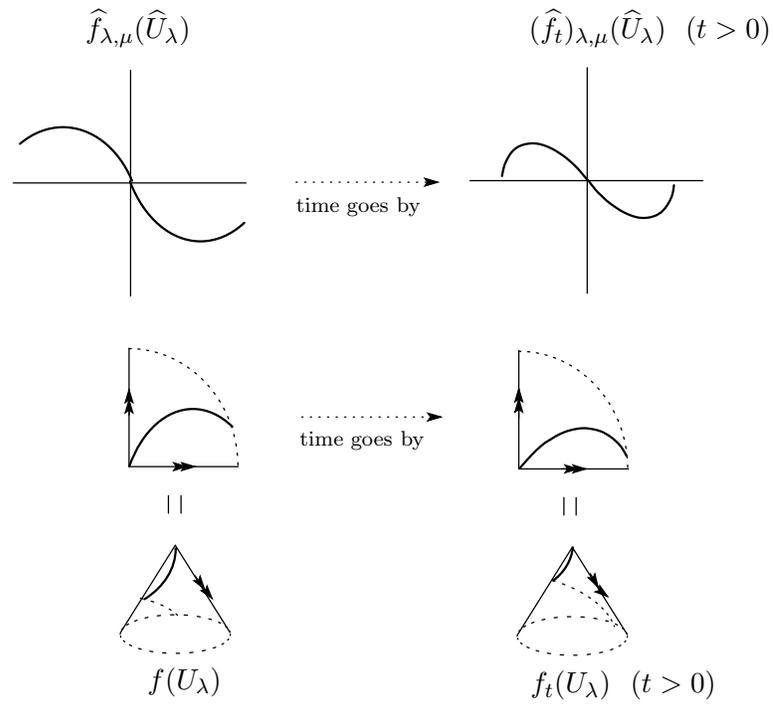
\centerline{{\bf Figure 2.}}

\vspace{0.5truecm}

\centerline{
\unitlength 0.1in
\begin{picture}( 36.0900, 34.2000)(  5.4000,-37.7000)
%
\special{pn 8}%
\special{pa 1146 2720}%
\special{pa 1146 2100}%
\special{fp}%
%
\special{pn 8}%
\special{pa 1146 2720}%
\special{pa 1718 2720}%
\special{fp}%
%
\special{pn 8}%
\special{ar 1146 2720 570 618  4.7123890 4.7326080}%
\special{ar 1146 2720 570 618  4.7932651 4.8134842}%
\special{ar 1146 2720 570 618  4.8741413 4.8943603}%
\special{ar 1146 2720 570 618  4.9550175 4.9752365}%
\special{ar 1146 2720 570 618  5.0358936 5.0561127}%
\special{ar 1146 2720 570 618  5.1167698 5.1369888}%
\special{ar 1146 2720 570 618  5.1976459 5.2178650}%
\special{ar 1146 2720 570 618  5.2785221 5.2987411}%
\special{ar 1146 2720 570 618  5.3593982 5.3796173}%
\special{ar 1146 2720 570 618  5.4402744 5.4604934}%
\special{ar 1146 2720 570 618  5.5211506 5.5413696}%
\special{ar 1146 2720 570 618  5.6020267 5.6222458}%
\special{ar 1146 2720 570 618  5.6829029 5.7031219}%
\special{ar 1146 2720 570 618  5.7637790 5.7839981}%
\special{ar 1146 2720 570 618  5.8446552 5.8648742}%
\special{ar 1146 2720 570 618  5.9255314 5.9457504}%
\special{ar 1146 2720 570 618  6.0064075 6.0266266}%
\special{ar 1146 2720 570 618  6.0872837 6.1075027}%
\special{ar 1146 2720 570 618  6.1681598 6.1883789}%
\special{ar 1146 2720 570 618  6.2490360 6.2692550}%
%
\special{pn 8}%
\special{pa 3186 2732}%
\special{pa 3186 2112}%
\special{fp}%
%
\special{pn 8}%
\special{pa 3186 2732}%
\special{pa 3760 2732}%
\special{fp}%
%
\special{pn 8}%
\special{ar 3186 2732 570 616  4.7123890 4.7326251}%
\special{ar 3186 2732 570 616  4.7933333 4.8135694}%
\special{ar 3186 2732 570 616  4.8742777 4.8945138}%
\special{ar 3186 2732 570 616  4.9552220 4.9754581}%
\special{ar 3186 2732 570 616  5.0361664 5.0564025}%
\special{ar 3186 2732 570 616  5.1171107 5.1373468}%
\special{ar 3186 2732 570 616  5.1980551 5.2182912}%
\special{ar 3186 2732 570 616  5.2789994 5.2992355}%
\special{ar 3186 2732 570 616  5.3599438 5.3801799}%
\special{ar 3186 2732 570 616  5.4408881 5.4611242}%
\special{ar 3186 2732 570 616  5.5218325 5.5420686}%
\special{ar 3186 2732 570 616  5.6027768 5.6230129}%
\special{ar 3186 2732 570 616  5.6837212 5.7039573}%
\special{ar 3186 2732 570 616  5.7646655 5.7849016}%
\special{ar 3186 2732 570 616  5.8456099 5.8658460}%
\special{ar 3186 2732 570 616  5.9265542 5.9467903}%
\special{ar 3186 2732 570 616  6.0074986 6.0277347}%
\special{ar 3186 2732 570 616  6.0884429 6.1086790}%
\special{ar 3186 2732 570 616  6.1693873 6.1896234}%
\special{ar 3186 2732 570 616  6.2503316 6.2705677}%
%
\special{pn 8}%
\special{ar 1388 3598 288 100  6.2831853 6.3453615}%
\special{ar 1388 3598 288 100  6.5318900 6.5940661}%
\special{ar 1388 3598 288 100  6.7805946 6.8427708}%
\special{ar 1388 3598 288 100  7.0292993 7.0914755}%
\special{ar 1388 3598 288 100  7.2780040 7.3401801}%
\special{ar 1388 3598 288 100  7.5267086 7.5888848}%
\special{ar 1388 3598 288 100  7.7754133 7.8375895}%
\special{ar 1388 3598 288 100  8.0241179 8.0862941}%
\special{ar 1388 3598 288 100  8.2728226 8.3349988}%
\special{ar 1388 3598 288 100  8.5215273 8.5837034}%
\special{ar 1388 3598 288 100  8.7702319 8.8324081}%
\special{ar 1388 3598 288 100  9.0189366 9.0811128}%
\special{ar 1388 3598 288 100  9.2676413 9.3298174}%
%
\special{pn 8}%
\special{ar 1388 3598 288 102  3.1415927 3.2034483}%
\special{ar 1388 3598 288 102  3.3890153 3.4508710}%
\special{ar 1388 3598 288 102  3.6364380 3.6982937}%
\special{ar 1388 3598 288 102  3.8838607 3.9457164}%
\special{ar 1388 3598 288 102  4.1312834 4.1931390}%
\special{ar 1388 3598 288 102  4.3787061 4.4405617}%
\special{ar 1388 3598 288 102  4.6261287 4.6879844}%
\special{ar 1388 3598 288 102  4.8735514 4.9354071}%
\special{ar 1388 3598 288 102  5.1209741 5.1828298}%
\special{ar 1388 3598 288 102  5.3683968 5.4302524}%
\special{ar 1388 3598 288 102  5.6158195 5.6776751}%
\special{ar 1388 3598 288 102  5.8632421 5.9250978}%
\special{ar 1388 3598 288 102  6.1106648 6.1725205}%
%
\special{pn 8}%
\special{pa 1388 3132}%
\special{pa 1102 3588}%
\special{fp}%
%
\special{pn 8}%
\special{pa 1388 3132}%
\special{pa 1674 3588}%
\special{fp}%
%
\special{pn 8}%
\special{ar 3466 3608 288 100  6.2831853 6.3450410}%
\special{ar 3466 3608 288 100  6.5306080 6.5924637}%
\special{ar 3466 3608 288 100  6.7780307 6.8398863}%
\special{ar 3466 3608 288 100  7.0254533 7.0873090}%
\special{ar 3466 3608 288 100  7.2728760 7.3347317}%
\special{ar 3466 3608 288 100  7.5202987 7.5821544}%
\special{ar 3466 3608 288 100  7.7677214 7.8295771}%
\special{ar 3466 3608 288 100  8.0151441 8.0769997}%
\special{ar 3466 3608 288 100  8.2625668 8.3244224}%
\special{ar 3466 3608 288 100  8.5099894 8.5718451}%
\special{ar 3466 3608 288 100  8.7574121 8.8192678}%
\special{ar 3466 3608 288 100  9.0048348 9.0666905}%
\special{ar 3466 3608 288 100  9.2522575 9.3141131}%
%
\special{pn 8}%
\special{ar 3466 3608 288 102  3.1415927 3.2032893}%
\special{ar 3466 3608 288 102  3.3883793 3.4500759}%
\special{ar 3466 3608 288 102  3.6351659 3.6968626}%
\special{ar 3466 3608 288 102  3.8819526 3.9436492}%
\special{ar 3466 3608 288 102  4.1287392 4.1904358}%
\special{ar 3466 3608 288 102  4.3755258 4.4372225}%
\special{ar 3466 3608 288 102  4.6223124 4.6840091}%
\special{ar 3466 3608 288 102  4.8690991 4.9307957}%
\special{ar 3466 3608 288 102  5.1158857 5.1775824}%
\special{ar 3466 3608 288 102  5.3626723 5.4243690}%
\special{ar 3466 3608 288 102  5.6094590 5.6711556}%
\special{ar 3466 3608 288 102  5.8562456 5.9179423}%
\special{ar 3466 3608 288 102  6.1030322 6.1647289}%
%
\special{pn 8}%
\special{pa 3466 3142}%
\special{pa 3180 3600}%
\special{fp}%
%
\special{pn 8}%
\special{pa 3466 3142}%
\special{pa 3754 3600}%
\special{fp}%
%
\special{pn 8}%
\special{pa 1146 2394}%
\special{pa 1146 2318}%
\special{fp}%
\special{sh 1}%
\special{pa 1146 2318}%
\special{pa 1126 2384}%
\special{pa 1146 2370}%
\special{pa 1166 2384}%
\special{pa 1146 2318}%
\special{fp}%
%
\special{pn 8}%
\special{pa 1390 2720}%
\special{pa 1446 2720}%
\special{fp}%
\special{sh 1}%
\special{pa 1446 2720}%
\special{pa 1380 2700}%
\special{pa 1394 2720}%
\special{pa 1380 2740}%
\special{pa 1446 2720}%
\special{fp}%
%
\special{pn 8}%
\special{pa 3186 2414}%
\special{pa 3186 2312}%
\special{fp}%
\special{sh 1}%
\special{pa 3186 2312}%
\special{pa 3166 2380}%
\special{pa 3186 2366}%
\special{pa 3206 2380}%
\special{pa 3186 2312}%
\special{fp}%
%
\special{pn 8}%
\special{pa 3444 2732}%
\special{pa 3536 2732}%
\special{fp}%
\special{sh 1}%
\special{pa 3536 2732}%
\special{pa 3470 2712}%
\special{pa 3484 2732}%
\special{pa 3470 2752}%
\special{pa 3536 2732}%
\special{fp}%
%
\special{pn 8}%
\special{pa 1352 2872}%
\special{pa 1352 2964}%
\special{fp}%
%
\special{pn 8}%
\special{pa 1418 2872}%
\special{pa 1418 2964}%
\special{fp}%
%
\special{pn 8}%
\special{pa 3422 2882}%
\special{pa 3422 2974}%
\special{fp}%
%
\special{pn 8}%
\special{pa 3486 2882}%
\special{pa 3486 2974}%
\special{fp}%
\put(12.4400,-39.2100){\makebox(0,0)[lb]{$f(U_{\lambda})$}}%
\put(32.4900,-39.4000){\makebox(0,0)[lb]{$f_t(U_{\lambda})\,\,\,\,(t>0)$}}%
%
\special{pn 8}%
\special{pa 2050 2460}%
\special{pa 2778 2460}%
\special{dt 0.045}%
\special{sh 1}%
\special{pa 2778 2460}%
\special{pa 2710 2440}%
\special{pa 2724 2460}%
\special{pa 2710 2480}%
\special{pa 2778 2460}%
\special{fp}%
%
\special{pn 8}%
\special{pa 1146 2450}%
\special{pa 1146 2372}%
\special{fp}%
\special{sh 1}%
\special{pa 1146 2372}%
\special{pa 1126 2438}%
\special{pa 1146 2424}%
\special{pa 1166 2438}%
\special{pa 1146 2372}%
\special{fp}%
\special{pa 1146 2372}%
\special{pa 1146 2372}%
\special{fp}%
%
\special{pn 8}%
\special{pa 1410 2720}%
\special{pa 1490 2720}%
\special{fp}%
\special{sh 1}%
\special{pa 1490 2720}%
\special{pa 1422 2700}%
\special{pa 1436 2720}%
\special{pa 1422 2740}%
\special{pa 1490 2720}%
\special{fp}%
%
\special{pn 8}%
\special{pa 3186 2466}%
\special{pa 3186 2366}%
\special{fp}%
\special{sh 1}%
\special{pa 3186 2366}%
\special{pa 3166 2434}%
\special{pa 3186 2420}%
\special{pa 3206 2434}%
\special{pa 3186 2366}%
\special{fp}%
%
\special{pn 8}%
\special{pa 3408 2732}%
\special{pa 3494 2732}%
\special{fp}%
\special{sh 1}%
\special{pa 3494 2732}%
\special{pa 3428 2712}%
\special{pa 3442 2732}%
\special{pa 3428 2752}%
\special{pa 3494 2732}%
\special{fp}%
%
\special{pn 8}%
\special{pa 540 1234}%
\special{pa 1758 1234}%
\special{fp}%
%
\special{pn 8}%
\special{pa 1154 1826}%
\special{pa 1154 644}%
\special{fp}%
%
\special{pn 8}%
\special{pa 2930 1222}%
\special{pa 4150 1222}%
\special{fp}%
%
\special{pn 8}%
\special{pa 3544 1812}%
\special{pa 3544 630}%
\special{fp}%
%
\special{pn 8}%
\special{pa 2020 1230}%
\special{pa 2762 1230}%
\special{dt 0.045}%
\special{sh 1}%
\special{pa 2762 1230}%
\special{pa 2694 1210}%
\special{pa 2708 1230}%
\special{pa 2694 1250}%
\special{pa 2762 1230}%
\special{fp}%
\put(9.1900,-5.2000){\makebox(0,0)[lb]{$\widehat f_{\lambda,\mu}(\widehat U_{\lambda})$}}%
\put(32.4000,-5.2000){\makebox(0,0)[lb]{$(\widehat f_t)_{\lambda,\mu}(\widehat U_{\lambda})\,\,\,\,(t>0)$}}%
%
\special{pn 13}%
\special{ar 1520 940 420 600  0.9939398 2.6435516}%
%
\special{pn 13}%
\special{ar 800 1480 412 538  4.1353393 5.7843726}%
%
\special{pn 13}%
\special{pa 3998 1246}%
\special{pa 3994 1278}%
\special{pa 3986 1308}%
\special{pa 3972 1336}%
\special{pa 3954 1364}%
\special{pa 3932 1388}%
\special{pa 3906 1404}%
\special{pa 3876 1414}%
\special{pa 3844 1418}%
\special{pa 3812 1416}%
\special{pa 3782 1408}%
\special{pa 3752 1396}%
\special{pa 3724 1382}%
\special{pa 3696 1364}%
\special{pa 3670 1346}%
\special{pa 3644 1326}%
\special{pa 3622 1304}%
\special{pa 3598 1282}%
\special{pa 3578 1258}%
\special{pa 3558 1232}%
\special{pa 3550 1222}%
\special{sp}%
%
\special{pn 13}%
\special{pa 3098 1198}%
\special{pa 3104 1166}%
\special{pa 3112 1136}%
\special{pa 3126 1108}%
\special{pa 3142 1080}%
\special{pa 3164 1056}%
\special{pa 3190 1040}%
\special{pa 3222 1030}%
\special{pa 3254 1026}%
\special{pa 3284 1028}%
\special{pa 3316 1036}%
\special{pa 3346 1048}%
\special{pa 3374 1062}%
\special{pa 3402 1080}%
\special{pa 3428 1098}%
\special{pa 3452 1118}%
\special{pa 3476 1140}%
\special{pa 3498 1162}%
\special{pa 3518 1186}%
\special{pa 3540 1212}%
\special{pa 3548 1222}%
\special{sp}%
%
\special{pn 13}%
\special{ar 1480 2990 378 572  3.6393456 5.2897678}%
%
\special{pn 8}%
\special{pa 1490 3290}%
\special{pa 1540 3360}%
\special{fp}%
\special{sh 1}%
\special{pa 1540 3360}%
\special{pa 1518 3294}%
\special{pa 1510 3318}%
\special{pa 1486 3318}%
\special{pa 1540 3360}%
\special{fp}%
%
\special{pn 8}%
\special{pa 1540 3360}%
\special{pa 1580 3410}%
\special{fp}%
\special{sh 1}%
\special{pa 1580 3410}%
\special{pa 1554 3346}%
\special{pa 1548 3368}%
\special{pa 1524 3370}%
\special{pa 1580 3410}%
\special{fp}%
%
\special{pn 8}%
\special{pa 3560 3280}%
\special{pa 3600 3350}%
\special{fp}%
\special{sh 1}%
\special{pa 3600 3350}%
\special{pa 3584 3282}%
\special{pa 3574 3304}%
\special{pa 3550 3302}%
\special{pa 3600 3350}%
\special{fp}%
%
\special{pn 8}%
\special{pa 3610 3370}%
\special{pa 3650 3400}%
\special{fp}%
\special{sh 1}%
\special{pa 3650 3400}%
\special{pa 3610 3344}%
\special{pa 3608 3368}%
\special{pa 3586 3376}%
\special{pa 3650 3400}%
\special{fp}%
%
\special{pn 13}%
\special{ar 950 3150 440 340  6.2831853 6.2831853}%
\special{ar 950 3150 440 340  0.0000000 0.8888526}%
%
\special{pn 8}%
\special{ar 1020 3580 420 190  5.1687753 5.2081196}%
\special{ar 1020 3580 420 190  5.3261524 5.3654966}%
\special{ar 1020 3580 420 190  5.4835294 5.5228737}%
\special{ar 1020 3580 420 190  5.6409065 5.6802507}%
\special{ar 1020 3580 420 190  5.7982835 5.8376278}%
%
\special{pn 13}%
\special{ar 3170 3120 300 260  0.0899527 0.8571312}%
%
\special{pn 8}%
\special{ar 3040 3630 660 360  5.2217428 5.2452722}%
\special{ar 3040 3630 660 360  5.3158604 5.3393898}%
\special{ar 3040 3630 660 360  5.4099781 5.4335075}%
\special{ar 3040 3630 660 360  5.5040957 5.5276251}%
\special{ar 3040 3630 660 360  5.5982134 5.6217428}%
\special{ar 3040 3630 660 360  5.6923310 5.7158604}%
\special{ar 3040 3630 660 360  5.7864487 5.8099781}%
\special{ar 3040 3630 660 360  5.8805663 5.9040957}%
\special{ar 3040 3630 660 360  5.9746840 5.9982134}%
\special{ar 3040 3630 660 360  6.0688016 6.0797874}%
%
\special{pn 13}%
\special{pa 3190 2730}%
\special{pa 3210 2706}%
\special{pa 3232 2682}%
\special{pa 3254 2658}%
\special{pa 3276 2636}%
\special{pa 3300 2614}%
\special{pa 3324 2594}%
\special{pa 3352 2578}%
\special{pa 3380 2560}%
\special{pa 3408 2546}%
\special{pa 3438 2534}%
\special{pa 3468 2526}%
\special{pa 3500 2520}%
\special{pa 3532 2518}%
\special{pa 3564 2522}%
\special{pa 3594 2528}%
\special{pa 3624 2538}%
\special{pa 3652 2554}%
\special{pa 3678 2572}%
\special{pa 3702 2594}%
\special{pa 3724 2618}%
\special{pa 3740 2646}%
\special{pa 3750 2664}%
\special{sp}%
\put(20.2000,-13.1000){\makebox(0,0)[lt]{{\scriptsize time goes by}}}%
\put(20.4000,-25.3000){\makebox(0,0)[lt]{{\scriptsize time goes by}}}%
\end{picture}%
\hspace{0.7truecm}}

\vspace{0.5truecm}

\centerline{{\bf Figure 2(Continued).}}

\section{Evolution equations}
Let $G\curvearrowright V$ be an isometric almost free action with minimal 
regularizable orbit of a Hilbert Lie group $G$ on a Hilbert space $V$ equipped 
with an inner product $\langle\,\,,\,\,\rangle$.  
The orbit space $V/G$ is a (finite dimensional) $C^{\infty}$-orbifold.  
Let $\phi:V\to V/G$ be the orbit map and set $N:=V/G$.  
Here we give an example of such an isometric almost free action of a Hilbert Lie group.  

\vspace{0.5truecm}

\noindent
{\it Example.} 
Let $G$ be a compact semi-simple Lie group, $K$ a closed subgroup of $G$ and $\Gamma$ 
a finite subgroup of $G$.  
Denote by $\mathfrak g$ and $\mathfrak k$ the Lie algebras of $G$ and $K$, respectively.  
Assume that a reductive decomposition $\mathfrak g=\mathfrak k+\mathfrak p$ exists.  
Let $B$ be the Killing form of $\mathfrak g$.  
Give $G$ the bi-invariant metric induced from $B$.  
Let $H^0([0,1],\mathfrak g)$ be the Hilbert space of all paths in the Lie 
algebra $\mathfrak g$ of $G$ which are $L^2$-integrable with respect to $B$.  
Also, let $H^1([0,1],G)$ the Hilbert Lie group of all paths in $G$ which are 
of class $H^1$ with respect to $g$.  This group $H^1([0,1],G)$ acts on 
$H^0([0,1],\mathfrak g)$ isometrically and transitively as a gauge action:
$$\begin{array}{r}
\displaystyle{(a\ast u)(t)={\rm Ad}_G(a(t))(u(t))-(R_{a(t)})_{\ast}^{-1}(a'(t))
}\\
\displaystyle{(a\in H^1([0,1],G),\,u\in H^0([0,1],\mathfrak g)),}
\end{array}$$
where ${\rm Ad}_G$ is the adjoint representation of $G$ and $R_{a(t)}$ is 
the right translation by $a(t)$ and $a'$ is the weak derivative of $a$.  
Set $P(G,\Gamma\times K):=\{a\in H^1([0,1],G)\,\vert\,(a(0),a(1))$\newline
$\in\Gamma\times K\}$.  
The group $P(G,\Gamma\times K)$ acts on $H^0([0,1],\mathfrak g)$ almost freely 
and isometrically, and the orbit space of this action is 
diffeomorphic to the orbifold $\Gamma\setminus G\,/\,K$.  
Furthermore, each orbit of this action is regularizable and minimal.  

\vspace{0.5truecm}

\noindent
Give $N$ the Riemannian orbimetric such that $\phi$ is a Riemannian orbisubmersion.  
Let $f:M\hookrightarrow V$ be a $G$-invariant submanifold immersion such that 
$(\phi\circ f)(M)$ is compact.  
For this immersion $f$, we can take an orbiimmesion $\overline f$ of a compact orbifold $\overline M$ 
into $N$ and an orbifold submersion $\phi_M:M\to \overline M$ with 
$\phi\circ f=\overline f\circ\phi_M$.  
Let $\overline f_t$ ($0\leq t<T$) be the mean curvature flow starting from $\overline f$.  
The existenceness and the uniqueness of this flow in short time is assured by Proposition 3.1.  
Define a map $\overline F:\overline M\times[0,T)\to N$ by 
$\overline F(x,t):=\overline f_t(x)$ ($(x,t)\in \overline M\times[0,T)$).  
Denote by $H$ the regularized mean curvature vector of $f$ and $\overline H$ that of 
$\overline f$.  Since $\phi$ has minimal regularizable fibres, $H$ is the horizontal lift of 
$\overline H$.   Take $x\in\overline M$ and $u\in\phi_M^{-1}(x)$.  
Define a curve $c_x:[0,T)\to N$ by 
$c_x(t):=\overline f_t(x)$ and let $(c_x)_u^L:[0,T)\to V$ be the horizontal lift of $c_x$ 
to $f(u)$ satisfying $((c_x)_u^L)'(0)=H_u$.  Define an immersion $f_t:M\hookrightarrow V$ by 
$f_t(u)=(c_x)_u^L(t)$ ($u\in\widetilde M$) and a map 
$F:M\times[0,T)\to V$ by $F(u,t)=f_t(u)$ ($(u,t)\in M\times[0,T)$).  

\vspace{0.4truecm}

\noindent
{\bf Proposition 4.1.} {\sl The flow $f_t$ ($0\leq t<T$) is the regularized mean curvature flow 
starting from $f$.}

\vspace{0.4truecm}

\noindent
{\it Proof.} 
Denote by $\overline H_t$ the mean curvature vector of 
$\overline f_t$ and $H_t$ the regularized mean curvature vector of $f_t$.  
Fix $(u,t)\in M\times[0,T)$.  
It is clear that $\phi\circ f_t=\overline f_t\circ\phi_M$.  
Hence, since each fibre of $\phi$ is regularizable and minimal, 
$(H_t)_u$ coincides with one 
of the horizontal lifts of $(\overline H_t)_{\phi(u)}$ to $(c_{\phi(u)})_u^L(t)$.  
On the other hand,  from the definition of $F$, we have 
$\displaystyle{\frac{\partial F}{\partial t}(u,t)=((c_{\phi(u)})_u^L)'(t)}$, which is 
one of the horizontal lifts of $(\overline H_t)_{\phi(u)}$ to $(c_{\phi(u)})_u^L(t)$.  
These facts together with $\displaystyle{\frac{\partial F}{\partial t}(u,0)=H_u}$ implies that 
$\displaystyle{\frac{\partial F}{\partial t}(u,t)=(H_t)_u}$.  
Thus $f_t$ ($0\leq t<T$) is the regularized mean curvature flow starting from $f$.  
This completes the proof.\hspace{9.5truecm}q.e.d.

\vspace{0.4truecm}

\centerline{
\unitlength 0.1in
\begin{picture}( 48.0500, 36.7000)(  3.5000,-37.1000)
%
\special{pn 8}%
\special{ar 3614 902 930 268  0.0217751 6.2831853}%
%
\special{pn 8}%
\special{pa 2686 910}%
\special{pa 4542 910}%
\special{fp}%
%
\special{pn 8}%
\special{pa 3336 1158}%
\special{pa 3906 650}%
\special{fp}%
%
\special{pn 13}%
\special{ar 2606 494 1606 526  0.3544844 1.4642637}%
%
\special{pn 13}%
\special{ar 4782 1342 1780 576  3.5185243 4.5261219}%
%
\special{pn 13}%
\special{ar 4636 866 1048 198  1.8169763 2.9073920}%
%
\special{pn 13}%
\special{ar 3290 1212 352 792  5.8918210 6.2265449}%
%
\special{pn 13}%
\special{ar 2566 940 1062 172  4.9345790 6.0907585}%
%
\special{pn 8}%
\special{pa 2686 538}%
\special{pa 2686 2016}%
\special{fp}%
%
\special{pn 8}%
\special{ar 3614 1654 930 266  0.0154253 6.2831853}%
%
\special{pn 8}%
\special{pa 2686 1658}%
\special{pa 4542 1658}%
\special{fp}%
%
\special{pn 8}%
\special{pa 3336 1906}%
\special{pa 3906 1400}%
\special{fp}%
%
\special{pn 13}%
\special{ar 2606 1246 1606 524  0.3537376 1.4643446}%
%
\special{pn 13}%
\special{ar 4782 2092 1780 576  3.5195978 4.5261219}%
%
\special{pn 13}%
\special{ar 4636 1616 1048 198  1.8162541 2.9041926}%
%
\special{pn 13}%
\special{ar 3290 1962 352 790  5.8918210 6.2265449}%
%
\special{pn 13}%
\special{ar 2566 1690 1062 174  4.9371324 6.1020132}%
%
\special{pn 8}%
\special{pa 4542 532}%
\special{pa 4542 2010}%
\special{fp}%
%
\special{pn 13}%
\special{pa 2772 576}%
\special{pa 2772 2052}%
\special{fp}%
%
\special{pn 13}%
\special{pa 4396 626}%
\special{pa 4396 2102}%
\special{fp}%
%
\special{pn 13}%
\special{pa 2798 440}%
\special{pa 2804 890}%
\special{fp}%
%
\special{pn 13}%
\special{pa 2804 934}%
\special{pa 2798 996}%
\special{fp}%
%
\special{pn 13}%
\special{pa 2804 1076}%
\special{pa 2798 1626}%
\special{fp}%
%
\special{pn 13}%
\special{pa 2798 1676}%
\special{pa 2804 1744}%
\special{fp}%
%
\special{pn 13}%
\special{pa 2804 1806}%
\special{pa 2812 1940}%
\special{fp}%
%
\special{pn 13}%
\special{pa 4112 342}%
\special{pa 4104 786}%
\special{fp}%
%
\special{pn 13}%
\special{pa 4112 842}%
\special{pa 4104 890}%
\special{fp}%
%
\special{pn 13}%
\special{pa 4104 934}%
\special{pa 4112 1014}%
\special{fp}%
%
\special{pn 13}%
\special{pa 4112 1064}%
\special{pa 4112 1112}%
\special{fp}%
%
\special{pn 13}%
\special{pa 4112 1158}%
\special{pa 4104 1532}%
\special{fp}%
%
\special{pn 13}%
\special{pa 4104 1590}%
\special{pa 4112 1638}%
\special{fp}%
%
\special{pn 13}%
\special{pa 4112 1682}%
\special{pa 4112 1762}%
\special{fp}%
%
\special{pn 13}%
\special{pa 4112 1806}%
\special{pa 4104 1856}%
\special{fp}%
%
\special{pn 13}%
\special{pa 4112 1892}%
\special{pa 4118 1936}%
\special{fp}%
%
\special{pn 13}%
\special{pa 4450 434}%
\special{pa 4444 890}%
\special{fp}%
%
\special{pn 13}%
\special{pa 4450 928}%
\special{pa 4444 1002}%
\special{fp}%
%
\special{pn 13}%
\special{pa 4444 1046}%
\special{pa 4450 1620}%
\special{fp}%
%
\special{pn 13}%
\special{pa 4450 1682}%
\special{pa 4450 1738}%
\special{fp}%
%
\special{pn 13}%
\special{pa 4450 1794}%
\special{pa 4450 1922}%
\special{fp}%
%
\special{pn 13}%
\special{pa 3614 266}%
\special{pa 3614 848}%
\special{fp}%
%
\special{pn 13}%
\special{pa 3608 978}%
\special{pa 3608 1150}%
\special{fp}%
%
\special{pn 13}%
\special{pa 3614 1188}%
\special{pa 3608 1194}%
\special{fp}%
\special{pa 3614 1502}%
\special{pa 3614 1502}%
\special{fp}%
%
\special{pn 13}%
\special{pa 3608 1188}%
\special{pa 3614 1502}%
\special{fp}%
%
\special{pn 13}%
\special{pa 3640 996}%
\special{pa 3634 2214}%
\special{fp}%
%
\special{pn 13}%
\special{pa 3136 612}%
\special{pa 3142 2158}%
\special{fp}%
%
\special{pn 8}%
\special{ar 694 958 350 96  3.3333559 5.1443363}%
%
\special{pn 8}%
\special{ar 698 1396 350 96  3.3435871 5.1394229}%
%
\special{pn 8}%
\special{pa 352 1438}%
\special{pa 352 1438}%
\special{fp}%
%
\special{pn 8}%
\special{pa 352 1438}%
\special{pa 352 852}%
\special{fp}%
%
\special{pn 8}%
\special{pa 844 1406}%
\special{pa 844 820}%
\special{fp}%
%
\special{pn 8}%
\special{ar 930 1576 350 96  3.3435871 5.1382089}%
%
\special{pn 8}%
\special{ar 936 2012 348 94  3.3450857 5.1373979}%
%
\special{pn 8}%
\special{pa 690 2040}%
\special{pa 690 2040}%
\special{fp}%
%
\special{pn 8}%
\special{pa 588 2056}%
\special{pa 588 1468}%
\special{fp}%
%
\special{pn 8}%
\special{pa 1078 2024}%
\special{pa 1078 1438}%
\special{fp}%
%
\special{pn 8}%
\special{ar 1276 2184 350 94  3.3456106 5.1375382}%
%
\special{pn 8}%
\special{ar 1282 2622 350 94  3.3456106 5.1372566}%
%
\special{pn 8}%
\special{pa 1234 2676}%
\special{pa 1234 2676}%
\special{fp}%
%
\special{pn 8}%
\special{pa 936 2664}%
\special{pa 936 2078}%
\special{fp}%
%
\special{pn 8}%
\special{pa 1426 2634}%
\special{pa 1426 2046}%
\special{fp}%
%
\special{pn 8}%
\special{ar 1906 2376 350 96  3.3333559 5.1415927}%
%
\special{pn 8}%
\special{ar 1914 2848 356 98  3.3475093 5.1366305}%
%
\special{pn 8}%
\special{pa 1566 2856}%
\special{pa 1566 2268}%
\special{fp}%
%
\special{pn 8}%
\special{pa 2056 2824}%
\special{pa 2056 2238}%
\special{fp}%
\put(5.3200,-10.4800){\makebox(0,0)[lt]{$M_1$}}%
\put(7.4600,-16.5900){\makebox(0,0)[lt]{$M_2$}}%
\put(10.9900,-22.6200){\makebox(0,0)[lt]{$M_3$}}%
\put(17.1600,-24.4800){\makebox(0,0)[lt]{$M_4$}}%
\put(13.1000,-14.6000){\makebox(0,0)[lt]{$\displaystyle{M=\mathop{\oplus}_{i=1}^4M_i\,\,\,\,\mathop{\hookrightarrow}^f}$}}%
%
\special{pn 4}%
\special{pa 3172 2156}%
\special{pa 3148 2114}%
\special{fp}%
\special{pa 3198 2148}%
\special{pa 3148 2058}%
\special{fp}%
\special{pa 3222 2138}%
\special{pa 3148 2000}%
\special{fp}%
\special{pa 3246 2130}%
\special{pa 3148 1944}%
\special{fp}%
\special{pa 3272 2122}%
\special{pa 3146 1888}%
\special{fp}%
\special{pa 3296 2112}%
\special{pa 3146 1830}%
\special{fp}%
\special{pa 3324 2104}%
\special{pa 3146 1774}%
\special{fp}%
\special{pa 3348 2094}%
\special{pa 3144 1718}%
\special{fp}%
\special{pa 3374 2086}%
\special{pa 3144 1662}%
\special{fp}%
\special{pa 3398 2076}%
\special{pa 3144 1606}%
\special{fp}%
\special{pa 3424 2068}%
\special{pa 3144 1548}%
\special{fp}%
\special{pa 3448 2060}%
\special{pa 3144 1492}%
\special{fp}%
\special{pa 3474 2050}%
\special{pa 3144 1436}%
\special{fp}%
\special{pa 3500 2044}%
\special{pa 3142 1380}%
\special{fp}%
\special{pa 3526 2036}%
\special{pa 3142 1322}%
\special{fp}%
\special{pa 3552 2028}%
\special{pa 3142 1266}%
\special{fp}%
\special{pa 3578 2022}%
\special{pa 3142 1210}%
\special{fp}%
\special{pa 3604 2014}%
\special{pa 3142 1152}%
\special{fp}%
\special{pa 3630 2008}%
\special{pa 3140 1096}%
\special{fp}%
\special{pa 3656 2000}%
\special{pa 3140 1040}%
\special{fp}%
\special{pa 3682 1994}%
\special{pa 3140 984}%
\special{fp}%
\special{pa 3708 1986}%
\special{pa 3140 926}%
\special{fp}%
\special{pa 3734 1980}%
\special{pa 3138 870}%
\special{fp}%
\special{pa 3762 1976}%
\special{pa 3138 814}%
\special{fp}%
\special{pa 3788 1970}%
\special{pa 3138 758}%
\special{fp}%
\special{pa 3816 1964}%
\special{pa 3136 702}%
\special{fp}%
\special{pa 3842 1958}%
\special{pa 3136 644}%
\special{fp}%
\special{pa 3870 1954}%
\special{pa 3142 600}%
\special{fp}%
\special{pa 3896 1948}%
\special{pa 3170 596}%
\special{fp}%
\special{pa 3924 1942}%
\special{pa 3198 592}%
\special{fp}%
%
\special{pn 4}%
\special{pa 3952 1938}%
\special{pa 3226 588}%
\special{fp}%
\special{pa 3978 1932}%
\special{pa 3254 584}%
\special{fp}%
\special{pa 4006 1928}%
\special{pa 3282 578}%
\special{fp}%
\special{pa 4036 1928}%
\special{pa 3308 576}%
\special{fp}%
\special{pa 4064 1926}%
\special{pa 3336 570}%
\special{fp}%
\special{pa 4094 1926}%
\special{pa 3364 568}%
\special{fp}%
\special{pa 4122 1926}%
\special{pa 3392 564}%
\special{fp}%
\special{pa 4152 1924}%
\special{pa 3418 558}%
\special{fp}%
\special{pa 4182 1924}%
\special{pa 3446 554}%
\special{fp}%
\special{pa 4212 1922}%
\special{pa 3474 550}%
\special{fp}%
\special{pa 4240 1922}%
\special{pa 3502 546}%
\special{fp}%
\special{pa 4270 1922}%
\special{pa 3530 542}%
\special{fp}%
\special{pa 4300 1920}%
\special{pa 3558 538}%
\special{fp}%
\special{pa 4330 1920}%
\special{pa 3586 536}%
\special{fp}%
\special{pa 4358 1920}%
\special{pa 3614 532}%
\special{fp}%
\special{pa 4388 1918}%
\special{pa 3642 530}%
\special{fp}%
\special{pa 4418 1918}%
\special{pa 3670 526}%
\special{fp}%
\special{pa 4444 1910}%
\special{pa 3698 522}%
\special{fp}%
\special{pa 4444 1856}%
\special{pa 3726 520}%
\special{fp}%
\special{pa 4444 1800}%
\special{pa 3754 516}%
\special{fp}%
\special{pa 4444 1746}%
\special{pa 3782 512}%
\special{fp}%
\special{pa 4446 1690}%
\special{pa 3810 510}%
\special{fp}%
\special{pa 4446 1636}%
\special{pa 3838 506}%
\special{fp}%
\special{pa 4446 1580}%
\special{pa 3864 500}%
\special{fp}%
\special{pa 4446 1524}%
\special{pa 3892 494}%
\special{fp}%
\special{pa 4446 1470}%
\special{pa 3918 488}%
\special{fp}%
\special{pa 4446 1414}%
\special{pa 3944 482}%
\special{fp}%
\special{pa 4446 1360}%
\special{pa 3972 476}%
\special{fp}%
\special{pa 4446 1304}%
\special{pa 3998 470}%
\special{fp}%
\special{pa 4446 1248}%
\special{pa 4026 464}%
\special{fp}%
%
\special{pn 4}%
\special{pa 4446 1194}%
\special{pa 4054 462}%
\special{fp}%
\special{pa 4446 1138}%
\special{pa 4082 460}%
\special{fp}%
\special{pa 4448 1082}%
\special{pa 4110 458}%
\special{fp}%
\special{pa 4448 1028}%
\special{pa 4140 454}%
\special{fp}%
\special{pa 4448 972}%
\special{pa 4168 452}%
\special{fp}%
\special{pa 4448 918}%
\special{pa 4196 450}%
\special{fp}%
\special{pa 4448 862}%
\special{pa 4226 446}%
\special{fp}%
\special{pa 4448 808}%
\special{pa 4254 444}%
\special{fp}%
\special{pa 4448 752}%
\special{pa 4282 442}%
\special{fp}%
\special{pa 4448 696}%
\special{pa 4310 440}%
\special{fp}%
\special{pa 4448 642}%
\special{pa 4338 438}%
\special{fp}%
\special{pa 4448 586}%
\special{pa 4368 434}%
\special{fp}%
\special{pa 4448 532}%
\special{pa 4396 432}%
\special{fp}%
\special{pa 4450 476}%
\special{pa 4424 430}%
\special{fp}%
%
\special{pn 4}%
\special{pa 364 1422}%
\special{pa 352 1400}%
\special{fp}%
\special{pa 388 1414}%
\special{pa 352 1346}%
\special{fp}%
\special{pa 414 1404}%
\special{pa 352 1290}%
\special{fp}%
\special{pa 440 1398}%
\special{pa 352 1234}%
\special{fp}%
\special{pa 468 1396}%
\special{pa 352 1178}%
\special{fp}%
\special{pa 496 1392}%
\special{pa 352 1122}%
\special{fp}%
\special{pa 524 1390}%
\special{pa 352 1068}%
\special{fp}%
\special{pa 554 1386}%
\special{pa 352 1010}%
\special{fp}%
\special{pa 582 1382}%
\special{pa 352 956}%
\special{fp}%
\special{pa 610 1380}%
\special{pa 352 900}%
\special{fp}%
\special{pa 638 1378}%
\special{pa 358 856}%
\special{fp}%
\special{pa 670 1380}%
\special{pa 386 854}%
\special{fp}%
\special{pa 700 1382}%
\special{pa 414 850}%
\special{fp}%
\special{pa 734 1386}%
\special{pa 442 846}%
\special{fp}%
\special{pa 766 1390}%
\special{pa 470 844}%
\special{fp}%
\special{pa 796 1396}%
\special{pa 498 840}%
\special{fp}%
\special{pa 828 1398}%
\special{pa 528 838}%
\special{fp}%
\special{pa 850 1382}%
\special{pa 556 834}%
\special{fp}%
\special{pa 850 1326}%
\special{pa 584 830}%
\special{fp}%
\special{pa 850 1272}%
\special{pa 612 828}%
\special{fp}%
\special{pa 850 1216}%
\special{pa 642 826}%
\special{fp}%
\special{pa 850 1160}%
\special{pa 670 824}%
\special{fp}%
\special{pa 850 1104}%
\special{pa 698 822}%
\special{fp}%
\special{pa 850 1048}%
\special{pa 728 822}%
\special{fp}%
\special{pa 850 994}%
\special{pa 760 824}%
\special{fp}%
\special{pa 850 938}%
\special{pa 790 828}%
\special{fp}%
\special{pa 850 882}%
\special{pa 822 830}%
\special{fp}%
\put(49.6000,-10.2000){\makebox(0,0)[lt]{$V$}}%
%
\special{pn 8}%
\special{pa 3918 178}%
\special{pa 3614 366}%
\special{dt 0.045}%
\special{sh 1}%
\special{pa 3614 366}%
\special{pa 3682 348}%
\special{pa 3660 338}%
\special{pa 3660 314}%
\special{pa 3614 366}%
\special{fp}%
\put(39.5000,-2.1000){\makebox(0,0)[lb]{${\rm Fix}\,G$}}%
%
\special{pn 8}%
\special{ar 3674 3532 848 434  4.8924697 6.2831853}%
%
\special{pn 8}%
\special{pa 3650 3530}%
\special{pa 4526 3530}%
\special{fp}%
%
\special{pn 8}%
\special{pa 3664 3532}%
\special{pa 3822 3100}%
\special{fp}%
%
\special{pn 13}%
\special{ar 4614 3710 1148 338  3.7016562 4.5879067}%
%
\special{pn 13}%
\special{ar 3348 2902 820 678  0.3893365 1.2012596}%
%
\special{pn 8}%
\special{pa 3760 3270}%
\special{pa 3798 3168}%
\special{fp}%
\special{sh 1}%
\special{pa 3798 3168}%
\special{pa 3756 3224}%
\special{pa 3780 3218}%
\special{pa 3794 3236}%
\special{pa 3798 3168}%
\special{fp}%
%
\special{pn 8}%
\special{pa 4320 3536}%
\special{pa 4424 3536}%
\special{fp}%
\special{sh 1}%
\special{pa 4424 3536}%
\special{pa 4356 3516}%
\special{pa 4370 3536}%
\special{pa 4356 3556}%
\special{pa 4424 3536}%
\special{fp}%
\put(26.4000,-32.4000){\makebox(0,0)[lt]{$\displaystyle{\overline M\qquad\,\,\mathop{\hookrightarrow}^{\overline f}}$}}%
\put(39.6900,-37.1000){\makebox(0,0)[lt]{$V/G$}}%
%
\special{pn 8}%
\special{pa 4036 2470}%
\special{pa 4036 2846}%
\special{fp}%
\special{sh 1}%
\special{pa 4036 2846}%
\special{pa 4056 2780}%
\special{pa 4036 2794}%
\special{pa 4016 2780}%
\special{pa 4036 2846}%
\special{fp}%
\put(41.8500,-25.7400){\makebox(0,0)[lt]{$\phi$}}%
%
\special{pn 13}%
\special{pa 3640 1670}%
\special{pa 3824 1826}%
\special{fp}%
\special{sh 1}%
\special{pa 3824 1826}%
\special{pa 3786 1768}%
\special{pa 3782 1792}%
\special{pa 3760 1798}%
\special{pa 3824 1826}%
\special{fp}%
%
\special{pn 13}%
\special{pa 3360 1750}%
\special{pa 3486 1892}%
\special{fp}%
\special{sh 1}%
\special{pa 3486 1892}%
\special{pa 3456 1828}%
\special{pa 3450 1852}%
\special{pa 3426 1854}%
\special{pa 3486 1892}%
\special{fp}%
%
\special{pn 13}%
\special{pa 4030 1580}%
\special{pa 4058 1770}%
\special{fp}%
\special{sh 1}%
\special{pa 4058 1770}%
\special{pa 4068 1700}%
\special{pa 4050 1716}%
\special{pa 4028 1706}%
\special{pa 4058 1770}%
\special{fp}%
%
\special{pn 13}%
\special{pa 4250 1550}%
\special{pa 4278 1748}%
\special{fp}%
\special{sh 1}%
\special{pa 4278 1748}%
\special{pa 4288 1678}%
\special{pa 4270 1694}%
\special{pa 4250 1684}%
\special{pa 4278 1748}%
\special{fp}%
%
\special{pn 13}%
\special{pa 3880 3470}%
\special{pa 3928 3594}%
\special{fp}%
\special{sh 1}%
\special{pa 3928 3594}%
\special{pa 3922 3524}%
\special{pa 3908 3544}%
\special{pa 3886 3538}%
\special{pa 3928 3594}%
\special{fp}%
%
\special{pn 13}%
\special{pa 4190 3400}%
\special{pa 4200 3560}%
\special{fp}%
\special{sh 1}%
\special{pa 4200 3560}%
\special{pa 4216 3492}%
\special{pa 4198 3508}%
\special{pa 4176 3496}%
\special{pa 4200 3560}%
\special{fp}%
%
\special{pn 13}%
\special{pa 3840 3430}%
\special{pa 3764 3330}%
\special{fp}%
\special{sh 1}%
\special{pa 3764 3330}%
\special{pa 3788 3394}%
\special{pa 3796 3372}%
\special{pa 3820 3370}%
\special{pa 3764 3330}%
\special{fp}%
%
\special{pn 13}%
\special{pa 4060 3240}%
\special{pa 3934 3182}%
\special{fp}%
\special{sh 1}%
\special{pa 3934 3182}%
\special{pa 3986 3228}%
\special{pa 3982 3204}%
\special{pa 4004 3192}%
\special{pa 3934 3182}%
\special{fp}%
%
\special{pn 13}%
\special{pa 3644 3522}%
\special{pa 4036 3372}%
\special{fp}%
\special{sh 1}%
\special{pa 4036 3372}%
\special{pa 3966 3378}%
\special{pa 3986 3392}%
\special{pa 3980 3416}%
\special{pa 4036 3372}%
\special{fp}%
%
\special{pn 13}%
\special{pa 3380 990}%
\special{pa 3506 1132}%
\special{fp}%
\special{sh 1}%
\special{pa 3506 1132}%
\special{pa 3476 1070}%
\special{pa 3470 1092}%
\special{pa 3446 1096}%
\special{pa 3506 1132}%
\special{fp}%
%
\special{pn 13}%
\special{pa 3616 910}%
\special{pa 3762 1100}%
\special{fp}%
\special{sh 1}%
\special{pa 3762 1100}%
\special{pa 3736 1036}%
\special{pa 3730 1058}%
\special{pa 3706 1060}%
\special{pa 3762 1100}%
\special{fp}%
%
\special{pn 13}%
\special{pa 4010 840}%
\special{pa 4038 1030}%
\special{fp}%
\special{sh 1}%
\special{pa 4038 1030}%
\special{pa 4048 960}%
\special{pa 4030 976}%
\special{pa 4008 966}%
\special{pa 4038 1030}%
\special{fp}%
%
\special{pn 13}%
\special{pa 4280 800}%
\special{pa 4308 998}%
\special{fp}%
\special{sh 1}%
\special{pa 4308 998}%
\special{pa 4318 928}%
\special{pa 4300 944}%
\special{pa 4280 934}%
\special{pa 4308 998}%
\special{fp}%
%
\special{pn 13}%
\special{pa 4950 1630}%
\special{pa 5100 1810}%
\special{fp}%
\special{sh 1}%
\special{pa 5100 1810}%
\special{pa 5072 1746}%
\special{pa 5066 1768}%
\special{pa 5042 1772}%
\special{pa 5100 1810}%
\special{fp}%
\put(51.5500,-16.7700){\makebox(0,0)[lt]{'s$\,\,:\,\,H$}}%
%
\special{pn 13}%
\special{pa 4772 3494}%
\special{pa 4920 3674}%
\special{fp}%
\special{sh 1}%
\special{pa 4920 3674}%
\special{pa 4894 3610}%
\special{pa 4886 3632}%
\special{pa 4862 3636}%
\special{pa 4920 3674}%
\special{fp}%
\put(49.7600,-35.4100){\makebox(0,0)[lt]{'s$\,\,:\,\,\overline H$}}%
%
\special{pn 13}%
\special{pa 3610 890}%
\special{pa 3586 870}%
\special{pa 3560 850}%
\special{pa 3538 828}%
\special{pa 3520 802}%
\special{pa 3506 774}%
\special{pa 3498 742}%
\special{pa 3504 712}%
\special{pa 3524 686}%
\special{pa 3548 668}%
\special{pa 3576 652}%
\special{pa 3608 642}%
\special{pa 3638 636}%
\special{pa 3670 632}%
\special{pa 3702 628}%
\special{pa 3704 628}%
\special{sp}%
%
\special{pn 13}%
\special{pa 3620 910}%
\special{pa 3642 934}%
\special{pa 3660 960}%
\special{pa 3674 990}%
\special{pa 3680 1020}%
\special{pa 3680 1052}%
\special{pa 3670 1084}%
\special{pa 3654 1110}%
\special{pa 3632 1134}%
\special{pa 3608 1152}%
\special{sp}%
%
\special{pn 13}%
\special{ar 3940 1520 470 184  2.3281168 4.0291079}%
%
\special{pn 13}%
\special{pa 3630 1670}%
\special{pa 3652 1694}%
\special{pa 3670 1720}%
\special{pa 3684 1748}%
\special{pa 3692 1780}%
\special{pa 3692 1812}%
\special{pa 3684 1842}%
\special{pa 3668 1870}%
\special{pa 3648 1894}%
\special{pa 3624 1912}%
\special{sp}%
%
\special{pn 8}%
\special{pa 4270 3536}%
\special{pa 4374 3536}%
\special{fp}%
\special{sh 1}%
\special{pa 4374 3536}%
\special{pa 4306 3516}%
\special{pa 4320 3536}%
\special{pa 4306 3556}%
\special{pa 4374 3536}%
\special{fp}%
%
\special{pn 8}%
\special{pa 3740 3320}%
\special{pa 3778 3218}%
\special{fp}%
\special{sh 1}%
\special{pa 3778 3218}%
\special{pa 3736 3274}%
\special{pa 3760 3268}%
\special{pa 3774 3286}%
\special{pa 3778 3218}%
\special{fp}%
%
\special{pn 13}%
\special{pa 4040 3420}%
\special{pa 4060 3560}%
\special{fp}%
\special{sh 1}%
\special{pa 4060 3560}%
\special{pa 4070 3492}%
\special{pa 4052 3508}%
\special{pa 4032 3498}%
\special{pa 4060 3560}%
\special{fp}%
%
\special{pn 13}%
\special{pa 3950 3360}%
\special{pa 3846 3280}%
\special{fp}%
\special{sh 1}%
\special{pa 3846 3280}%
\special{pa 3888 3336}%
\special{pa 3888 3314}%
\special{pa 3912 3306}%
\special{pa 3846 3280}%
\special{fp}%
%
\special{pn 13}%
\special{ar 2720 3870 1780 576  3.5195978 4.5261219}%
\end{picture}%
\hspace{1.5truecm}}

\vspace{0.5truecm}

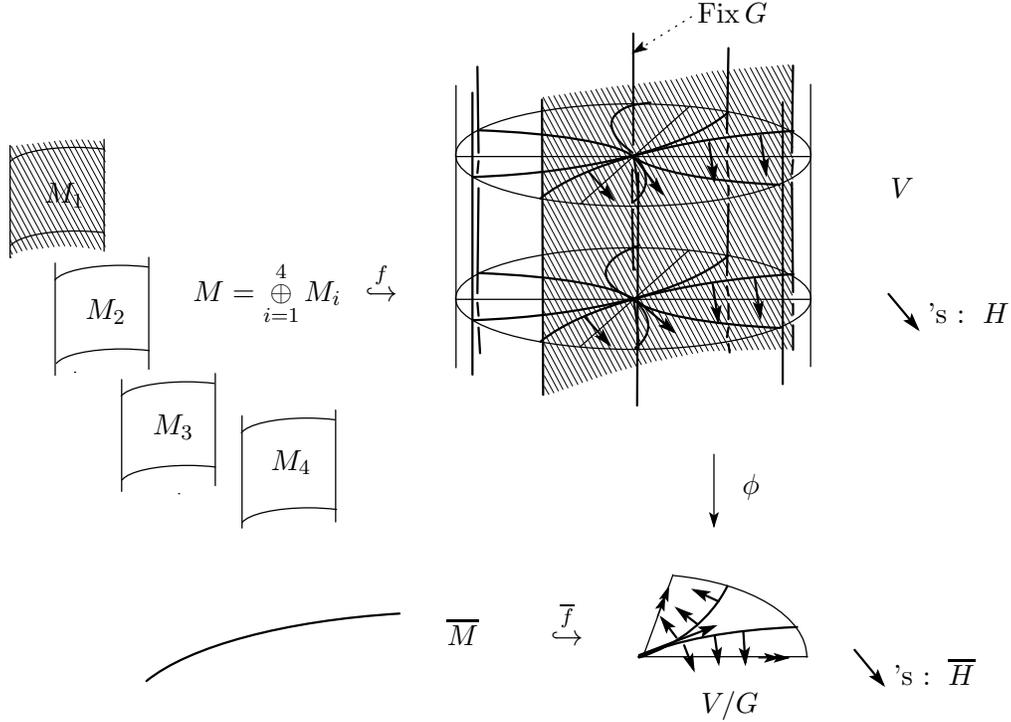
\centerline{{\bf Figure 3.}}

\vspace{0.4truecm}

\noindent
Assume that the codimension of $M$ is equal to one.  
Denote by $\widetilde{{\cal H}}$ (resp. $\widetilde{{\cal V}}$) the horizontal 
(resp. vertical) distribution of $\phi$.  
Denote by ${\rm pr}_{\widetilde{\cal H}}$ (resp. 
${\rm pr}_{\widetilde{\cal V}}$) the orthogonal projection of $TV$ onto 
$\widetilde{\cal H}$ (resp. $\widetilde{\cal V}$).  For simplicity, for 
$X\in TV$, we denote ${\rm pr}_{\widetilde{\cal H}}(X)$ (resp. 
${\rm pr}_{\widetilde{\cal V}}(X)$) by $X_{\widetilde{\cal H}}$ 
(resp. $X_{\widetilde{\cal V}}$).  
Define a distribution ${\cal H}_t$ on 
$M$ by $f_{t\ast}(({\cal H}_t)_u)=f_{t\ast}(T_uM)\cap
\widetilde{{\cal H}}_{f_t(u)}$ ($u\in M$) and a distribution ${\cal V}_t$ on 
$M$ by $f_{t\ast}(({\cal V}_t)_u)=\widetilde{{\cal V}}_{f_t(u)}$ ($u\in M$).  
Note that ${\cal V}_t$ is independent of the choice of $t\in[0,T)$.  
Denote by $g_t,h_t,A_t,H_t$ and $\xi_t$ 
the induced metric, the second fundamental form, the shape tensor and the 
regularized mean curvature vector and the unit normal vector field of $f_t$, 
respectively.  
The group $G$ acts on $M$ through $f_t$.  
Since $\phi:V\to V/G$ is a $G$-orbibundle and $\widetilde{\cal H}$ is a connection 
of the orbibundle, it follows from Proposition 4.1 that this action $G\curvearrowright M$ 
is independent of the choice of $t\in[0,T)$.  
It is clear that quantities $g_t,h_t,A_t$ and $H_t$ are $G$-invariant.  
Also, let $\nabla^t$ 
be the Riemannian connection of $g_t$.  
Let $\pi_M$ be the projection of $M\times[0,T)$ onto $M$.  
For a vector bundle $E$ over $M$, denote by $\pi_M^{\ast}E$ the induced bundle 
of $E$ by $\pi_M$.  Also denote by $\Gamma(E)$ the space of all sections of 
$E$.  Define a section $g$ of $\pi_M^{\ast}(T^{(0,2)}M)$ by 
$g(u,t)=(g_t)_u$ ($(u,t)\in M\times[0,T)$), where $T^{(0,2)}M$ is the 
$(0,2)$-tensor bundle of $M$.  Similarly, we define a section 
$h$ of $\pi_M^{\ast}(T^{(0,2)}M)$, a section $A$ of $\pi_M^{\ast}
(T^{(1,1)}M)$, a map $H:M\times[0,T)\to TV$ and 
a map $\xi:M\times[0,T)\to TV$.  
We regard $H$ and $\xi$ as $V$-valued functions over 
$M\times[0,T)$ under the identification of $T_uV$'s ($u\in V$) and $V$.  
Define a subbundle ${\cal H}$ (resp. ${\cal V}$) of 
$\pi_M^{\ast}TM$ by ${\cal H}_{(u,t)}:=({\cal H}_t)_u$ (resp. 
${\cal V}_{(u,t)}:=({\cal V}_t)_u$).  Denote by ${\rm pr}_{\cal H}$ (resp. 
${\rm pr}_{\cal V}$) the orthogonal projection of $\pi_M^{\ast}(TM)$ onto 
${\cal H}$ (resp. ${\cal V}$).  For simplicity, for $X\in \pi_M^{\ast}(TM)$, 
we denote ${\rm pr}_{\cal H}(X)$ (resp. ${\rm pr}_{\cal V}(X)$) by 
$X_{\cal H}$ (resp. $X_{\cal V}$).  
The bundle $\pi_M^{\ast}(TM)$ 
is regarded as a subbundle of $T(M\times[0,T))$.  
For a section $B$ of $\pi_M^{\ast}(T^{(r,s)}M)$, 
we define $\displaystyle{\frac{\partial B}{\partial t}}$ by 
$\displaystyle{\left(\frac{\partial B}{\partial t}\right)_{(u,t)}
:=\frac{d B_{(u,t)}}{dt}}$, where the right-hand side of this relation is 
the derivative of the vector-valued function $t\mapsto B_{(u,t)}\,(\in 
T^{(r,s)}_uM)$.  Also, we define a section $B_{\cal H}$ of 
$\pi_M^{\ast}(T^{(r,s)}M)$ by 
$$B_{\cal H}=\mathop{({\rm pr}_{\cal H}\otimes\cdots
\otimes{\rm pr}_{\cal H})}_{(r{\rm -times})}\circ B\circ
\mathop{({\rm pr}_{\cal H}\otimes\cdots
\otimes{\rm pr}_{\cal H})}_{(s{\rm -times})}.$$
The restriction of $B_{\cal H}$ to 
${\cal H}\times\cdots\times{\cal H}$ ($s$-times) is regarded as a section of 
the $(r,s)$-tensor bundle ${\cal H}^{(r,s)}$ of $\cal H$.  This restriction 
also is denoted by the same symbol $B_{\cal H}$.  
For a tangent vector field $X$ on $M$ (or an open set $U$ of $M$), we define 
a section $\bar X$ of 
$\pi_M^{\ast}TM$ (or $\pi_M^{\ast}TM\vert_U$) by $\bar X_{(u,t)}:=X_u$ 
($(u,t)\in M\times[0,T)$).  Denote by $\widetilde{\nabla}$ 
the Riemannian connection of $V$.  
Define a connection $\nabla$ of $\pi_M^{\ast}TM$ by 
$$(\nabla_XY)_{(\cdot,t)}:=\nabla^t_XY_{(\cdot,t)}\,\,\,{\rm and}\,\,\,
\nabla_{\frac{\partial}{\partial t}}Y:=\frac{dY_{(u,\cdot)}}{dt}$$
for $X\in T_{(u,t)}(M\times\{t\})$ and $Y\in\Gamma(\pi_M^{\ast}TM)$, where 
$\displaystyle{\frac{dY_{(u,t)}}{dt}}$ is the derivative of the 
vector-valued function $t\mapsto Y_{(u,t)}\,(\in T_uM)$.  
Define a connection $\nabla^{\cal H}$ of ${\cal H}$ by 
$\nabla^{\cal H}_XY:=(\nabla_XY)_{\cal H}$ for $X\in T(M\times[0,T))$ and 
$Y\in\Gamma({\cal H})$.  Similarly, define a connection 
$\nabla^{\cal V}$ of ${\cal V}$ by $\nabla^{\cal V}_XY:=(\nabla_XY)_{\cal V}$ 
for $X\in T(M\times[0,T))$ and $Y\in\Gamma({\cal V})$.  
Now we shall derive the evolution equations for some geometric quantities.  
First we derive the following evolution equation for $g_{\cal H}$.  

\vspace{0.5truecm}

\noindent
{\bf Lemma 4.2.} {\sl The sections $(g_{\cal H})_t$'s of $\pi_M^{\ast}(T^{(0,2)}M)$ 
satisfy the following evolution equation:
$$\frac{\partial g_{\cal H}}{\partial t}
=-2\vert\vert H\vert\vert h_{\cal H},$$
where $\vert\vert H\vert\vert:=\sqrt{g(H,H)}$.}

\vspace{0.5truecm}

\noindent
{\it Proof.} Take $X,Y\in \Gamma(TM)$.  We have 
$$\begin{array}{l}
\hspace{0.6truecm}
\displaystyle{
\frac{\partial g_{\cal H}}{\partial t}
(\bar X,\bar Y)
=\frac{\partial}{\partial t}g_{\cal H}(\bar X,\bar Y)
=\frac{\partial}{\partial t}g({\bar X}_{\cal H},{\bar Y}_{\cal H})
=\frac{\partial}{\partial t}\langle F_{\ast}{\bar X}_{\cal H},
F_{\ast}{\bar Y}_{\cal H}\rangle}\\
\displaystyle{=\langle\frac{\partial}{\partial t}({\bar X}_{\cal H}F),
{\bar Y}_{\cal H}F\rangle
+\langle{\bar X}_{\cal H}F,\frac{\partial}{\partial t}({\bar Y}_{\cal H}F)
\rangle}\\
\displaystyle{=\langle{\bar X}_{\cal H}
\left(\frac{\partial F}{\partial t}\right)
+[\frac{\partial}{\partial t},{\bar X}_{\cal H}]F,{\bar Y}_{\cal H}F\rangle
+\langle{\bar X}_{\cal H}F,
{\bar Y}_{\cal H}\left(\frac{\partial F}{\partial t}\right)
+[\frac{\partial}{\partial t},{\bar Y}_{\cal H}]F\rangle}\\
\displaystyle{=\langle{\bar X}_{\cal H}(\vert\vert H\vert\vert\xi),
{\bar Y}_{\cal H}F\rangle
+\langle{\bar X}_{\cal H}F,{\bar Y}_{\cal H}(\vert\vert H\vert\vert\xi)
\rangle}\\
\displaystyle{=-\vert\vert H\vert\vert g(A{\bar X}_{\cal H},{\bar Y}_{\cal H})
-\vert\vert H\vert\vert g({\bar X}_{\cal H},A{\bar Y}_{\cal H})
=-2\vert\vert H\vert\vert h_{\cal H}(\bar X,\bar Y),}
\end{array}$$
where we use 
$\displaystyle{\left[\frac{\partial}{\partial t},{\bar X}_{\cal H}\right]\in{\cal V}}$ 
and 
$\displaystyle{\left[\frac{\partial}{\partial t},{\bar Y}_{\cal H}\right]\in{\cal V}}$.  
Thus we obtain the desired evolution equation.  
\hspace{11.4truecm}q.e.d.

\vspace{0.5truecm}

Next we derive the following evolution equation for $\xi$.  

\vspace{0.5truecm}

\noindent
{\bf Lemma 4.3.} {\sl The unit normal vector fields $\xi_t$'s 
satisfy the following evolution equation:
$$\frac{\partial\xi}{\partial t}=-F_{\ast}
({\rm grad}_g\vert\vert H\vert\vert),$$
where ${\rm grad}_g(\vert\vert H\vert\vert)$ 
is the element of $\pi_M^{\ast}(TM)$ such that 
$d\vert\vert H\vert\vert(X)=g({\rm grad}_g\vert\vert H\vert\vert,X)$ for 
any $X\in\pi_M^{\ast}(TM)$.}

\vspace{0.5truecm}

\noindent
{\it Proof.} Since $\langle\xi,\xi\rangle=1$, we have 
$\langle\frac{\partial\xi}{\partial t},\xi\rangle=0$.  Hence 
$\frac{\partial\xi}{\partial t}$ is tangent to 
$f_t(M)$.  Take any $(u_0,t_0)\in M\times[0,T)$.  
Let $\{e_i\}_{i=1}^{\infty}$ be an orthonormal base of $T_{u_0}M$ 
with respect to $g_{(u_0,t_0)}$.  By the Fourier expanding 
$\frac{\partial\xi}{\partial t}\left\vert_{t=t_0}\right.$, we have 
$$\begin{array}{l}
\hspace{0.6truecm}\displaystyle{\left.\frac{\partial\xi}{\partial t}
\right\vert_{t=t_0}
=\sum\limits\langle\left.\frac{\partial\xi}{\partial t}\right\vert_{t=t_0},
f_{t_0\ast}(\bar e_i\vert_{t=t_0})\rangle
f_{t_0\ast}(\bar e_i\vert_{t=t_0})}\\
\displaystyle{=-\sum\limits\langle\xi_{t_0},
\left.\frac{\partial f_{t\ast}(\bar e_i)}{\partial t}\right\vert_{t=t_0}
\rangle f_{t_0\ast}(\bar e_i\vert_{t=t_0})}\\
\displaystyle{=-\sum\limits\langle\xi_{t_0},
\left.\frac{\partial}{\partial t}(\bar e_i F)\right\vert_{t=t_0}\rangle
f_{t_0\ast}(\bar e_i\vert_{t=t_0})}\\
\displaystyle{=-\sum\limits\langle\xi_{t_0},
\bar e_i\left(\left.\frac{\partial F}{\partial t}\right\vert_{t=t_0}\right)
\rangle f_{t_0\ast}(\bar e_i\vert_{t=t_0})}\\
\displaystyle{=-\sum\limits\langle\xi_{t_0},
(\bar e_iH)\vert_{t=t_0}\rangle f_{t_0\ast}(\bar e_i\vert_{t=t_0})}\\
\displaystyle{=-\sum\limits(\bar e_i\vert\vert H\vert\vert)
\vert_{t=t_0}f_{t_0\ast}(\bar e_i\vert_{t=t_0})}\\
\displaystyle{=-\sum\limits g_{t_0}({\rm grad}_{g_{t_0}}\vert\vert H_{t_0}
\vert\vert,\bar e_i\vert_{t=t_0})f_{t_0\ast}(\bar e_i\vert_{t=t_0})}\\
\displaystyle{=-f_{t_0\ast}({\rm grad}_{g_{t_0}}\vert\vert H_{t_0}\vert\vert)
=-F_{\ast}({\rm grad}_g\vert\vert H\vert\vert)\vert_{t=t_0}}
\end{array}$$
on $U$, where we use $\displaystyle{\left[\frac{\partial}{\partial t},\bar e_i\right]=0}$.  
Here we note that $\sum(\cdot)_i$ means 
$\lim\limits_{k\to\infty}\sum_{i\in S_k}(\cdot)_i$ as 
$S_k:=\{i\,\,\vert\,\,\vert(\cdot)_i\vert>\frac1k\}$ ($k\in{\Bbb N}$).  
In particular, we obtain 
$$\left(\frac{\partial\xi}{\partial t}\right)_{(u_0,t_0)}
=-(F_{\ast}(({\rm grad}_g\vert\vert H\vert\vert))_{(u_0,t_0)}.$$
This completes the proof.  \hspace{8.6truecm}q.e.d.

\vspace{0.5truecm}

Let $S_t$ ($0\leq t<T$) be a $C^{\infty}$-family of a $(r,s)$-tensor fields on 
$M$ and $S$ a section of $\pi_M^{\ast}(T^{(r,s)}M)$ defined by 
$S_{(u,t)}:=(S_t)_u$.  We define a section $\triangle_{{\cal H}}S$ 
of $\pi_M^{\ast}(T^{(r,s)}M)$ by 
$$(\triangle_{{\cal H}}S)_{(u,t)}:=\sum_{i=1}^n\nabla_{e_i}\nabla_{e_i}S,$$
where $\nabla$ is the connection of $\pi_M^{\ast}(T^{(r,s)}M)$ (or 
$\pi_M^{\ast}(T^{(r,s+1)}M)$) induced from $\nabla$ and $\{e_1,\cdots,e_n\}$ 
is an orthonormal base of ${\cal H}_{(u,t)}$ with respect to 
$(g_{\cal H})_{(u,t)}$.  
Also, we define a section $\bar{\triangle}_{{\cal H}}S_{\cal H}$ 
of ${\cal H}^{(r,s)}$ by 
$$(\bar{\triangle}_{{\cal H}}S_{\cal H})_{(u,t)}:=
\sum_{i=1}^n\nabla^{\cal H}_{e_i}\nabla^{\cal H}_{e_i}S_{\cal H},$$
where $\nabla^{\cal H}$ is the connection of ${\cal H}^{(r,s)}$ (or 
${\cal H}^{(r,s+1)}$) induced from $\nabla^{\cal H}$ and $\{e_1,\cdots,e_n\}$ 
is as above.  Let ${\cal A}^{\phi}$ be the section of $T^{\ast}V\otimes 
T^{\ast}V\otimes TV$ defined by 
$${\cal A}^{\phi}_XY:=(\widetilde{\nabla}_{X_{\widetilde{\cal H}}}
Y_{\widetilde{\cal H}})_{\widetilde{\cal V}}
+(\widetilde{\nabla}_{X_{\widetilde{\cal H}}}Y_{\widetilde{\cal V}})_{\widetilde{\cal H}}
\,\,\,\,(X,Y\in TV).$$
Also, let ${\cal T}^{\phi}$ be the section of 
$T^{\ast}V\otimes T^{\ast}V\otimes TV$ defined by 
$${\cal T}^{\phi}_XY:=(\widetilde{\nabla}_{X_{\widetilde{\cal V}}}
Y_{\widetilde{\cal H}})_{\widetilde{\cal V}}
+(\widetilde{\nabla}_{X_{\widetilde{\cal V}}}
Y_{\widetilde{\cal V}})_{\widetilde{\cal H}}
\,\,\,\,(X,Y\in TV).$$
Also, let ${\cal A}_t$ be the section of $T^{\ast}M\otimes 
T^{\ast}M\otimes TM$ defined by 
$$({\cal A}_t)_XY:=(\nabla^t_{X_{{\cal H}_t}}
Y_{{\cal H}_t})_{{\cal V}_t}+(\nabla^t_{X_{{\cal H}_t}}
Y_{{\cal V}_t})_{{\cal H}_t}\,\,\,\,(X,Y\in TM).$$
Also let ${\cal A}$ be the section of $\pi_M^{\ast}(T^{\ast}M\otimes T^{\ast}M
\otimes TM)$ defined in terms of ${\cal A}_t$'s ($t\in[0,T)$).  
Also, let ${\cal T}_t$ be the section of $T^{\ast}M\otimes 
T^{\ast}M\otimes TM$ defined by 
$$({\cal T}_t)_XY:=(\nabla^t_{X_{{\cal V}_t}}Y_{{\cal V}_t})_{{\cal H}_t}
+(\nabla^t_{X_{{\cal V}_t}}Y_{{\cal H}_t})_{{\cal V}_t}\,\,\,\,(X,Y\in TM).$$
Also let ${\cal T}$ be the section of $\pi_M^{\ast}(T^{\ast}M\otimes T^{\ast}M
\otimes TM)$ defined in terms of ${\cal T}_t$'s ($t\in[0,T)$).  
Clearly we have 
$$F_{\ast}({\cal A}_XY)={\cal A}^{\phi}_{F_{\ast}X}F_{\ast}Y$$
for $X,Y\in{\cal H}$ and 
$$F_{\ast}({\cal T}_WX)={\cal T}^{\phi}_{F_{\ast}W}F_{\ast}X$$
for $X\in{\cal H}$ and $W\in{\cal V}$.  
Let $E$ be a vector bundle over $M$.  For a section $S$ of 
$\pi_M^{\ast}(T^{(0,r)}M\otimes E)$, we define 
$\displaystyle{{\rm Tr}_{g_{\cal H}}^{\bullet}\,S(\cdots,\mathop{\bullet}^j,
\cdots,\mathop{\bullet}^k,\cdots)}$ by 
$$({\rm Tr}_{g_{\cal H}}^{\bullet}\,S(\cdots,\mathop{\bullet}^j,\cdots,
\mathop{\bullet}^k,\cdots))
_{(u,t)}=\sum_{i=1}^nS_{(u,t)}(\cdots,\mathop{e_i}^j,\cdots,\mathop{e_i}^k,
\cdots)$$
$((u,t)\in M\times[0,T))$, where $\{e_1,\cdots,e_n\}$ is an 
orthonormal base of ${\cal H}_{(u,t)}$ with respect to 
$(g_{\cal H})_{(u,t)}$, $\displaystyle{S(\cdots,\mathop{\bullet}^j,
\cdots,\mathop{\bullet}^k,\cdots)}$ means that $\bullet$ is entried into the 
$j$-th component and the $k$-th component of $S$ and 
$\displaystyle{S_{(u,t)}(\cdots,\mathop{e_i}^j,\cdots,\mathop{e_i}^k,\cdots)}$ 
means that $e_i$ is entried into the $j$-th component and the $k$-th component 
of $S_{(u,t)}$.  

Then we have the following relation.  

\vspace{0.5truecm}

\noindent
{\bf Lemma 4.4.} {\sl Let $S$ be a section of 
$\pi_M^{\ast}(T^{(0,2)}M)$ which is symmetric with respect to $g$.  
Then we have 
$$\begin{array}{l}
\displaystyle{(\triangle_{\cal H}S)_{\cal H}(X,Y)=
(\triangle_{\cal H}^{\cal H}S_{\cal H})(X,Y)}\\
\hspace{3.2truecm}
\displaystyle{-2{\rm Tr}^{\bullet}_{g_{\cal H}}
((\nabla_{\bullet}S)({\cal A}_{\bullet}X,Y))
-2{\rm Tr}^{\bullet}_{g_{\cal H}}((\nabla_{\bullet}S)({\cal A}_{\bullet}Y,X))
}\\
\hspace{3.2truecm}
\displaystyle{-{\rm Tr}^{\bullet}_{g_{\cal H}}
S({\cal A}_{\bullet}({\cal A}_{\bullet}X),Y)
-{\rm Tr}^{\bullet}_{g_{\cal H}}S({\cal A}_{\bullet}({\cal A}_{\bullet}Y),X)}\\
\hspace{3.2truecm}
\displaystyle{-{\rm Tr}^{\bullet}_{g_{\cal H}}
S((\nabla_{\bullet}{\cal A})_{\bullet}X,Y)
-{\rm Tr}^{\bullet}_{g_{\cal H}}S((\nabla_{\bullet}{\cal A})_{\bullet}Y,X)}\\
\hspace{3.2truecm}
\displaystyle{-2{\rm Tr}^{\bullet}_{g_{\cal H}}
S({\cal A}_{\bullet}X,{\cal A}_{\bullet}Y)}
\end{array}$$
for $X,Y\in{\cal H}$, where $\nabla$ is the connection of 
$\pi_M^{\ast}(T^{(1,2)}M)$ induced from $\nabla$.}

\vspace{0.5truecm}

\noindent
{\it Proof.} Take any $(u_0,t_0)\in M\times[0,T)$.  
Let $\{e_1,\cdots,e_n\}$ be an orthonormal base of 
${\cal H}_{(u_0,t_0)}$ with respect to 
$(g_{\cal H})_{(u_0,t_0)}$.  Take any $X,Y\in{\cal H}_{(u_0,t_0)}$.  
Let $\widetilde X$ be a section of ${\cal H}$ on a neighborhood of 
$(u_0,t_0)$ with $\widetilde X_{(u_0,t_0)}=X$ and 
$(\nabla^{\cal H} X)_{(u_0,t_0)}=0$.  Similarly we define $\widetilde Y$ and 
$\widetilde e_i$.  Let $W=X,Y$ or $e_i$.  Then, it follows from 
$(\nabla^{\cal H}_{e_i}\widetilde W)_{(u_0,t_0)}=0,\,
(\nabla_{e_i}\widetilde W)_{(u_0,t_0)}={\cal A}_{e_i}W$ and the 
skew-symmetricness of ${\cal A}\vert_{{\cal H}\times{\cal H}}$ that 
$$\begin{array}{l}
\hspace{0.6truecm}\displaystyle{(\triangle_{\cal H}S)_{\cal H}(X,Y)=
\sum_{i=1}^n(\nabla_{e_i}\nabla_{e_i}S)(X,Y)}\\
\displaystyle{=\sum_{i=1}^n(\nabla^{\cal H}_{e_i}\nabla^{\cal H}_{e_i}
S_{\cal H})(X,Y)-2\sum_{i=1}^n\left((\nabla_{e_i}S)({\cal A}_{e_i}X,Y)
+(\nabla_{e_i}S)({\cal A}_{e_i}Y,X)\right)}\\
\hspace{0.6truecm}
\displaystyle{-\sum_{i=1}^n\left(S({\cal A}_{e_i}({\cal A}_{e_i}X),Y)
+S({\cal A}_{e_i}({\cal A}_{e_i}Y),X)\right)}\\
\hspace{0.6truecm}
\displaystyle{-\sum_{i=1}^n\left(S((\nabla_{e_i}{\cal A})_{e_i}X,Y)
+S((\nabla_{e_i}{\cal A})_{e_i}Y,X)\right)}\\
\hspace{0.6truecm}
\displaystyle{-2\sum_{i=1}^nS({\cal A}_{e_i}X,{\cal A}_{e_i}Y).}
\end{array}$$
The right-hand side of this relation is equal to the right-hand side of the relation 
in the statement.  This completes the proof.  
\hspace{6.5truecm}q.e.d.


\vspace{0.5truecm}

Also we have the following Simons-type identity.  

\vspace{0.5truecm}

\noindent
{\bf Lemma 4.5.} {\sl We have 
$$\triangle_{{\cal H}}h=\nabla d\vert\vert H\vert\vert
+\vert\vert H\vert\vert(A^2)_{\sharp}
-({\rm Tr}\,(A^2)_{{\cal H}})h,$$
where $(A^2)_{\sharp}$ is the element of $\Gamma(\pi_M^{\ast}T^{(0,2)}M)$ 
defined by $(A^2)_{\sharp}(X,Y):=g(A^2X,Y)$ ($X,Y\in\pi_M^{\ast}TM$).}

\vspace{0.5truecm}

\noindent
{\it Proof.} Take $X,Y,Z,W\in\pi_M^{\ast}(TM)$.  
Since the ambient space $V$ is flat, it follows from the Ricci's 
identity, the Gauss equation and the Codazzi equation that 
$$\begin{array}{l}
\hspace{0.6truecm}
\displaystyle{(\nabla_X\nabla_Yh)(Z,W)
-(\nabla_Z\nabla_Wh)(X,Y)}\\
\displaystyle{=(\nabla_X\nabla_Zh)(Y,W)
-(\nabla_Z\nabla_Xh)(Y,W)}\\
\displaystyle{=h(X,Y)h(AZ,W)-h(Z,Y)h(AX,W)}\\
\hspace{0.6truecm}
\displaystyle{+h(X,W)h(AZ,Y)-h(Z,W)h(AX,Y).}
\end{array}$$
By using this relation, we obtain the desired relation.  
\hspace{3.75truecm}q.e.d.

\vspace{0.3truecm}

\noindent
{\bf Note.} In the sequel, we omit the notation $F_{\ast}$ 
for simplicity.  

\vspace{0.3truecm}

Define a section ${\cal R}$ of $\pi_M^{\ast}({\cal H}^{(0,2)})$ by 
$$\begin{array}{l}
\displaystyle{{\cal R}(X,Y):=
{\rm Tr}^{\bullet}_{g_{\cal H}}h({\cal A}_{\bullet}({\cal A}_{\bullet}X),Y)
+{\rm Tr}^{\bullet}_{g_{\cal H}}h({\cal A}_{\bullet}({\cal A}_{\bullet}Y),X)}\\
\hspace{2.25truecm}\displaystyle{
+{\rm Tr}^{\bullet}_{g_{\cal H}}h((\nabla_{\bullet}{\cal A})_{\bullet}X,
Y)+{\rm Tr}^{\bullet}_{g_{\cal H}}h((\nabla_{\bullet}{\cal A})_{\bullet}Y,X)}\\
\hspace{2.25truecm}\displaystyle{
+2{\rm Tr}^{\bullet}_{g_{\cal H}}(\nabla_{\bullet}h)({\cal A}_{\bullet}X,Y)
+2{\rm Tr}^{\bullet}_{g_{\cal H}}(\nabla_{\bullet}h)({\cal A}_{\bullet}Y,X)}\\
\hspace{2.25truecm}\displaystyle{
+2{\rm Tr}^{\bullet}_{g_{\cal H}}h({\cal A}_{\bullet}X,{\cal A}_{\bullet}Y)
\qquad\qquad\quad(X,Y\in{\cal H}).}
\end{array}$$
From Lemmas 4.3, 4.4 and 4.5, we derive the following evolution equation for $(h_{\cal H})_t$)s.  

\vspace{0.3truecm}

\noindent
{\bf Theorem 4.6.} {\sl The sections $(h_{\cal H})_t$'s of $\pi_M^{\ast}(T^{(0,2)}M)$ 
satisfies the following evolution equation:
$$\begin{array}{l}
\displaystyle{\frac{\partial h_{\cal H}}{\partial t}(X,Y)
=(\triangle_{\cal H}^{\cal H}h_{\cal H})(X,Y)-2\vert\vert H\vert\vert
((A_{\cal H})^2)_{\sharp}(X,Y)-2\vert\vert H\vert\vert
(({\cal A}^{\phi}_{\xi})^2)_{\sharp}(X,Y)}\\
\hspace{2.4truecm}\displaystyle{
+{\rm Tr}\left((A_{\cal H})^2-(({\cal A}^{\phi}_{\xi})^2)_{\cal H}\right)
h_{\cal H}(X,Y)
-{\cal R}(X,Y)}
\end{array}$$
for $X,Y\in{\cal H}$.}

\vspace{0.3truecm}

\noindent
{\it Proof.} Take $X,Y\in{\cal H}_{(u,t)}$.  
Easily we have 
$$AX=A_{\cal H}X+{\cal A}^{\phi}_{\xi}X,\leqno{(4.1)}$$
and
$$(A^2)_{\cal H}X=(A_{\cal H})^2X-({\cal A}^{\phi}_{\xi})^2X,\leqno{(4.2)}$$
where we use 
$$\left(\widetilde{\nabla}_W\xi\right)_{\widetilde{\cal H}}=
\left(\widetilde{\nabla}_{\xi}W+[W,\xi]\right)_{\widetilde{\cal H}}=
\left(\widetilde{\nabla}_{\xi}W\right)_{\widetilde{\cal H}}={\cal A}_{\xi}W$$
for $W\in\Gamma(\widetilde{\cal V})$ because of $[W,\xi]\in\Gamma
(\widetilde{\cal V})$.  
Also, since $\displaystyle{\left[\frac{\partial}{\partial t},\bar X_{\cal H}
\right]\in{\cal V}}$, we have 
$$\left[\frac{\partial}{\partial t},{\bar X}_{\cal H}\right]
=2\vert\vert H\vert\vert{\cal A}^{\phi}_{\xi}{\bar X}_{\cal H}.\leqno{(4.3)}$$
From Lemma 4.3, $(4.1), (4.2)$ and $(4.3)$, we have 
$$\begin{array}{l}
\hspace{0.6truecm}
\displaystyle{\frac{\partial h_{\cal H}}{\partial t}(X,Y)
=\frac{\partial}{\partial t}(h_{\cal H}(\bar X,\bar Y))
=\frac{\partial}{\partial t}\langle\xi,{\bar X}_{\cal H}({\bar Y}_{\cal H}F)
\rangle}\\
\displaystyle{
=\langle\frac{\partial\xi}{\partial t},{\bar X}_{\cal H}({\bar Y}_{\cal H}F)
\rangle+\langle\xi,\frac{\partial}{\partial t}
\left({\bar X}_{\cal H}({\bar Y}_{\cal H}F)\right)\rangle}\\
\displaystyle{
=-\langle F_{\ast}({\rm grad}_g\vert\vert H\vert\vert),
\widetilde{\nabla}_XF_{\ast}{\bar Y}_{\cal H}\rangle
+\langle\xi,
X\left({\bar Y}_{\cal H}\left(
\frac{\partial F}{\partial t}\right)\right)\rangle}\\
\hspace{0.6truecm}\displaystyle{
+\langle\xi,X([\frac{\partial}{\partial t},{\bar Y}_{\cal H}]F)
\rangle
+\langle\xi,[\frac{\partial}{\partial t},{\bar X}_{\cal H}]({\bar Y}_{\cal H}F)
\rangle}\\
\displaystyle{=-g({\rm grad}_g\vert\vert H\vert\vert,\nabla_X
{\bar Y}_{\cal H})+X({\bar Y}_{\cal H}\vert\vert H\vert\vert)
-\vert\vert H\vert\vert\langle\xi,\widetilde{\nabla}_X
F_{\ast}(A({\bar Y}_{\cal H})\rangle}\\
\hspace{0.6truecm}\displaystyle{
+\langle\xi,\widetilde{\nabla}_XF_{\ast}([\frac{\partial}{\partial t},
{\bar Y}_{\cal H}])\rangle
+\langle\xi,\widetilde{\nabla}_{[\frac{\partial}{\partial t},
{\bar X}_{\cal H}]}F_{\ast}{\bar Y}_{\cal H}\rangle}\\
\displaystyle{=(\nabla d\vert\vert H\vert\vert)(X,Y)
-\vert\vert H\vert\vert h_{\cal H}(X,A_{\cal H}Y)
+\vert\vert H\vert\vert h(X,{\cal A}^{\phi}_{\xi}Y)
+2\vert\vert H\vert\vert h({\cal A}^{\phi}_{\xi}X,Y)}\\
\displaystyle{=(\nabla d\vert\vert H\vert\vert)(X,Y)
-\vert\vert H\vert\vert g_{\cal H}((A_{\cal H})^2X,Y)
-3\vert\vert H\vert\vert g(({\cal A}^{\phi}_{\xi})^2X,Y)}
\end{array}$$
From this relation and the Simons-type identity in Lemma 4.5, we have 
$$\begin{array}{l}
\displaystyle{\frac{\partial h_{\cal H}}{\partial t}
=\triangle_{\cal H}h-2\vert\vert H\vert\vert((A_{\cal H})^2)_{\sharp}
-2\vert\vert H\vert\vert(({\cal A}^{\phi}_{\xi})^2)_{\sharp}}\\
\hspace{1.35truecm}\displaystyle{+{\rm Tr}\left((A_{\cal H})^2
-(({\cal A}^{\phi}_{\xi})^2)_{\cal H}\right)h_{\cal H}.}
\end{array}\leqno{(4.4)}$$
Substituting the relation in Lemma 4.4 into $(4.4)$, 
we obtain the desired relation.  
\begin{flushright}q.e.d.\end{flushright}

\vspace{0.15truecm}

From Lemma 4.2, we derive the following relation.  

\vspace{0.3truecm}

\noindent
{\bf Lemma 4.7.} {\sl Let $X$ and $Y$ be local sections 
of ${\cal H}$ such that $g(X,Y)$ is constant.  
Then we have 
$g(\nabla_{\frac{\partial}{\partial t}}X,Y)
+g(X,\nabla_{\frac{\partial}{\partial t}}Y)
=2\vert\vert H\vert\vert h(X,Y)$.}

\vspace{0.3truecm}

\noindent
{\it Proof.} 
From Lemma 4.2, we have 
$$\begin{array}{l}
\displaystyle{\frac{\partial}{\partial t}g(X,Y)
=\frac{\partial g}{\partial t}(X,Y)
+g(\nabla_{\frac{\partial}{\partial t}}X,Y)
+g(X,\nabla_{\frac{\partial}{\partial t}}Y)}\\
\hspace{1.9truecm}\displaystyle{=-2\vert\vert H\vert\vert h(X,Y)
+g(\nabla_{\frac{\partial}{\partial t}}X,Y)+g(X,\nabla_{\frac{\partial}
{\partial t}}Y).}
\end{array}$$
Hence the desired relation follows from the constancy of $g(X,Y)$.  
\hspace{2truecm}q.e.d.

\vspace{0.5truecm}

Next we prepare the following lemma for ${\cal R}$.  

\vspace{0.5truecm}

\noindent
{\bf Lemma 4.8.} {\sl For $X,Y\in{\cal H}$, we have 
$$\begin{array}{l}
\displaystyle{{\cal R}(X,Y)=2{\rm Tr}^{\bullet}_{g_{\cal H}}
\left(\langle({\cal A}^{\phi}_{\bullet}X,{\cal A}^{\phi}_{\bullet}
(A_{\cal H}Y)\rangle
+\langle({\cal A}^{\phi}_{\bullet}Y,{\cal A}^{\phi}_{\bullet}(A_{\cal H}X)
\rangle\right)}\\
\hspace{1.95truecm}\displaystyle{+2{\rm Tr}^{\bullet}_{g_{\cal H}}
\left(\langle({\cal A}^{\phi}_{\bullet}X,{\cal A}^{\phi}_Y(A_{\cal H}\bullet)
\rangle+\langle({\cal A}^{\phi}_{\bullet}Y,{\cal A}^{\phi}_X
(A_{\cal H}\bullet)\rangle\right)}\\
\hspace{1.95truecm}\displaystyle{+2{\rm Tr}^{\bullet}_{g_{\cal H}}
\left(\langle(\widetilde{\nabla}_{\bullet}{\cal A}^{\phi})_{\xi}Y,
{\cal A}^{\phi}_{\bullet}
X\rangle+\langle(\widetilde{\nabla}_{\bullet}{\cal A}^{\phi})_{\xi}X,
{\cal A}^{\phi}_{\bullet}Y
\rangle\right)}\\
\hspace{1.95truecm}\displaystyle{+{\rm Tr}^{\bullet}_{g_{\cal H}}
\left(\langle(\widetilde{\nabla}_{\bullet}{\cal A}^{\phi})_{\bullet}X,
{\cal A}^{\phi}_{\xi}Y\rangle
+\langle(\widetilde{\nabla}_{\bullet}{\cal A}^{\phi})_{\bullet}Y,
{\cal A}^{\phi}_{\xi}X\rangle\right)}\\
\hspace{1.95truecm}\displaystyle{+2{\rm Tr}^{\bullet}_{g_{\cal H}}
\langle{\cal T}^{\phi}_{{\cal A}^{\phi}_{\bullet}X}\xi,
{\cal A}^{\phi}_{\bullet}Y\rangle,}
\end{array}\leqno{(4.5)}$$
where we omit $F_{\ast}$.
}

\vspace{0.5truecm}

\noindent
{\it Proof.} Take $e,X,Y\in{\cal H}$.  
Easily we have 
$$\begin{array}{l}
\displaystyle{(\nabla_eh)({\cal A}_eX,Y)=e(\langle{\cal A}^{\phi}_eX,
{\cal A}^{\phi}_{\xi}Y\rangle)
-h({\cal A}_e({\cal A}_eX),Y)}\\
\hspace{3.2truecm}\displaystyle{-h((\nabla_e{\cal A})_eX,Y)
-h({\cal A}_eX,{\cal A}_eY).}
\end{array}\leqno{(4.6)}$$
On the other hand, by simple calculation, we have 
$((\widetilde{\nabla}_e{\cal A}^{\phi})_X\xi)_{\widetilde{\cal V}}
=-((\widetilde{\nabla}_e{\cal A}^{\phi})_{\xi}X)_{\widetilde{\cal V}}$.  
By using this relation, we can show 
$$e(\langle{\cal A}^{\phi}_eX,{\cal A}^{\phi}_{\xi}Y\rangle)
=\langle(\widetilde{\nabla}_e{\cal A}^{\phi})_eX,{\cal A}^{\phi}_{\xi}Y)
+\langle(\widetilde{\nabla}_e{\cal A}^{\phi})_{\xi}Y,{\cal A}^{\phi}_eX\rangle
+h({\cal A}_Y{\cal A}_Xe,e\rangle.\leqno{(4.7)}$$
Also, by simple calculations, we have 
$$\begin{array}{l}
\displaystyle{h({\cal A}_e({\cal A}_eX),Y)=-\langle{\cal A}^{\phi}_eX,
{\cal A}^{\phi}_e(A_{\cal H}Y)\rangle}\\
\displaystyle{h({\cal A}_Y({\cal A}_Xe),e)=\langle{\cal A}^{\phi}_eX,
{\cal A}^{\phi}_Y(A_{\cal H}e)\rangle,}\\
\displaystyle{h((\nabla_e{\cal A})_eX,Y)=\langle
(\widetilde{\nabla}_e{\cal A}^{\phi})_eX,{\cal A}^{\phi}_{\xi}Y\rangle,}\\
\displaystyle{h({\cal A}_eX,{\cal A}_eY)=-\langle
{{\cal T}^{\phi}}_{{\cal A}^{\phi}_eX}\xi,{\cal A}^{\phi}_eY\rangle.}
\end{array}\leqno{(4.8)}$$
From $(4.6),(4.7)$ and the relations in $(4.8)$, we have the desired 
relation.\hspace{1.5truecm}q.e.d.

\vspace{0.5truecm}

Also, we prepare the following lemma.  

\vspace{0.5truecm}

\noindent
{\bf Lemma 4.9.} {\sl For $X,Y,Z\in{\cal H}$, we have 
$$\begin{array}{l}
\displaystyle{2\langle{\cal T}^{\phi}_{{\cal A}^{\phi}_XY}\xi,
{\cal A}^{\phi}_XZ\rangle
=-\langle{\cal A}^{\phi}_XZ,(\widetilde{\nabla}_X{\cal A}^{\phi})_{\xi}Y
\rangle+\langle{\cal A}^{\phi}_XZ,(\widetilde{\nabla}_Y{\cal A}^{\phi})_{\xi}X
\rangle.}
\end{array}$$
}

\vspace{0.5truecm}

\noindent
{\it Proof.} 
Fix $(u_0,t_0)\in M\times[0,T)$.  Let $\widetilde X$ be an element of 
$\Gamma({\cal H})$ satisfying $\widetilde X_{(u_0,t_0)}=X$ and 
$(\widetilde{\nabla}^{\cal H}\widetilde X)_{(u_0,t_0)}=0$.  
Let $\widetilde Y$ and $\widetilde Z$ be similar elements of 
$\Gamma({\cal H})$ for $Y$ and $Z$, respectively.  
At $(u_0,t_0)$, we have 
$$\begin{array}{l}
\displaystyle{\langle{\cal A}^{\phi}_XZ,{\cal A}^{\phi}_Y(A_{\cal H}X)\rangle
=-\langle{\cal A}^{\phi}_XZ,\widetilde{\nabla}_Y(\widetilde{\nabla}_X\xi)
-\widetilde{\nabla}_Y({\cal A}^{\phi}_X\xi)\rangle}\\
\hspace{3.2truecm}\displaystyle{
=-\langle{\cal A}^{\phi}_XZ,\widetilde{\nabla}_X(\widetilde{\nabla}_Y\xi)
+\widetilde{\nabla}_{[Y,X]}\xi-\widetilde{\nabla}_Y({\cal A}^{\phi}_X\xi)
\rangle}\\
\hspace{3.2truecm}\displaystyle{
=\langle{\cal A}^{\phi}_XZ,\widetilde{\nabla}_X(AY)\rangle
-2\langle{\cal A}^{\phi}_XZ,\widetilde{\nabla}_{{\cal A}^{\phi}_YX}\xi\rangle
}\\
\hspace{3.7truecm}\displaystyle{
+\langle{\cal A}^{\phi}_XZ,(\widetilde{\nabla}_Y{\cal A}^{\phi})_X\xi\rangle
-\langle{\cal A}^{\phi}_XZ,{\cal A}^{\phi}_X(AY)\rangle}\\
\hspace{3.2truecm}\displaystyle{
=\langle{\cal A}^{\phi}_XZ,{\cal A}^{\phi}_Y(A_{\cal H}X)\rangle
-\langle{\cal A}^{\phi}_XZ,(\widetilde{\nabla}_X{\cal A}^{\phi})_Y\xi\rangle
}\\
\hspace{3.7truecm}\displaystyle{
-2\langle{\cal A}^{\phi}_XZ,{\cal T}^{\phi}_{{\cal A}^{\phi}_YX}\xi\rangle
+\langle{\cal A}^{\phi}_XZ,(\widetilde{\nabla}_Y{\cal A}^{\phi})_X\xi\rangle,}
\end{array}\leqno{(4.9)}
$$
where we use $(\widetilde{\nabla}^{\cal H}\widetilde X)_{(u_0,t_0)}
=(\widetilde{\nabla}^{\cal H}\widetilde Y)_{(u_0,t_0)}
=(\widetilde{\nabla}^{\cal H}\widetilde Z)_{(u_0,t_0)}=0$.  
Also we have
$$\langle{\cal A}^{\phi}_XZ,(\widetilde{\nabla}_X{\cal A}^{\phi})_Y\xi\rangle
=-\langle{\cal A}^{\phi}_XZ,(\widetilde{\nabla}_X{\cal A}^{\phi})_{\xi}Y\rangle
$$
and
$$\langle{\cal A}^{\phi}_XZ,(\widetilde{\nabla}_Y{\cal A}^{\phi})_X\xi\rangle
=-\langle{\cal A}^{\phi}_XZ,(\widetilde{\nabla}_Y{\cal A}^{\phi})_{\xi}X
\rangle.$$
Form $(4.9)$ and these relations, we obtain the desired relation.  
\hspace{2.1truecm}q.e.d.

\vspace{0.5truecm}

\noindent
{\bf Lemma 4.10.} 
{\sl For $X\in{\cal H}$, we have 
$$\begin{array}{l}
\displaystyle{{\cal R}(X,X)=4{\rm Tr}^{\bullet}_{g_{\cal H}}
\langle{\cal A}^{\phi}_{\bullet}X,{\cal A}^{\phi}_{\bullet}(A_{\cal H}X)
\rangle
+4{\rm Tr}^{\bullet}_{g_{\cal H}}\langle{\cal A}^{\phi}_{\bullet}X,
{\cal A}^{\phi}_X(A_{\cal H}\bullet)\rangle}\\
\hspace{1.95truecm}\displaystyle{
+3{\rm Tr}^{\bullet}_{g_{\cal H}}\langle(\widetilde{\nabla}_{\bullet}
{\cal A}^{\phi})_{\xi}X,{\cal A}^{\phi}_{\bullet}X\rangle
+2{\rm Tr}^{\bullet}_{g_{\cal H}}\langle(\widetilde{\nabla}_{\bullet}
{\cal A}^{\phi})_{\bullet}X,{\cal A}^{\phi}_{\xi}X
\rangle}\\
\hspace{1.95truecm}\displaystyle{+{\rm Tr}^{\bullet}_{g_{\cal H}}
\langle{\cal A}^{\phi}_{\bullet}X,
(\widetilde{\nabla}_X{\cal A}^{\phi})_{\xi}\bullet\rangle}
\end{array}$$
and hence 
$$
{\rm Tr}^{\bullet}_{g_{\cal H}}{\cal R}(\bullet,\bullet)=0.
$$
}

\vspace{0.5truecm}

\noindent
{\it Proof.} 
The first relation follows from the relations in Lemmas 4.8 and 4.9 directly.  
Also, the second relation follows from the frist relation directly.  
\hspace{2.5truecm}q.e.d.

\vspace{0.5truecm}

By using Theorem 4.6 and Lemmas 4.7 and 4.10, 
we can show the following evolution equation for $\vert\vert H_t\vert\vert$'s.  

\vspace{0.5truecm}

\noindent
{\bf Corollary 4.11.} {\sl The norms $\vert\vert H_t\vert\vert$'s of $H_t$ 
satisfy the following evolution equation:
$$
\frac{\partial\vert\vert H\vert\vert}{\partial t}=
\triangle_{{\cal H}}\vert\vert H\vert\vert+\vert\vert H\vert\vert
{\rm Tr}(A_{\cal H})^2-3\vert\vert H\vert\vert
{\rm Tr}(({\cal A}^{\phi}_{\xi})^2)_{\cal H}
.$$
}

\vspace{0.5truecm}

\noindent
{\it Proof.} Fix $(u_0,t_0)\in M\times[0,T)$.  
Take a local orthonormal frame field $\{e_1,\cdots,e_n\}$ of ${\cal H}$ (with respect to $g$) over 
a neighborhood $U$ of $(u_0,t_0)$ consisting of the eigenvectors of 
$A_{\cal H}$.  Since the fibres of $\phi$ are minimal regularizable submanifolds, we have 
$\vert\vert H\vert\vert=\sum\limits_{i=1}^nh(e_i,e_i)$ on $U$.  
Clearly we have 
$$\frac{\partial\vert\vert H\vert\vert}{\partial t}=
\sum_{i=1}^n\left(\frac{\partial h_{\cal H}}{\partial t}(e_i,e_i)
+2h_{\cal H}(\nabla_{\frac{\partial}{\partial t}}e_i,e_i)\right).
\leqno{(4.10)}$$
On the other hand, it follows from Theorem 4.6 that 
$$\begin{array}{l}
\displaystyle{\sum_{i=1}^n\frac{\partial h_{\cal H}}{\partial t}(e_i,e_i)=
\triangle_{\cal H}\vert\vert H\vert\vert
-\vert\vert H\vert\vert{\rm Tr}(A_{\cal H})^2-3\vert\vert H\vert\vert
{\rm Tr}(({\cal A}^{\phi}_{\xi})^2)_{\cal H},}
\end{array}
\leqno{(4.11)}$$
where we use 
$\sum\limits_{i=1}^n(\triangle_{\cal H}^{\cal H}h_{\cal H})(e_i,e_i)
=\triangle_{\cal H}\vert\vert H\vert\vert$ and 
${\rm Tr}_{g_{\cal H}}^{\bullet}{\cal R}(\bullet,\bullet)=0$ (by Lemma 4.10).  
Since each $e_i$ is an eigenvector of $A_{\cal H}$, we have 
$h(e_i,e_j)=0\,\,(i\not=j)$.  By using Lemma 4.7, we can show 
$$\sum_{i=1}^nh_{\cal H}(\nabla_{\frac{\partial}{\partial t}}e_i,e_i)
=\sum_{i=1}^ng(\nabla_{\frac{\partial}{\partial t}}e_i,e_i)h(e_i,e_i)
=\vert\vert H\vert\vert{\rm Tr}(A_{\cal H})^2.\leqno{(4.12)}$$
From $(4.10),\,(4.11)$ and $(4.12)$, we obtain the desired relation.  
\hspace{2.3truecm}q.e.d.

\vspace{0.5truecm}

From we derive the following evolution equation  for ${\rm Tr}(A_{\cal H})_t^2$.  

\vspace{0.5truecm}

\noindent
{\bf Corollary 4.12.} {\sl The quantities ${\rm Tr}(A_{\cal H})_t^2$'s 
satisfy the following evolution equation:
$$\begin{array}{l}
\displaystyle{\frac{\partial{\rm Tr}(A_{\cal H})^2}{\partial t}
=\triangle_{\cal H}({\rm Tr}(A_{\cal H})^2)
-2{\rm Tr}{\rm Tr}^{\bullet}_{g_{\cal H}}
(\nabla_{\bullet}^{\cal H}A_{\cal H}\circ\nabla_{\bullet}^{\cal H}A_{\cal H})
}\\
\hspace{2.2truecm}\displaystyle{
+2{\rm Tr}((A_{\cal H})^2)\left(
{\rm Tr}((A_{\cal H})^2)-{\rm Tr}(({\cal A}^{\phi}_{\xi})^2)_{\cal H}\right)}\\
\hspace{2.2truecm}\displaystyle{
-4\vert\vert H\vert\vert{\rm Tr}\left((({\cal A}^{\phi}_{\xi})^2)\circ 
A_{\cal H}\right)
-2{\rm Tr}^{\bullet}_{g_{\cal H}}{\cal R}(A_{\cal H}\bullet,\bullet).}
\end{array}$$
}

\vspace{0.5truecm}

\noindent
{\it Proof.} Fix $(u_0,t_0)\in M\times[0,T)$.  
Take a local orthonormal frame field $\{e_1,\cdots,e_n\}$ of ${\cal H}$ 
(with respect to $g_{\cal H}$) over a neighborhood $U$ of $(u_0,t_0)$ consisting of 
the eigenvectors of $A_{\cal H}$.  From Lemma 4.2, we have 
$$\begin{array}{l}
\displaystyle{\frac{\partial h_{\cal H}}{\partial t}(X,Y)
=\frac{\partial g_{\cal H}}{\partial t}(A_{\cal H}X,Y)
+g_{\cal H}(\frac{\partial A_{\cal H}}{\partial t}(X),Y)}\\
\hspace{2.1truecm}\displaystyle{=-2\vert\vert H\vert\vert h_{\cal H}
(A_{\cal H}X,Y)+g_{\cal H}(\frac{\partial A_{\cal H}}{\partial t}(X),Y)}
\end{array}
\leqno{(4.13)}$$
for any $X,Y\in\pi_M^{\ast}TM$.  
Since $\{e_1,\cdots,e_n\}$ consists of the eigenvectors of $A_{\cal H}$, 
it follows from Lemma 4.7 that 
$$g(\nabla_{\frac{\partial}{\partial t}}e_i,e_i)=\vert\vert H\vert\vert 
h(e_i,e_i).\leqno{(4.14)}$$
From these relations, Lemmas 4.2 and 4.7, 
we have 
$$\begin{array}{l}
\hspace{0.6truecm}
\displaystyle{\frac{\partial{\rm Tr}(A_{\cal H})^2}{\partial t}
=\sum_{i=1}^n\frac{\partial}{\partial t}(h_{\cal H}(A_{\cal H}e_i,e_i))}\\
\displaystyle{=\sum_{i=1}^n\left(\frac{\partial h_{\cal H}}{\partial t}
(A_{\cal H}e_i,e_i)+h_{\cal H}(\frac{\partial A_{\cal H}}{\partial t}(e_i),
e_i)+2h_{\cal H}(A_{\cal H}e_i,\nabla_{\frac{\partial}{\partial t}}e_i)\right)
}\\
\displaystyle{=\sum_{i=1}^n\left(\frac{\partial h_{\cal H}}{\partial t}
(A_{\cal H}e_i,e_i)+g_{\cal H}((\frac{\partial A_{\cal H}}{\partial t}(e_i),
A_{\cal H}e_i)\right.}\\
\hspace{1.4truecm}\displaystyle{
\left.+2\vert\vert H\vert\vert h(e_i,e_i)h_{\cal H}(A_{\cal H}e_i,e_i)
\right)}\\
\displaystyle{=\sum_{i=1}^n\left(2\frac{\partial h_{\cal H}}{\partial t}
(A_{\cal H}e_i,e_i)+2\vert\vert H\vert\vert g((A_{\cal H})^3e_i,e_i)\right.}\\
\hspace{1.4truecm}\displaystyle{
\left.+2\vert\vert H\vert\vert h(e_i,e_i)g((A_{\cal H})^2e_i,e_i)\right)}\\
\displaystyle{=\sum_{i=1}^n\left(2\frac{\partial h_{\cal H}}{\partial t}
(A_{\cal H}e_i,e_i)+4\vert\vert H\vert\vert 
g((A_{\cal H})^3e_i,e_i)\right).}
\end{array}\leqno{(4.15)}$$
Also we have 
$$\sum_{i=1}^n(\triangle_{\cal H}^{\cal H}h_{\cal H})(A_{\cal H}e_i,e_i)
=\frac12\triangle_{\cal H}{\rm Tr}((A_{\cal H})^2)
-{\rm Tr}{\rm Tr}^{\bullet}_{g_{\cal H}}\left(
\nabla^{\cal H}_{\bullet}A_{\cal H}\circ\nabla^{\cal H}_{\bullet}A_{\cal H}
\right).\leqno{(4.16)}$$
From Theorem 4.6, $(4.15)$ and $(4.16)$, 
we obtain the desired relation.  
\hspace{1.5truecm}q.e.d.

\vspace{0.5truecm}

By using Corollaries 4.11 and 4.12, we can show the following evolution equation.  

\vspace{0.5truecm}

\noindent
{\bf Corollary 4.13.} 
{\sl The quantities ${\rm Tr}\,(A_{\cal H})_t^2
-\frac{\vert\vert H_t\vert\vert^2}{n}$'s satisfy the following evolution equation:
$$\begin{array}{l}
\displaystyle{
\frac{\partial({\rm Tr}\,(A_{\cal H})^2
-\frac{\vert\vert H\vert\vert^2}{n})}{\partial t}=
\triangle_{\cal H}\left({\rm Tr}(A_{\cal H})^2
-\frac{\vert\vert H\vert\vert^2}{n}\right)
+\frac2n\left\vert\left\vert{\rm grad}
\vert\vert H\vert\vert\right\vert\right\vert^2}\\
\hspace{3.8truecm}\displaystyle{
+2{\rm Tr}(A_{\cal H})^2\times\left({\rm Tr}(A_{\cal H})^2
-\frac{\vert\vert H\vert\vert^2}{n}\right)}\\
\hspace{3.8truecm}\displaystyle{
-2{\rm Tr}{\rm Tr}_{g_{\cal H}}^{\bullet}\left(\nabla^{\cal H}A_{\cal H}\circ
\nabla^{\cal H}A_{\cal H}\right)}\\
\hspace{3.8truecm}\displaystyle{
-2{\rm Tr}(({\cal A}^{\phi}_{\xi})^2)_{\cal H}\times
\left({\rm Tr}(A_{\cal H})^2
-\frac{\vert\vert H\vert\vert^2}{n}\right)}\\
\hspace{3.8truecm}\displaystyle{
-4\vert\vert H\vert\vert\left({\rm Tr}\left(({\cal A}^{\phi}_{\xi})^2\circ
\left(A_{\cal H}-\frac{\vert\vert H\vert\vert}{n}{\rm id}\right)\right)\right)
}\\
\hspace{3.8truecm}\displaystyle{
-2{\rm Tr}^{\bullet}_{g_{\cal H}}{\cal R}\left(\left(
A_{\cal H}-\frac{\vert\vert H\vert\vert}{n}{\rm id}\right)\bullet,\bullet
\right),}
\end{array}$$
where ${\rm grad}\vert\vert H\vert\vert$ is the gradient vector field of 
$\vert\vert H\vert\vert$ with respect to $g$ and 
$\left\vert\left\vert{\rm grad}\vert\vert H\vert\vert\right\vert\right\vert$ 
is the norm of ${\rm grad}\vert\vert H\vert\vert$ with respect to $g$.}

\vspace{0.5truecm}

\noindent
{\it Proof.} This relation follows directly from Corollaries 4.11, 4.12 and 
$\triangle_{\cal H}\vert\vert H\vert\vert^2=$\newline
$2\vert\vert H\vert\vert\triangle_{\cal H}\vert\vert H\vert\vert
+2\left\vert\left\vert{\rm grad}\vert\vert H\vert\vert
\right\vert\right\vert^2$.  
\hspace{7.35truecm}q.e.d.

\vspace{0.5truecm}

\noindent
{\it Remark 4.1.} From the evolution equations obtained in this section, 
the evolution equations for the corresponding geometric quantities of 
$\overline f_t(:\overline M\hookrightarrow V/G)$ are derived, respectively. 
In the case where the $G$-action is free and hence $V/G$ is a (complete) Riemannian manifold, these 
derived evolution equations coincide with the evolution equations for the corresponding geometric 
quantities along the mean curvature flow in a complete Riemannian manifold which were given 
by Huisken [Hu2].  
That is, the discussion in this section give a new proof of the evolution equations in [Hu2] 
in the case where the ambient complete Riemannian manifold occurs as $V/G$.  
In the proof of [Hu2], one need to take local coordinates of the ambient space to derive the evolution 
equations.  
On the other hand, in our proof, one need not take local coordinates of the ambient space because 
the ambient space is a Hilbert space.  This is an advantage of our proof.  

\section{A maximum principle} 
Let $M$ be a Hilbert manifold and $g_t$ ($0\leq t<T$) a $C^{\infty}$-family of 
Riemannian metrics on $M$ and $G\curvearrowright M$ 
a almost free action which is isometric with respect to $g_t$'s ($t\in[0,T)$).  
Assume that the orbit space $M/G$ is compact.  
Let ${\cal H}_t$ ($0\leq t<T$) be the horizontal distribution of the $G$-action and 
define a subbundle ${\cal H}$ of $\pi_M^{\ast}TM$ by 
${\cal H}_{(x,t)}:=({\cal H}_t)_x$.  
For a tangent vector field $X$ on $M$ (or an open set $U$ of $M$), we define 
a section $\bar X$ of $\pi_M^{\ast}TM$ 
(or $\pi_M^{\ast}TM\vert_U$) by $\bar X_{(x,t)}:=X_x$ ($(x,t)\in M\times[0,T)$).  
Let $\nabla^t$ ($0\leq t<T$) be the Riemannian connection of $g_t$ and 
$\nabla$ the connection of $\pi_M^{\ast}TM$ defined in terms of $\nabla^t$'s ($t\in[0,T)$).  
Define a connection $\nabla^{\cal H}$ of $\cal H$ by 
$\nabla^{\cal H}_XY={\rm pr}_{\cal H}(\nabla_XY)$ 
for any $X\in T(M\times[0,T))$ and any $Y\in\Gamma({\cal H})$.  
For $B\in\Gamma(\pi_M^{\ast}T^{(r_0,s_0)}M)$, 
we define maps $\psi_{B\otimes}$ and $\psi_{\otimes B}$ 
from $\Gamma(\pi_M^{\ast}T^{(r,s)}M)$ to 
$\Gamma(\pi_M^{\ast}T^{(r+r_0,s+s_0)}M)$ by 
$$\psi_{B\otimes}(S):=B\otimes S,\,\,\,\,{\rm and}\,\,\,\,
\psi_{\otimes B}(S):=S\otimes B\,\,\,\,(S\in\Gamma(\pi_M^{\ast}T^{(r,s)}M),$$
respectively.  
Also, we define a map $\psi_{\otimes^k}$ of 
$\Gamma(\pi_M^{\ast}T^{(r,s)}M)$ to 
$\Gamma(\pi_M^{\ast}T^{(kr,ks)}M)$ by 
$$\psi_{\otimes^k}(S):=S\otimes\cdots\otimes S\,\,(k{\rm -times})\,\,\,\,
(S\in\Gamma(\pi_M^{\ast}T^{(r,s)}M).$$
Also, we define a map $\psi_{g_{\cal H},ij}\,\,(i<j)$ from 
$\Gamma(\pi_M^{\ast}T^{(0,s)}M)$ (or $\Gamma(\pi_M^{\ast}T^{(1,s)}M)$) to 
$\Gamma(\pi_M^{\ast}T^{(0,s-2)}M)$ (or $\Gamma(\pi_M^{\ast}T^{(1,s-2)}M)$) by 
$$\begin{array}{l}
\hspace{0.4truecm}\displaystyle{(\psi_{g_{\cal H},ij}(S))_{(x,t)}(X_1,\cdots,
X_{s-2})}\\
\displaystyle{:=\sum_{k=1}^nS_{(x,t)}(X_1,\cdots,X_{i-1},e_k,X_{i+1},\cdots,
X_{j-1},e_k,X_{j+1},\cdots,X_{s-2})}
\end{array}$$
and define a map $\psi_{{\cal H},i}$ from $\Gamma(\pi_M^{\ast}T^{(1,s)}M)$ 
to $\Gamma(\pi_M^{\ast}T^{(0,s-1)}M)$ by 
$$(\psi_{{\cal H},i}(S))_{(x,t)}(X_1,\cdots,X_{s-1})
:={\rm Tr}\,S_{(x,t)}(X_1,\cdots,X_{i-1},\bullet,X_i,\cdots,X_{s-1}),$$
where $X_i\in T_xM$ ($i=1,\cdots,s-1$) and $\{e_1,\cdots,e_n\}$ is an 
orthonormal base of $({\cal H}_t)_x$ with respect to $g_t$.  
We call a map $P$ from $\Gamma(\pi_M^{\ast}T^{(0,s)}M)$ to oneself given by 
the composition of the above maps of five type {\it a map of polynomial type}.  

In this section, we prove the following maximum principle for 
a $C^{\infty}$-family of $G$-invariant symmetric $(0,2)$-tensor fields on $M$.  

\vspace{0.5truecm}

\noindent
{\bf Theorem 5.1.} {\sl Let $S\in\Gamma(\pi_M^{\ast}(T^{(0,2)}M))$ 
such that, for each $t\in[0,T)$, $S_t(:=S_{(\cdot,t)})$ is a $G$-invariant 
symmetric $(0,2)$-tensor field on $M$.  Assume that $S_t$'s 
($0\leq t<T$) satisfy the following evolution equation:
$$\frac{\partial S_{\cal H}}{\partial t}=\triangle_{\cal H}^{\cal H}
S_{\cal H}+\nabla^{\cal H}_{\bar X_0}S_{\cal H}+P(S)_{\cal H},\leqno{(5.1)}$$
where $X_0\in\Gamma(TM)$ and $P$ is a map of polynomial type from 
$\Gamma(\pi_M^{\ast}(T^{(0,2)}M))$ to oneself.  

${\rm(i)}$ Assume that $P$ satisfies the following condition:
$$\begin{array}{r}
\displaystyle{X\in{\rm Ker}((S+\varepsilon g)_{\cal H})_{(x,t)}\,\,\,
\Rightarrow\,\,\,P(S+\varepsilon g)_{(x,t)}(X,X)\geq0}\\
\displaystyle{(\forall\,\varepsilon>0,\,\,\,\,(x,t)\in M\times[0,T)).}
\end{array}\leqno{(\ast_{S_{\cal H}}^+)}$$
If $(S_{\cal H})_{(\cdot,0)}\geq0$ (resp. $>0$), 
then $(S_{\cal H})_{(\cdot,t)}\geq0$ (resp. $>0$) holds for all $t\in[0,T)$.  

${\rm(ii)}$ Assume that $P$ satisfies the following condition:
$$\begin{array}{r}
\displaystyle{X\in{\rm Ker}((S+\varepsilon g)_{\cal H})_{(x,t)}\,\,\,
\Rightarrow\,\,\,P(S+\varepsilon g)_{(x,t)}(X,X)\leq0}\\
\displaystyle{(\forall\,\varepsilon>0,\,\,\,\,(x,t)\in M\times[0,T)).}
\end{array}\leqno{(\ast_{S_{\cal H}}^-)}$$
If $(S_{\cal H})_{(\cdot,0)}\leq0$ (resp. $<0$), then $(S_{\cal H})_{(\cdot,t)}\leq0$ (resp. $<0$) 
holds for all $t\in[0,T)$.}

\vspace{0.5truecm}

\noindent
{\it Proof.} First we shall show the part of 
``If $(S_{\cal H})_{(\cdot,0)}\geq0$, 
then $(S_{\cal H})_{(\cdot,t)}\geq0$ holds for all $t\in(0,T)$" in the statement (i).  
For positive numbers $\varepsilon$ and $\delta$, we define 
$S_{\varepsilon,\delta}$ by 
$(S_{\varepsilon,\delta})_{(x,t)}:=S_{(x,t)}+\varepsilon(\delta+t)
g_{(x,t)}$.  

(Step I) In this step, we show the following statement:
$$(\ast)\,\,\,\,\exists\,\delta>0\,\,{\rm s.t.}\,\,{\rm``}
((S_{\varepsilon,\delta})_{\cal H})_{(x,t)}>0\,\,
(\forall\,(x,t)\in M\times[0,\delta),\forall\,\varepsilon>0){\rm"}.$$
Suppose that such a positive number $\delta$ does not exists.  
Fix a sufficiently small positive number $\delta$.  
For some $\varepsilon_0>0$, 
there exists $(x_0,t_0)\in M\times[0,\delta)$ such that 
$((S_{\varepsilon_0,\delta})_{\cal H})_{(x_0,t_0)}=0$.  
Here we take $t_0$ as smally as possible.  We have ${\rm Ker}
((S_{\varepsilon_0,\delta})_{\cal H})_{(x_0,t_0)}$\newline
$\not=\{0\}$ and 
$((S_{\varepsilon_0,\delta})_t)_{{\cal H}_t}\geq0\,\,
(\forall\,t\in[0,t_0])$.  Take $v_1\in{\rm Ker}
((S_{\varepsilon_0,\delta})_{\cal H})_{(x_0,t_0)})$ with \newline
$g_{(x_0,t_0)}(v_1,v_1)=1$.  From the assumption $(\ast_{S_{\cal H}}^+)$ for $P$, 
we have 
$$P((S_{\varepsilon_0,\delta})_{(x_0,t_0)})(v_1,v_1)\geq 0.
\leqno{(5.2)}$$
The map $P$ is of polynomial type, $M/G$ is compact and $S_t$ is $G$-invariant.  
Hence, for each $t\in[0,T)$, there exists a positive constant $C_{\delta,t}$ 
(depending only on $\vert\vert(S_{\cal H})_{(\cdot,t)}\vert\vert$ and 
$\vert\vert((S_{(\varepsilon_0,\delta)})_{\cal H})_{(\cdot,t)}\vert\vert$) such that 
$$\vert\vert((P(S_{\varepsilon_0,\delta}))_{\cal H})_{(\cdot,t)}
-((P(S))_{\cal H})_{(\cdot,t)}\vert\vert\leq C_{\delta,t}\vert\vert
((S_{\varepsilon_0,\delta})_{\cal H})_{(\cdot,t)}-(S_{\cal H})_{(\cdot,t)}
\vert\vert\leqno{(5.3)}$$
on $M$, where $\vert\vert\cdot\vert\vert$ is the pointwise norm of a tensor field $(\cdot)$.  
We take $C_{\delta,t}$ as smally as possible.  
Since $P$ is of polynomial type, 
$\lim_{\delta\to+0}C_{\delta,t}$ exists and $\lim_{\delta\to+0}C_{\delta,t}
>0$.  Denote by $C_t$ this limit.  Fix $T_1\in(t_0,T)$.  Set 
$$C_{\delta}:=\max
\left\{
\mathop{\max}_{0\leq t\leq T_1}\,C_{\delta,t},\quad\,\,\,\,
\mathop{\max}_
{\begin{array}{c}
\displaystyle{(x,t)\in M\times[0,T_1]}\\
\displaystyle{v\in TM\,\,{\rm s.t.}\,\,g_t(v,v)=1}
\end{array}}
\left\vert\left(\frac{\partial g_{\cal H}}{\partial t}\right)_{(x,t)}(v,v)
\right\vert\right\}$$
and 
$$C:=\max
\left\{
\mathop{\max}_{0\leq t\leq T_1}\,C_t,\quad\,\,\,\,
\mathop{\max}_
{\begin{array}{c}
\displaystyle{(x,t)\in M\times[0,T_1]}\\
\displaystyle{v\in TM\,\,{\rm s.t.}\,\,g_t(v,v)=1}
\end{array}}
\left\vert\left(\frac{\partial g_{\cal H}}{\partial t}\right)_{(x,t)}(v,v)
\right\vert\right\}.$$
Since $C$ is independent of the choice of $\delta$, we may assume that 
$C\delta<\frac14$ by replacing $\delta$ to a smaller positive number 
if necessary.  
Furthermore, since $\delta\mapsto C_{\delta}$ is upper semi-continuous and 
$\lim_{\delta\to+0}C_{\delta,t}<\infty$, 
we may assume that $C_{\delta}\delta<\frac14$ by replacing $\delta$ to 
a smaller positive number if necessary.  
From $(5.2)$ and $(5.3)$, we have 
$$P(S)_{(x_0,t_0)}(v_1,v_1)\geq-2C_{\delta}\varepsilon_0\delta.\leqno{(5.4)}$$
Let $X_1$ be a section of ${\cal H}$ on a normal neighborhood $U$ of 
$(x_0,t_0)$ in $M\times[0,T)$ such that $(X_1)_{(x_0,t_0)}=v_1$ and that 
$\nabla^{\cal H}X_1=0$ at $(x_0,t_0)$.  
Define a function $\rho$ on $U$ by 
$\rho(x,t):=(S_{\varepsilon_0,\delta})_{(x,t)}((\bar X_1)_{(x,t)},
(\bar X_1)_{(x,t)})$ ($(x,t)\in U$).  
Since we take $(x_0,t_0)$ and $v_1$ as above, we have 
$(\frac{\partial\rho}{\partial t})_{(x_0,t_0)}\leq0$ (see Figure 4).  
Also, we have 
$$\left(\frac{\partial\rho}{\partial t}\right)_{(x_0,t_0)}
=\left(\frac{\partial S_{\cal H}}{\partial t}\right)_{(x_0,t_0)}(v_1,v_1)
+\varepsilon_0(\delta+t_0)\left(\frac{\partial g_{\cal H}}{\partial t}\right)_{(x_0,t_0)}
(v_1,v_1)+\varepsilon_0.$$
Hence we have 
$$\left(\frac{\partial S_{\cal H}}{\partial t}\right)_{(x_0,t_0)}(v_1,v_1)
\leq-\varepsilon_0(\delta+t_0)\left(\frac{\partial g_{\cal H}}{\partial t}\right)
_{(x_0,t_0)}(v_1,v_1)-\varepsilon_0.\leqno{(5.5)}$$
Take $w\in T_{x_0}(M\times\{t_0\})$.  
Clearly we have $d\rho_{(x_0,t_0)}(w)=0$.  
Also we have 
$d\rho_{(x_0,t_0)}(w)=(\nabla^{\cal H}_w(S_{\varepsilon_0,\delta})_{\cal H})
_{(x_0,t_0)}(v_1,v_1)$.  Hence we have 
$$(\nabla^{\cal H}_w(S_{\varepsilon_0,\delta})_{\cal H})_{(x_0,t_0)}(v_1,v_1)
=0.\leqno{(5.6)}$$
Clearly we have 
$(\triangle_{t_0}\,\rho_{t_0})_{x_0}\geq0$, where $\triangle_{t_0}$ is the 
Laplacian operator with respect to $g_{t_0}$.  
Also, we have $(\triangle_{t_0}\,\rho_{t_0})_{x_0}
=(\triangle_{\cal H}^{\cal H}(S_{\varepsilon_0,\delta})_{\cal H})(v_1,v_1)$.  
Hence we have 
$$(\triangle_{\cal H}^{\cal H}(S_{\varepsilon_0,\delta})_{\cal H})(v_1,v_1)\geq0.
\leqno{(5.7)}$$
From $(5.1),\,(5.5),\,(5.6)$ and $(5.7)$, we have 
$$\begin{array}{l}
\displaystyle{P(S)_{(x_0,t_0)}(v_1,v_1)\leq-\varepsilon_0
+\varepsilon_0(\delta+t_0)
\left\vert(\frac{\partial g_{\cal H}}{\partial t})_{(x_0,t_0)}(v_1,v_1)
\right\vert}\\
\hspace{3.1truecm}\displaystyle{
\leq-\varepsilon_0+2\varepsilon_0 C_{\delta}\delta.}
\end{array}
\leqno{(5.8)}$$
From $(5.4)$ and $(5.8)$, we have $C_{\delta}\delta\geq\frac14$.  
This contradicts $C_{\delta}\delta<\frac14$.  
Therefore the statement $(\ast)$ is true.  

(Step II) Let $\delta$ be a positive number as in the statement $(\ast)$.  
Then, for any $(x,t)\in M\times[0,\delta)$ and any $\varepsilon>0$, we have 
$((S_{\varepsilon,\delta})_{\cal H})_{(x,t)}>0$.  
Hence we have $\lim\limits_{\varepsilon\to+0}
((S_{\varepsilon,\delta})_{\cal H})_{(x,t)}
=(S_{\cal H})_{(x,t)}\geq0$ for any 
$(x,t)\in M\times[0,\delta)$.  
Set 
$$T_1:=\sup\{t_1\,\vert\,
(S_{\cal H})_{(x,t)}\geq0\,\,(\forall\,(x,t)\in M\times[0,t_1]\}.$$
Suppose that $T_1<T$.  Then, by the similar discussion for 
$(S_{\cal H})_{(\cdot,T_1)}$ instead of 
$(S_{\cal H})_{(\cdot,0)}$, we can show that 
$(S_{\cal H})_{(x,t)}\geq0$ for any $t\in[T_1,T_1+\delta']$ 
and any $x\in M$, where $\delta'$ is some positive number.  This contradicts 
the definition of $T_1$.  Therefore we have $T_1=T$.  Thus we obtain 
$(S_{\cal H})_{(\cdot,t)}\geq0$ for any $t\in[0,T)$.  

Similarly we can show the part of ``If $(S_{\cal H})_{(\cdot,0)}>0$, 
then $(S_{\cal H})_{(\cdot,t)}>0$ holds for all $t\in(0,T)$" in the statement (i) as follows.  
The map $P$ is of polynomial type, $M/G$ is compact and $S_t$ is $G$-invariant.  
Hence it follows from $(S_{\cal H})_{(\cdot,0)}>0$ that 
$(S_{\cal H})_{(\cdot,0)}\geq b(g_{\cal H})_{(\cdot,0)}$ holds for some positive coonstant $b$.  
Set $\overline S:=S-bg$.  
Then it is easy to show that $\overline S$ satisfies 
$$\frac{\partial\overline S_{\cal H}}{\partial t}=\triangle_{\cal H}^{\cal H}
\overline S_{\cal H}+\nabla^{\cal H}_{\bar X_0}\overline S_{\cal H}+\overline P(\overline S)_{\cal H}$$
for some map $\overline P$ of polynomial type satisfying the condition $(\ast_{\overline S_{\cal H}}^+)$ :
$$\begin{array}{r}
\displaystyle{X\in{\rm Ker}((\overline S+\varepsilon g)_{\cal H})_{(x,t)}\,\,\,
\Rightarrow\,\,\,\overline P(\overline S+\varepsilon g)_{(x,t)}(X,X)\geq0}\\
\displaystyle{(\forall\,\varepsilon>0,\,\,\,\,(x,t)\in M\times[0,T)).}
\end{array}$$
Hence, it follows from Theorem 5.1 that $(\overline S_{\cal H})_{(\cdot,t)}\geq0$ 
(hence $(S_{\cal H})_{(\cdot,t)}>0$) holds for all $t\in[0,T)$.  
The statement (ii) also are derived by the similar discussion.  
\begin{flushright}q.e.d.\end{flushright}

\vspace{0.5truecm}

\centerline{
\unitlength 0.1in
\begin{picture}( 30.2400, 26.8400)(  4.0300,-30.1400)
%
\special{pn 8}%
\special{pa 1938 2180}%
\special{pa 1508 2636}%
\special{pa 2962 2636}%
\special{pa 3428 2180}%
\special{pa 3428 2180}%
\special{pa 1938 2180}%
\special{fp}%
%
\special{pn 8}%
\special{pa 1486 2950}%
\special{pa 2946 2950}%
\special{fp}%
%
\special{pn 8}%
\special{pa 1132 2640}%
\special{pa 1630 2166}%
\special{fp}%
%
\special{pn 8}%
\special{sh 1}%
\special{ar 1488 2954 10 10 0  6.28318530717959E+0000}%
\special{sh 1}%
\special{ar 1488 2954 10 10 0  6.28318530717959E+0000}%
%
\special{pn 13}%
\special{sh 1}%
\special{ar 1488 2954 10 10 0  6.28318530717959E+0000}%
\special{sh 1}%
\special{ar 1488 2954 10 10 0  6.28318530717959E+0000}%
%
\special{pn 13}%
\special{sh 1}%
\special{ar 2630 2954 10 10 0  6.28318530717959E+0000}%
\special{sh 1}%
\special{ar 2630 2954 10 10 0  6.28318530717959E+0000}%
%
\special{pn 13}%
\special{sh 1}%
\special{ar 1384 2402 10 10 0  6.28318530717959E+0000}%
\special{sh 1}%
\special{ar 1384 2402 10 10 0  6.28318530717959E+0000}%
%
\special{pn 8}%
\special{pa 1390 2402}%
\special{pa 2488 2402}%
\special{dt 0.045}%
%
\special{pn 8}%
\special{pa 2286 2636}%
\special{pa 2286 2950}%
\special{dt 0.045}%
%
\special{pn 13}%
\special{sh 1}%
\special{ar 2286 2954 10 10 0  6.28318530717959E+0000}%
\special{sh 1}%
\special{ar 2286 2954 10 10 0  6.28318530717959E+0000}%
%
\special{pn 13}%
\special{sh 1}%
\special{ar 2494 2402 10 10 0  6.28318530717959E+0000}%
\special{sh 1}%
\special{ar 2494 2402 10 10 0  6.28318530717959E+0000}%
%
\special{pn 8}%
\special{pa 1384 2050}%
\special{pa 1384 346}%
\special{fp}%
\special{sh 1}%
\special{pa 1384 346}%
\special{pa 1364 412}%
\special{pa 1384 398}%
\special{pa 1404 412}%
\special{pa 1384 346}%
\special{fp}%
%
\special{pn 13}%
\special{sh 1}%
\special{ar 1384 1684 10 10 0  6.28318530717959E+0000}%
\special{sh 1}%
\special{ar 1384 1684 10 10 0  6.28318530717959E+0000}%
%
\special{pn 8}%
\special{pa 1384 1684}%
\special{pa 2494 1684}%
\special{dt 0.045}%
%
\special{pn 8}%
\special{pa 2494 1684}%
\special{pa 2494 2394}%
\special{dt 0.045}%
%
\special{pn 13}%
\special{sh 1}%
\special{ar 2494 1684 10 10 0  6.28318530717959E+0000}%
\special{sh 1}%
\special{ar 2494 1684 10 10 0  6.28318530717959E+0000}%
%
\special{pn 8}%
\special{pa 2930 2950}%
\special{pa 2958 2950}%
\special{fp}%
\special{sh 1}%
\special{pa 2958 2950}%
\special{pa 2890 2930}%
\special{pa 2904 2950}%
\special{pa 2890 2970}%
\special{pa 2958 2950}%
\special{fp}%
%
\special{pn 20}%
\special{sh 1}%
\special{ar 2494 1684 10 10 0  6.28318530717959E+0000}%
\special{sh 1}%
\special{ar 2494 1684 10 10 0  6.28318530717959E+0000}%
\put(14.3900,-30.0600){\makebox(0,0)[lt]{$0$}}%
\put(25.9800,-30.1400){\makebox(0,0)[lt]{$\delta$}}%
\put(22.1600,-30.1400){\makebox(0,0)[lt]{$t_0$}}%
\put(10.5000,-26.7700){\makebox(0,0)[lt]{$M$}}%
\put(13.2800,-24.0200){\makebox(0,0)[rb]{$x_0$}}%
\put(13.3400,-16.2200){\makebox(0,0)[rt]{$0$}}%
\put(13.0300,-3.3000){\makebox(0,0)[rt]{${\Bbb R}$}}%
\put(23.1000,-7.0900){\makebox(0,0)[lb]{The graph of $\rho$}}%
%
\special{pn 8}%
\special{pa 1692 1432}%
\special{pa 1714 1410}%
\special{pa 1744 1406}%
\special{pa 1776 1414}%
\special{pa 1804 1428}%
\special{pa 1822 1454}%
\special{pa 1824 1474}%
\special{sp}%
%
\special{pn 8}%
\special{pa 1964 1506}%
\special{pa 1942 1528}%
\special{pa 1910 1530}%
\special{pa 1880 1524}%
\special{pa 1852 1508}%
\special{pa 1832 1482}%
\special{pa 1832 1464}%
\special{sp}%
%
\special{pn 8}%
\special{pa 1964 1516}%
\special{pa 1986 1496}%
\special{pa 2016 1492}%
\special{pa 2048 1498}%
\special{pa 2074 1516}%
\special{pa 2094 1540}%
\special{pa 2096 1560}%
\special{sp}%
%
\special{pn 8}%
\special{pa 2228 1586}%
\special{pa 2206 1606}%
\special{pa 2174 1610}%
\special{pa 2144 1602}%
\special{pa 2116 1586}%
\special{pa 2098 1562}%
\special{pa 2096 1542}%
\special{sp}%
%
\special{pn 8}%
\special{pa 2228 1588}%
\special{pa 2252 1572}%
\special{pa 2284 1570}%
\special{pa 2314 1578}%
\special{pa 2342 1592}%
\special{pa 2358 1618}%
\special{pa 2358 1622}%
\special{sp}%
%
\special{pn 8}%
\special{pa 2462 1672}%
\special{pa 2432 1662}%
\special{pa 2402 1648}%
\special{pa 2376 1632}%
\special{pa 2358 1610}%
\special{pa 2358 1604}%
\special{sp}%
%
\special{pn 8}%
\special{ar 2738 1586 358 130  1.7889727 2.4864515}%
%
\special{pn 8}%
\special{pa 2988 1592}%
\special{pa 2962 1610}%
\special{pa 2934 1626}%
\special{pa 2906 1640}%
\special{pa 2876 1654}%
\special{pa 2846 1664}%
\special{pa 2816 1674}%
\special{pa 2786 1684}%
\special{pa 2754 1692}%
\special{pa 2724 1698}%
\special{pa 2692 1702}%
\special{pa 2660 1708}%
\special{pa 2642 1708}%
\special{sp}%
%
\special{pn 8}%
\special{pa 1510 968}%
\special{pa 1532 948}%
\special{pa 1562 942}%
\special{pa 1594 946}%
\special{pa 1624 956}%
\special{pa 1650 974}%
\special{pa 1666 1002}%
\special{pa 1666 1018}%
\special{sp}%
%
\special{pn 8}%
\special{pa 1830 1054}%
\special{pa 1810 1074}%
\special{pa 1778 1080}%
\special{pa 1746 1076}%
\special{pa 1716 1066}%
\special{pa 1692 1048}%
\special{pa 1674 1022}%
\special{pa 1674 1006}%
\special{sp}%
%
\special{pn 8}%
\special{pa 1830 1064}%
\special{pa 1852 1042}%
\special{pa 1884 1036}%
\special{pa 1914 1042}%
\special{pa 1944 1052}%
\special{pa 1970 1070}%
\special{pa 1988 1096}%
\special{pa 1988 1112}%
\special{sp}%
%
\special{pn 8}%
\special{pa 2144 1140}%
\special{pa 2122 1162}%
\special{pa 2090 1168}%
\special{pa 2060 1164}%
\special{pa 2030 1154}%
\special{pa 2004 1136}%
\special{pa 1988 1110}%
\special{pa 1988 1092}%
\special{sp}%
%
\special{pn 8}%
\special{pa 2142 1142}%
\special{pa 2168 1124}%
\special{pa 2200 1122}%
\special{pa 2230 1128}%
\special{pa 2262 1136}%
\special{pa 2288 1154}%
\special{pa 2294 1174}%
\special{sp}%
%
\special{pn 8}%
\special{pa 2420 1236}%
\special{pa 2390 1226}%
\special{pa 2360 1216}%
\special{pa 2332 1200}%
\special{pa 2306 1184}%
\special{pa 2296 1162}%
\special{sp}%
%
\special{pn 8}%
\special{ar 2746 1142 422 144  1.7880735 2.4873061}%
%
\special{pn 8}%
\special{ar 2754 1116 250 188  0.8415189 2.0574284}%
%
\special{pn 8}%
\special{ar 2026 1028 104 466  1.5707963 3.0786622}%
%
\special{pn 8}%
\special{ar 2026 838 166 654  6.2073556 6.3536971}%
\special{ar 2026 838 166 654  6.4415020 6.5878434}%
\special{ar 2026 838 166 654  6.6756483 6.8219898}%
\special{ar 2026 838 166 654  6.9097947 7.0561361}%
\special{ar 2026 838 166 654  7.1439410 7.2902825}%
\special{ar 2026 838 166 654  7.3780873 7.5244288}%
\special{ar 2026 838 166 654  7.6122337 7.7585751}%
\special{ar 2026 838 166 654  7.8463800 7.8539816}%
%
\special{pn 8}%
\special{ar 2508 1140 144 542  1.5707963 2.9734285}%
%
\special{pn 8}%
\special{ar 2492 786 204 898  0.2222553 0.3311482}%
\special{ar 2492 786 204 898  0.3964840 0.5053769}%
\special{ar 2492 786 204 898  0.5707127 0.6796056}%
\special{ar 2492 786 204 898  0.7449413 0.8538342}%
\special{ar 2492 786 204 898  0.9191700 1.0280629}%
\special{ar 2492 786 204 898  1.0933987 1.2022916}%
\special{ar 2492 786 204 898  1.2676274 1.3765203}%
\special{ar 2492 786 204 898  1.4418560 1.5507490}%
%
\special{pn 8}%
\special{pa 1764 710}%
\special{pa 1786 690}%
\special{pa 1816 684}%
\special{pa 1848 688}%
\special{pa 1878 698}%
\special{pa 1904 716}%
\special{pa 1920 744}%
\special{pa 1920 760}%
\special{sp}%
%
\special{pn 8}%
\special{pa 2084 796}%
\special{pa 2062 816}%
\special{pa 2032 822}%
\special{pa 2000 818}%
\special{pa 1970 808}%
\special{pa 1944 788}%
\special{pa 1928 764}%
\special{pa 1928 746}%
\special{sp}%
%
\special{pn 8}%
\special{pa 2084 806}%
\special{pa 2106 784}%
\special{pa 2136 778}%
\special{pa 2168 782}%
\special{pa 2198 794}%
\special{pa 2224 812}%
\special{pa 2242 838}%
\special{pa 2240 854}%
\special{sp}%
%
\special{pn 8}%
\special{pa 2396 882}%
\special{pa 2376 904}%
\special{pa 2344 910}%
\special{pa 2312 906}%
\special{pa 2282 896}%
\special{pa 2258 878}%
\special{pa 2240 850}%
\special{pa 2240 834}%
\special{sp}%
%
\special{pn 8}%
\special{pa 2396 886}%
\special{pa 2420 868}%
\special{pa 2450 864}%
\special{pa 2482 870}%
\special{pa 2512 880}%
\special{pa 2538 896}%
\special{pa 2548 924}%
\special{pa 2548 924}%
\special{sp}%
%
\special{pn 8}%
\special{pa 2672 978}%
\special{pa 2642 968}%
\special{pa 2614 956}%
\special{pa 2586 942}%
\special{pa 2560 924}%
\special{pa 2550 904}%
\special{sp}%
%
\special{pn 8}%
\special{ar 3000 882 422 144  1.7886281 2.4860847}%
%
\special{pn 8}%
\special{ar 3008 858 250 188  0.8415189 2.0551544}%
%
\special{pn 8}%
\special{pa 2286 2636}%
\special{pa 2698 2180}%
\special{dt 0.045}%
\put(30.8700,-28.8600){\makebox(0,0)[lt]{${\Bbb R}$}}%
\end{picture}%
\hspace{2.5truecm}}

\vspace{1truecm}

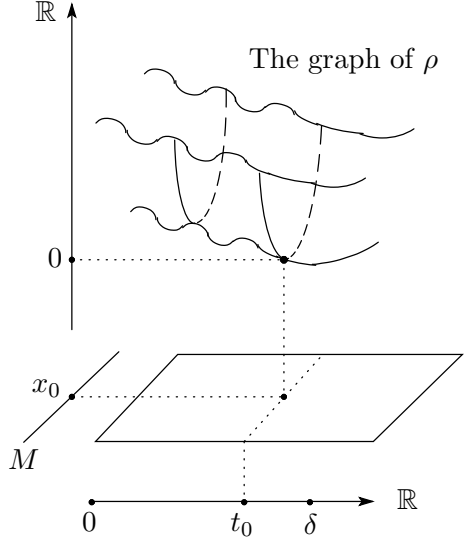
\centerline{{\bf Figure 4.}}

\vspace{0.5truecm}

\noindent
{\it Remark 5.1.} (i) 
According to the proof of the maximum principle by R.S. Hamilton (Theorem 9.1 of [Ha]), we can improve 
the statement of his maximum principle as follows.  

\vspace{0.15truecm}

{\sl Let $g_t\,\,\,(0\leq t<T)$ be a $C^{\infty}$-family of Riemannian metrics on a compact manifold 

$M$ and $S_t\,\,\,(0\leq t<T)$ be a $C^{\infty}$-family of symmetric $(0,2)$-tensor field on $M$.  

Assume that $S_t$'s ($0\leq t<T$) satisfy the following evolution equation:
$$\frac{\partial S}{\partial t}=\triangle S+\nabla_{\bar X_0}S+P(S),$$

where $\triangle S$ is the Laplacian of $S$ with respect to the connection of 
$\pi^{\ast}TM$ defined 

by the Levi-Civita connections $\nabla^t$'s of $g_t$, 
$X_0\in\Gamma(TM)$ and $P$ is a map of poly-

nomial type from $\Gamma(\pi_M^{\ast}(T^{(0,2)}M))$ to oneself.  
Assume that $P$ satisfies the follo-

wing condition:
$$\begin{array}{r}
\displaystyle{(\ast_S)\qquad\,\,X\in{\rm Ker}(S+\varepsilon g)_{(x,t)}\,\,\,
\Rightarrow\,\,\,P(S+\varepsilon g)_{(x,t)}(X,X)\geq0\qquad\qquad}\\
\displaystyle{(\forall\,\varepsilon>0,\,\,\,\,\forall\,(x,t)\in M\times[0,T)).\qquad\qquad}
\end{array}$$

If $S_0\geq0$ (resp. $>0$), then $S_t\geq0$ (resp. $>0$) holds 
for all $t\in[0,T)$.}

\vspace{0.3truecm}

\noindent
The null-eigenvector condition in [Ha] means the following condition:
$$\begin{array}{c}
\displaystyle{X\in{\rm Ker}(\widehat S)_{(x,t)}\,\,\,
\Rightarrow\,\,\,P(\widehat S)_{(x,t)}(X,X)\geq0}\\
\displaystyle{(\forall\,\widehat S\,:\,{\rm symmetric}\,\,(0,2){\rm -tensor}\,\,{\rm field}\,\,{\rm on}\,\,M,
\,\,\,\,\forall\,(x,t)\in M\times[0,T)).}
\end{array}$$
This condition is stronger than the above condition $(\ast_S)$.  
In [Hu1], G. Huisken proved the statement of Theroem 4.3 in [Hu1] by showing that the family 
$S=(S_{ij}:=\frac{h_{ij}}{H}-\varepsilon g_{ij})$ of symmetric $(0,2)$-tensor fields satisfies the above 
condition $(\ast_S)$ and applying the maximum principle of R.S. Hamilton.  In his proof, it is not shown that 
the family $S$ satisfies the null-eigenvector condition.  
The statement of Theorem 4.2 in [Hu2] also was proved by showing that some another family $S$ of symmetric 
$(0,2)$-tensor fields satisfies the above condition $(\ast_S)$.  

(ii) The constant $C_{\delta}$ in this proof corresponds to 
the constant $C$ in the proof of Theorem 9.1 in [Ha].  

\vspace{0.5truecm}

Similarly we obtain the following maximal principle for a $C^{\infty}$-family 
of $G$-invariant functions on $M$.  

\vspace{0.5truecm}

\noindent
{\bf Theorem 5.2.} {\sl Let $\rho$ be a $C^{\infty}$-function over 
$M\times[0,T)$ such that, for each $t\in [0,T)$, $\rho_t(:=\rho(\cdot,t))$ 
is a $G$-invariant function on $M$.  Assume that $\rho_t$'s ($0\leq t<T$) 
satisfy the following evolution equation:
$$\frac{\partial\rho}{\partial t}=\triangle_{\cal H}\rho
+d\rho(\bar X_0)+P(\rho),$$
where $X_0\in\Gamma(TM)$ and $P$ is a map of polynomial type from 
$C^{\infty}(M\times[0,T))$ to oneself.  

${\rm(i)}$ Assume that $P$ satisfies the following condition:
$$\begin{array}{r}
\displaystyle{(\rho+\varepsilon)_{(x,t)}=0\,\,\Rightarrow\,\,
P(\rho+\varepsilon)_{(x,t)}\geq0}\\
\displaystyle{(\forall\,\varepsilon>0,\,\,\,\,(x,t)\in M\times[0,T)).}
\end{array}$$
If $\rho_0\geq0$ (resp. $>0$), then $\rho_t\geq0$ 
(resp. $>0$) holds for all $t\in[0,T)$.

${\rm(ii)}$ Assume that $P$ satisfies the following condition:
$$\begin{array}{r}
\displaystyle{(\rho+\varepsilon)_{(x,t)}=0\,\,\Rightarrow\,\,
P(\rho+\varepsilon)_{(x,t)}\leq0}\\
\displaystyle{(\forall\,\varepsilon>0,\,\,\,\,(x,t)\in M\times[0,T)).}
\end{array}$$
If $\rho_0\leq0$ (resp. $<0$), then 
$\rho_t\leq0$ (resp. $<0$) holds for all $t\in[0,T)$.}

\section{Horizontally strongly convexity preservability theorem} 
Let $G\curvearrowright V$ be an isometric almost free action with minimal regularizable 
orbit of a Hilbert Lie group $G$ on a Hilbert space $V$ equipped with an inner 
product $\langle\,\,,\,\,\rangle$ and $\phi:V\to V/G$ the orbit map.  
Denote by $\widetilde{\nabla}$ the Riemannian connection of $V$.  
Set $n:={\rm dim}\,V/G-1$.  
Let $M(\subset V)$ be a $G$-invariant hypersurface in $V$ such that $\phi(M)$ is compact.  
Let $f$ be an inclusion map of $M$ into $V$ and 
$f_t\,\,(0\leq t<T)$ the regularized mean curvature flow starting from $f$.  
We use the notations in Section 4.  In the sequel, we omit the notation $f_{t\ast}$ for simplicity.  
For each $u\in V$, we set 
$$L:=\mathop{\sup}_{u\in V}\mathop{\max}_{(X_1,\cdots,X_5)\in({\widetilde{\cal H}}_1)_u^5}
\vert\langle{\cal A}^{\phi}_{X_1}
((\widetilde{\nabla}_{X_2}{\cal A}^{\phi})_{X_3}X_4),\,X_5\rangle\vert,$$
where $\widetilde{\cal H}_1:=\{X\in\widetilde{\cal H}\,\vert\,\,\,
\vert\vert X\vert\vert=1\}$.  
Assume that $L<\infty$.  Note that $L<\infty$ in the case where $V/G$ is compact.  
In this section, we prove the following horizontally strongly 
convexity preservability theorem by using results stated in Section 4 and Theorem 5.1.  

\vspace{0.5truecm}

\noindent
{\bf Theorem 6.1.} {\sl 
If $M$ satisfies $\vert\vert H_0\vert\vert^2(h_{\cal H})_{(\cdot,0)}
>2n^2L(g_{\cal H})_{(\cdot,0)}$, then $T<\infty$ holds and 
$\vert\vert H_t\vert\vert^2(h_{\cal H})_{(\cdot,t)}>2n^2L
(g_{\cal H})_{(\cdot,t)}$ holds for all $t\in[0,T)$.}

\vspace{0.5truecm}

\noindent
{\it Proof.} 
Since ${\cal A}^{\phi}_{\xi}$ is skew-symmetric, we have 
$${\rm Tr}(({\cal A}^{\phi}_{\xi})^2)_{\cal H}\leq0.\leqno{(6.1)}$$
From Corollary 4.11, 
${\rm Tr}(A_{\cal H})^2\geq\frac{\vert\vert H\vert\vert^2}{n}$ and $(6.1)$, 
we have 
$$\frac{\partial\vert\vert H\vert\vert}{\partial t}\geq
\triangle_{\cal H}\vert\vert H\vert\vert+\frac{\vert\vert H\vert\vert^3}{n}.
\leqno{(6.2)}$$
Define a function $\rho$ over $[0,T)$ by $\rho(t):=\min\,\vert\vert H_t
\vert\vert$.  Form $(6.2)$, we have $\frac{d\rho}{dt}\geq\frac1n\rho^3$.  
Also we have $\rho(0)>0$ by the assumption.  
Hence we obtain $T\leq\frac{n}{2\rho(0)^2}$.  
Set 
$\displaystyle{S:=\frac{1}{\vert\vert H\vert\vert}h
-\frac{2n^2L}{\vert\vert H\vert\vert^3}g}$ and 
$S_{\varepsilon}:=S+\varepsilon g$, where $\varepsilon$ is a positive constant.  
Take $X,Y\in{\cal H}$.  
By using Lemma 4.2, Theorem 4.6, Corollary 4.8 and Lemma 4.10, we can show 
$$\begin{array}{l}
\hspace{0.5truecm}\displaystyle{
\frac{\partial(S_{\varepsilon})_{\cal H}}{\partial t}(X,Y)}\\
\displaystyle{=\frac{1}{\vert\vert H\vert\vert}
(\triangle_{\cal H}^{\cal H}h_{\cal H})(X,Y)
-2((A_{\cal H})^2)_{\sharp}(X,Y)
-2(({\cal A}^{\phi}_{\xi})^2)_{\sharp}(X,Y)}\\
\hspace{0.5truecm}\displaystyle{-\frac{1}{\vert\vert H\vert\vert^2}\left(
\triangle_{\cal H}^{\cal H}\vert\vert H\vert\vert
-2\vert\vert H\vert\vert{\rm Tr}(({\cal A}^{\phi}_{\xi})^2)_{\cal H}
-4n^2L\right)h_{\cal H}(X,Y)}\\
\hspace{0.5truecm}\displaystyle{
-\frac{1}{\vert\vert H\vert\vert}{\cal R}(X,Y)-2\varepsilon\vert\vert H\vert\vert h_{\cal H}(X,Y)}\\
\hspace{0.5truecm}\displaystyle{
+\frac{3n^2L}{\vert\vert H\vert\vert^4}\left(
\triangle_{\cal H}^{\cal H}\vert\vert H\vert\vert
+\vert\vert H\vert\vert{\rm Tr}(A_{\cal H})^2
-3\vert\vert H\vert\vert{\rm Tr}(({\cal A}^{\phi}_{\xi})^2)_{\cal H}\right)
g_{\cal H}(X,Y)}
\end{array}\leqno{(6.3)}$$
Also, we have 
$$\begin{array}{l}
\displaystyle{(\nabla^{\cal H}_{{\rm grad}\vert\vert H\vert\vert}(S_{\varepsilon})_{\cal H})(X,Y)
=\frac{1}{\vert\vert H\vert\vert}
(\nabla^{\cal H}_{{\rm grad}\vert\vert H\vert\vert}h_{\cal H})(X,Y)}\\
\hspace{4.65truecm}\displaystyle{
-\frac{\vert\vert{\rm grad}\vert\vert H\vert\vert\,\vert\vert^2}
{\vert\vert H\vert\vert^2}h(X,Y)}\\
\hspace{4.65truecm}\displaystyle{+\frac{3n^2L}{\vert\vert H\vert\vert^4}
\vert\vert{\rm grad}\vert\vert H\vert\vert\,\vert\vert^2g(X,Y)}
\end{array}\leqno{(6.4)}$$
and 
$$\begin{array}{l}
\displaystyle{(\triangle_{\cal H}^{\cal H}(S_{\varepsilon})_{\cal H})(X,Y)=
\frac{1}{\vert\vert H\vert\vert}(\triangle_{\cal H}^{\cal H}h_{\cal H})(X,Y)
-\frac{2}{\vert\vert H\vert\vert^2}
(\nabla^{\cal H}_{{\rm grad}\vert\vert H\vert\vert}h_{\cal H})(X,Y)}\\
\hspace{3.3truecm}\displaystyle{+\frac{1}{\vert\vert H\vert\vert^3}
\left(2\vert\vert{\rm grad}\vert\vert H\vert\vert\,\vert\vert^2
-\vert\vert H\vert\vert\triangle_{\cal H}^{\cal H}\vert\vert H\vert\vert\right)
h_{\cal H}(X,Y)}\\
\hspace{3.3truecm}\displaystyle{+\frac{3n^2L}{\vert\vert H\vert\vert^5}\left(
-4\vert\vert{\rm grad}\vert\vert H\vert\vert\,\vert\vert^2
+\vert\vert H\vert\vert\triangle_{\cal H}^{\cal H}\vert\vert H\vert\vert\right)
g_{\cal H}(X,Y).}
\end{array}\leqno{(6.5)}$$
From $(6.3),\,(6.4)$ and $(6.5)$, we have 
$$\begin{array}{l}
\displaystyle{\frac{\partial(S_{\varepsilon})_{\cal H}}{\partial t}(X,Y)
=\triangle_{\cal H}^{\cal H}(S_{\varepsilon})_{\cal H}(X,Y)
+\frac{2}{\vert\vert H\vert\vert}
(\nabla^{\cal H}_{{\rm grad}\vert\vert H\vert\vert}
(S_{\varepsilon})_{\cal H})(X,Y)}\\
\hspace{3.25truecm}\displaystyle{+P(S_{\varepsilon})(X,Y),}
\end{array}\leqno{(6.6)}$$
where $P(S_{\varepsilon})$ is defined by 
$$\begin{array}{l}
\displaystyle{P(S_{\varepsilon})(Z,W)
:=-2((A_{\cal H})^2)_{\sharp}(Z,W)
-2(({\cal A}^{\phi}_{\xi})^2)_{\sharp}(Z,W)
-\frac{1}{\vert\vert H\vert\vert}{\cal R}(Z,W)}\\
\hspace{2.2truecm}\displaystyle{
+\frac{1}{\vert\vert H\vert\vert^2}\left(2\vert\vert H\vert\vert{\rm Tr}
(({\cal A}^{\phi}_{\xi})^2)_{\cal H}+4n^2L
\right)h_{\cal H}(Z,W)}\\
\hspace{2.2truecm}\displaystyle{
+\frac{6n^2L}{\vert\vert H\vert\vert^3}
\left({\rm Tr}(A_{\cal H})^2-3{\rm Tr}(({\cal A}^{\phi}_{\xi})^2)_{\cal H}
+\frac{2\vert\vert{\rm grad}\,\vert\vert H\vert\vert\,\vert\vert^2}
{\vert\vert H\vert\vert^2}\right)g_{\cal H}(Z,W)}\\
\hspace{2.2truecm}\displaystyle{
-2\varepsilon\vert\vert H\vert\vert h_{\cal H}(Z,W)}
\end{array}$$
for $Z,W\in\pi_M^{\ast}TM$.  
Fix any positive constant $\varepsilon_0$ and any $(x_0,t_0)\in M\times[0,T)$.  
Assume that ${\rm Ker}((S_{\varepsilon_0})_{\cal H})_{(x_0,t_0)}\not=\{0\}$.  
Take $X_0\in{\rm Ker}
((S_{\varepsilon_0})_{\cal H})_{(x_0,t_0)})$ with $g(X_0,X_0)=1$.  
Since 
$$h(X_0,Y)=\left(\frac{2n^2L}{\vert\vert H\vert\vert^2}-\varepsilon_0\vert\vert H\vert\vert\right)g(X_0,Y)
\,\,\,\,\,(\forall\,Y\in{\cal H}),$$
we have 
$$A_{\cal H}X_0=\left(\frac{2n^2L}{\vert\vert H\vert\vert^2}
-\varepsilon_0\vert\vert H\vert\vert\right)X_0.$$
For simplicity, we set 
$\lambda_1:=\frac{2n^2L}{\vert\vert H\vert\vert^2}-\varepsilon_0\vert\vert H\vert\vert$.  
By using the first relation in Lemma 4.10, we have 
$$\begin{array}{l}
\displaystyle{P(S_{\varepsilon_0})(X_0,X_0)
=\frac{6n^2L}{\vert\vert H\vert\vert^3}{\rm Tr}(A_{\cal H})^2
+\frac{12n^2L}{\vert\vert H\vert\vert^5}
\vert\vert{\rm grad}\vert\vert H\vert\vert\,\vert\vert^2}\\
\hspace{1.5truecm}\displaystyle{
-2(({\cal A}^{\phi}_{\xi})^2)_{\sharp}(X_0,X_0)
-\left(\frac{14n^2L}{\vert\vert H\vert\vert^3}
+2\varepsilon_0\right){\rm Tr}(({\cal A}^{\phi}_{\xi})^2)_{\cal H}
}\\
\hspace{1.5truecm}\displaystyle{
+\frac{4}{\vert\vert H\vert\vert}{\rm Tr}^{\bullet}_{g_{\cal H}}
\langle{\cal A}^{\phi}_{X_0}\bullet,{\cal A}^{\phi}_{X_0}(A_{\cal H}\bullet)\rangle
-\frac{4}{\vert\vert H\vert\vert}{\rm Tr}^{\bullet}_{g_{\cal H}}
\langle{\cal A}^{\phi}_{\bullet}X_0,{\cal A}^{\phi}_{\bullet}(A_{\cal H}X_0)\rangle}\\
\hspace{1.5truecm}\displaystyle{
+\frac{3}{\vert\vert H\vert\vert}{\rm Tr}^{\bullet}_{g_{\cal H}}
\langle(\widetilde{\nabla}_{\bullet}{\cal A}^{\phi})_{X_0}\xi,
{\cal A}^{\phi}_{\bullet}X_0\rangle
-\frac{1}{\vert\vert H\vert\vert}{\rm Tr}^{\bullet}_{g_{\cal H}}
\langle(\widetilde{\nabla}_{X_0}{\cal A}^{\phi})_{\bullet}\xi,
{\cal A}^{\phi}_{X_0}\bullet\rangle}\\
\hspace{1.5truecm}\displaystyle{
+\frac{2}{\vert\vert H\vert\vert}{\rm Tr}^{\bullet}_{g_{\cal H}}
\langle(\widetilde{\nabla}_{\bullet}{\cal A}^{\phi})_{\bullet}X_0,
{\cal A}^{\phi}_{X_0}\xi\rangle.}
\end{array}\leqno{(6.7)}$$
Hence, since ${\rm Tr}(A_{\cal H})^2\geq\frac{\vert\vert H\vert\vert^2}{n}$, 
$(({\cal A}^{\phi}_{\xi})^2)_{\sharp}(X,X)\leq0$, 
${\rm Tr}(({\cal A}^{\phi}_{\xi})^2)_{\cal H}\leq0$ and the definition of $L$, 
we have 
$$\begin{array}{l}
\displaystyle{P(S_{\varepsilon_0})(X_0,X_0)
>\frac{4}{\vert\vert H\vert\vert}{\rm Tr}^{\bullet}_{g_{\cal H}}
\langle{\cal A}^{\phi}_{X_0}\bullet,{\cal A}^{\phi}_{X_0}(A_{\cal H}\bullet)\rangle}\\
\hspace{3.3truecm}\displaystyle{-\frac{4}{\vert\vert H\vert\vert}
{\rm Tr}^{\bullet}_{g_{\cal H}}\langle{\cal A}^{\phi}_{\bullet}X_0,
{\cal A}^{\phi}_{\bullet}(A_{\cal H}X_0)\rangle.}
\end{array}\leqno{(6.8)}$$
Since $A_{\cal H}X_0=\lambda_1X_0$, 
$X\in{\rm Ker}((S_{\varepsilon_0})_{\cal H})_{(x_0,t_0)}$ and 
$((S_{\varepsilon_0})_{\cal H})_{(x_0,t_0)}\geq0$, we may assume that $\lambda_1$ is 
the smallest eigenvalue of $(A_{\cal H})_{(x_0,t_0)}$.  Let 
$\{\lambda_i\,\vert\,i=1,\cdots,n\}\,\,(\lambda_1\leq\cdots
\leq\lambda_n)$ be the set of all eigenvalues of $(A_{\cal H})_{(x_0,t_0)}$.  
Let $\{e_1,\cdots,e_n\}$ be an orthonormal base of $T_{x_0}M$ with respect to 
$(g_{\cal H})_{(x_0,t_0)}$ satisfying 
$e_1=X_0$ and $A_{\cal H}e_i=\lambda_ie_i\,\,\,(i=2,\cdots,n)$.  
Then we have 
$$\begin{array}{l}
\hspace{0.4truecm}\displaystyle{
\frac{4}{\vert\vert H\vert\vert}{\rm Tr}^{\bullet}_{g_{\cal H}}
\langle{\cal A}^{\phi}_{X_0}\bullet,{\cal A}^{\phi}_{X_0}(A_{\cal H}\bullet)\rangle
-\frac{4}{\vert\vert H\vert\vert}{\rm Tr}^{\bullet}_{g_{\cal H}}
\langle{\cal A}^{\phi}_{\bullet}X_0,{\cal A}^{\phi}_{\bullet}(A_{\cal H}X_0)
\rangle}\\
\displaystyle{=\frac{4}{\vert\vert H\vert\vert}\sum_{i=1}^n
(\lambda_i-\lambda_1)\langle{\cal A}^{\phi}_{X_0}e_i,{\cal A}^{\phi}_{X_0}e_i\rangle
\geq0.}
\end{array}$$
From $(6.8)$ and this inequality, we obtain 
$P(S_{\varepsilon_0})(X_0,X_0)\geq0$.  
Hence it follows from the arbitrarinesses of $\varepsilon_0$ and $(x_0,t_0)$ that $P$ satisfies the condition 
$(\ast_{S_{\cal H}}^+)$.  Therefore it follows from Theorem 5.1 that 
$(S_{\cal H})_{(\cdot,t)}>0$ holds for all $t\in[0,T)$.  
\hspace{1.5truecm}q.e.d.

\section{Strongly convex preservability theorem in the orbit space} 
Let $V,\,G$ and $\phi$ be as in the previous section.  
Set $N:=V/G$ and $n:={\rm dim}\,V/G-1$.  Denote by $g_N$ and $R_N$ 
the Riemannian orbimetric and the curvature orbitensor of $N$.  
Also, $\nabla^N$ the Riemannian connection of $g_N\vert_{N\setminus{\rm Sing}(N)}$.  
Denote by $\vert\vert\nabla^NR_N\vert\vert$ the norm of $\nabla^NR_N$ 
with respect to $g_N$.  
Set $L_N:=\sup\,\vert\vert\nabla^NR_N\vert\vert$.  
Assume that $L_N\,<\,\infty$.  
Let $\overline M$ be a compact suborbifold of codimension one in $N$ immersed by 
$\overline f$ and $\overline f_t$ ($t\in[0,T)$) the mean curvature flow starting from 
$\overline f$.  Denote by $\overline g_t,\overline h_t,\overline A_t$ and $\overline H_t$ 
be the induced orbimetric, the second fundamental orbiform, the shape orbitensor and 
the mean curvature orbifunction of $\overline f_t$, respectively, and 
$\overline{\xi}_t$ the unit normal vector field of 
$\overline f_t\vert_{\overline M\setminus{\rm Sing}(\overline M)}$.  

From Theorem 6.1, we obtain the following strongly convexity preservability theorem for 
compact suborbifolds in $N$.  

\vspace{0.5truecm}

\noindent
{\bf Theorem 7.1.} {\sl 
If $\overline f$ satisfies $\vert\vert\overline H_0\vert\vert^2\overline h_0
>n^2L_N\overline g_0$, then $T<\infty$ holds and 
$\vert\vert\overline H_t\vert\vert^2\overline h_t>n^2L_N\overline g_t$ holds for all $t\in[0,T)$.}

\vspace{0.5truecm}

\noindent
{\it Proof.} Set $M:=\{(x,u)\in\overline M\times V\,\vert\,\overline f(x)=\phi(u)\}$ and 
define $f:M\to V$ by $f(x,u)=u$ ($(x,u)\in M$).  
It is clear that $f$ is an immersion.  
Denote by $H_0$ the regularized mean curvature vector of $f$.  
Define a curve $c_x:[0,T)\to N$ by $c_x(t):=\overline f_t(x)$ ($t\in[0,T)$) and 
let $(c_x)_u^L$ be the horizontal lift of $c_x$ with $(c_x)_u^L(0)=u$ and 
$((c_x)_u^L)'(0)=(H_0)_{(x,u)}$, where $(x,u)\in M$.  
Define an immersion $f_t:M\hookrightarrow V$ by 
$f_t(x,u):=(c_x)_u^L(t)$ ($(x,u)\in M$).  
Then $f_t$ ($t\in[0,T)$) is the regularized mean curvature flow starting from $f$ 
(see the proof of Theorem 4.1).  
Denote by $g_t,\,h_t,\,A^t$ and $H_t$ the induced metric, the second fundamental form, 
the shape tensor and the mean curvature vector of $f_t$, respectively.  
By the assumption, $\overline f_0$ satisfies $\vert\vert\overline H_0\vert\vert^2\overline h_0
>n^2L_N\overline g_0$.  
Also, we can show $L_N=2L$ by long calculation, where $L$ is as in the prevoius section.  
From these facts, we can show that $f_0$ satisfies 
$\vert\vert H_0\vert\vert^2(h_{\cal H})_0>2n^2L(g_{\cal H})_0$.  
Hence, it follows from Theorem 6.1 that $f_t$ ($t\in[0,T)$) satisfies 
$\vert\vert H_t\vert\vert^2(h_{\cal H})_t>2n^2L(g_{\cal H})_t$.  
Furthermore, it follows from this fact that $\overline f_t$ ($t\in[0,T)$) satisfies 
$\vert\vert\overline H_t\vert\vert^2\overline h_t>n^2L_N\overline g_t$.  
\hspace{6.9truecm}q.e.d.

\vspace{0.5truecm}

\noindent
{\it Remark 7.1.} In the case where the $G$-action is free and hence $N$ is a (complete) 
Riemannian manifold, Theorem 7.1 implies the strongly convexity preservability theorem by G. Huisken 
(see [Hu2, Theorem 4.2]).  

\vspace{1truecm}


\centerline{{\bf References}}

\vspace{0.5truecm}

{\small 
\noindent
[AK] A. Adem and M. Klaus, 
Lectures on orbifolds and group cohomology, Advanced 

Lectures in Mathematics {\bf 16}, Transformation groups 
and Moduli Spaces of Curves 

(2010), 1-17.

\noindent
[BB] J. E. Borzelliono and V. Brunsden, 
Orbifold homeomorphism and Diffeomorphism 

groups, Report in International Conference 
on Infinite Dimensional Lie Groups in 

Geometry and Representation Theory at 
Howard University, Washington DC, 2000.

\noindent
[Ge] C. Gerhart, Curvature Problems, 
Series in Geometry and Topology, {\bf 39}, Internatio-

nal Press 2006.

\noindent
[GKP] P. Gilkey, H. J. Kim and J. H. Park, 
Eigenforms of the Laplacian for Riemannian 

$V$-submersions, Tohoku Math. J. 
{\bf 57} (2005) 505-519.

\noindent
[Ha] R. S. Hamilton, Three-manifolds with 
positive Ricci curvature, J. Differential Geom. 

{\bf 17} (1982) 255-306.

\noindent
[Hu1] G. Huisken, Flow by mean curvature of 
convex surfaces into spheres, J. Differential 

Geom. {\bf 20} (1984) 237-266.

\noindent
[Hu2] G. Huisken, Contracting convex hypersurfaces 
in Riemannian manifolds by their mean 

curvature, Invent. math. {\bf 84} (1986) 463-480.

\noindent
[HLO] E. Heintze, X. Liu and C. Olmos, Isoparametric submanifolds and a 
Chevalley-type 

restriction theorem, Integrable systems, geometry, and topology, 151-190, 
AMS/IP Stud. 

Adv. Math. 36, Amer. Math. Soc., Providence, RI, 2006.

%

\noindent
[HP] G. Huisken and A. Polden, Geometric evolution equations for 
hypersurfaces, Calculus 

of Variations and Geometric Evolution Problems, CIME Lectures 
at Cetraro of 1966 

(S. Hildebrandt and M. Struwe, eds.) Springer, 1999.

\noindent
[KT] C. King and C. L. Terng, Minimal submanifolds in path spaces, 
Global Analysis and 

Modern Mathematics, edited by K. Uhlenbeck, Publish or Perish (1993) 253-282.

\noindent
[Koi1] N. Koike, On proper Fredholm submanifolds in a Hilbert space arising 
from subma-

nifolds in a symmetric space, Japan. J. Math. {\bf 28} (2002) 61-80.

\noindent
[Koi2] N. Koike, Collapse of the mean curvature flow for equifocal 
submanifolds, Asian J. 

Math. {\bf 15} (2011) 101-128.

\noindent
[O'N] B. O'Neill, The fundamental equations of a submersion, Michigan Math. J. 
{\bf 13} (1966) 

459-469.

\noindent
[P] R. S. Palais, 
Morse theory on Hilbert manifolds, Topology {\bf 2} (1963) 299-340.

\noindent
[PT] R. S. Palais and C. L. Terng, Critical point theory and submanifold 
geometry, Lecture 

Notes in Math. {\bf 1353}, Springer, Berlin, 1988.

\noindent
[Sa] I. Satake, 
The Gauss-Bonnet theorem fro $V$-manifolds, 
J. Math. Soc. Japan {\bf 9} (1957) 

464-492.

\noindent
[Sh] T. Shioya, 
Eigenvalues and suspension structure of compact 
Riemannian orbifolds with 

positive Ricci curvature, Manuscripta Math. {\bf 99} (1999) 509-516.

\noindent
[Th] W. P. Thurston, 
Three-dimensional geometry and topology, Princeton 
Mathematical 

Series 35, Princeton University Press, Princeon, N. J., 1997.

\noindent
[Te] C. L. Terng, Proper Fredholm submanifolds of Hilbert space, J. Differential Geometry 

{\bf 29} (1989) 9-47.

\noindent
[TT] C. L. Terng and G. Thorbergsson, Submanifold geometry in symmetric spaces, J. Diff-

erential Geometry {\bf 42} (1995) 665-718.

\noindent
[Z] X. P. Zhu, Lectures on mean curvature flows, Studies in Advanced Math., 
AMS/IP, 2002.

\vspace{0.5truecm}

{\small 
\rightline{Department of Mathematics, Faculty of Science}
\rightline{Tokyo University of Science, 1-3 Kagurazaka}
\rightline{Shinjuku-ku, Tokyo 162-8601 Japan}
\rightline{(koike@ma.kagu.tus.ac.jp)}
}

\end{document}